\documentclass[10pt,a4paper]{article} 

\usepackage[left=2cm,right=2cm, top=2cm, bottom=2cm]{geometry}

\usepackage{tikz}
\usepackage{pgfplots}
\usepackage{hyperref}
\usepackage{mathtools}
\usepackage{amssymb}
\usepackage{amsmath}
\usepackage{amsthm}
\usepackage{cleveref}
\usepackage{booktabs}
\usepackage{float}
\usepackage{threeparttable}
\usepackage{multirow}
\usepackage{subcaption}
\usepackage{booktabs}
\usepackage{array}
\usepackage{float}
\usepackage{bm}
\usepackage[dvipsnames]{xcolor}
\usepackage{adjustbox}
\usepackage{xspace}

\definecolor{EM_color}{RGB}{255, 150, 150}
\definecolor{EI_color}{RGB}{180, 0, 0}
\definecolor{SRA1_color}{RGB}{0, 0, 180}
\definecolor{Chang_color}{RGB}{150, 150, 255}

\newcolumntype{R}{>{\footnotesize}l}
\pgfplotsset{compat=newest}

\renewcommand{\d}{\,\mathrm d}

\newcommand{\R}{\mathbb{R}} 
\newcommand{\N}{\mathbb{N}} 
\renewcommand{\O}{\mathcal{O}}

\newcommand{\BackwardNormalstartingNetworkDiscretizedGeneralStrong}{\ensuremath{\bar y^{ \mathcal{N},\theta}}}
\newcommand{\BackwardNormalstartingNetworkDiscretizedItoTaylor}{\ensuremath{\hat y^{ \mathcal{N},\theta}}}

\newcommand{\BackwardNormalstartingNetworkDiscretized}{\ensuremath{\bar y^{\mathcal{N},\theta}}}

\newcommand{\BackwardNormalstartingNetwork}{\ensuremath{ y^{\mathcal{N},\theta}}}
\newcommand{\BackwardNormalstarting}{\ensuremath{ y^{\mathcal{N}}}}
\newcommand{\BackwardNetwork}{\ensuremath{ y^{\theta}}}

\newcommand{\Backward}{ y}

\newcommand{\Forward}{ x}
\newcommand{\ReferenceSDE}{ x}
\newcommand{\Discrete}{ \bar  x}

\newcommand{\SchemeClassName}{iterative-wise general strong scheme\xspace}
\newcommand{\SchemeClassShortHand}{IGSS\xspace}

\newtheorem{thm}{Theorem}[section]
\newtheorem{lem}[thm]{Lemma}

\newtheorem{cor}[thm]{Corollary}

\newtheorem{defi}[thm]{Definition}
\newtheorem{rem}[thm]{Remark}

\newtheorem{prob}{Problem}

\usepackage{enumitem}

\newtheorem{assumption}{Assumption}

\newlist{asslist}{enumerate}{1}
\setlist[asslist,1]{
  label=\theassumption.\Alph*.,
  ref=\theassumption.\Alph*,
  leftmargin=3.5em,
  itemsep=0.1\baselineskip
}

\crefname{assumption}{Assumption}{Assumptions}
\crefalias{asslisti}{assumption}

\crefname{prob}{Problem}{Problems}
\Crefname{prob}{Problem}{Problems}

\newcounter{probitem}[prob]           

\crefname{thm}{Theorem}{Theorems}
\Crefname{thm}{Theorem}{Theorems}

\crefname{lem}{Lemma}{Lemmas}
\Crefname{lem}{Lemma}{Lemmas}

\crefname{section}{Section}{Sections}
\Crefname{section}{Section}{Sections}

\crefname{subsection}{Subsection}{Subsections}
\Crefname{subsection}{Subsection}{Subsections}

\crefname{probitem}{Problem}{Problems}
\Crefname{probitem}{Problem}{Problems}

\newsavebox{\LatPlotBox}
\title{Stochastic Sampling for Generative Diffusion Models: From Euler-Maruyama to Higher-Order Schemes\thanks{This work was partially supported by the Alexander von Humboldt Foundation.}}

\author{Emanuel Pfarr\thanks{Institute of Mathematics, Julius-Maximilians-Universit\"at W\"urzburg, W\"urzburg, Germany (emanuel.pfarr@uni-wuerzburg.de)} \quad Radu Timofte\thanks{Center for Artificial Intelligence and Data Science, Julius-Maximilians-Universit\"at W\"urzburg, W\"urzburg, Germany (radu.timofte@uni-wuerzburg.de)} \quad Frank Werner\thanks{Institute of Mathematics, Julius-Maximilians-Universit\"at W\"urzburg, W\"urzburg, Germany (frank.werner@uni-wuerzburg.de).}}

\begin{document}
\maketitle
\allowdisplaybreaks
\begin{abstract}
We develop a convergence analysis for generative diffusion models that simultaneously accounts for the three principal error sources in stochastic sampling: initialization, score matching, and discretization of the reverse-time stochastic differential equation (SDE). Our central tool is the notion of a \emph{\SchemeClassName}, a broad class of discretization methods for the reverse dynamics defined via explicit, index-wise tolerances on their Itô--Taylor coefficients. This notion refines Kloeden and Platen's general strong schemes into an iterate-wise formulation whose per-step conditions imply the classical remainder bound. We prove a 2-Wasserstein convergence theorem that applies to this entire class, reducing the analysis of any concrete sampler to a finite checklist, and covers general forward processes with time-dependent, spatially linear drift and spatially independent diffusion coefficient. We focus on dissipative reverse-time SDEs, which allows us to relax the $T$-dependence of the discretization and score-matching errors from exponential to polynomial. We instantiate our theorem for Euler--Maruyama, the exponential integrator, and, as our main application, derivative-free stochastic Runge--Kutta schemes of strong order $1.5$. To our knowledge, these are the first full convergence guarantees for solvers whose strong order strictly exceeds that of Euler--Maruyama. We further derive the resulting iteration complexity and a parameter-selection rule for the terminal time, score accuracy, and step size. Numerical experiments on Gaussian toy models and CIFAR-10 are consistent with the predicted convergence orders, and on ImageNet, out-of-the-box higher-order SDE solvers can yield substantial speed-ups, requiring at least four times fewer score evaluations than Euler--Maruyama.
\end{abstract}

\section{Introduction}

Generative diffusion models (GDMs) have, within a few years, become one of the dominant paradigms for generative modeling of complex, high-dimensional data, underlying state-of-the-art methods for image \cite{song2020score,ho2020denoising,dhariwal2021diffusion,rombach2022high}, audio \cite{kong2021diffwave}, and video \cite{ho2022video} synthesis. Their success rests on a deceptively simple principle: rather than modeling the data distribution directly, one learns to reverse a stochastic process that gradually transforms it into a tractable reference distribution, so that generation reduces to the numerical simulation of a stochastic or ordinary differential equation. Improving the quality and speed of GDMs therefore proceeds on two fronts: the learning front, concerning how well the neural network approximates the score, and the numerical front, concerning the efficiency and theoretical guarantees of the schemes used to simulate the reverse-time dynamics. While both fronts have been explored extensively in practice, the theory of the numerical front is still taking shape \cite{chen2023improved,benton2023nearly,wu2024stochastic}, and this paper aims to advance it. Specifically, we develop a general convergence theory for discretizations of the reverse-time SDE that simultaneously accounts for all three sources of error in diffusion-based generation: score approximation, time discretization, and initialization at the reference distribution. This yields the first complete convergence proof for an entire class of stochastic samplers whose convergence order provably exceeds that of the Euler--Maruyama scheme.
Within this theory, we pay particular attention to the case of a dissipative drift, where tracking the time horizon $T$ explicitly throughout our proofs yields error bounds that do not deteriorate exponentially in $T$. 
Since long-time error control for higher-order schemes is of interest well beyond generative modeling, we believe these results are of independent interest for the numerical analysis of SDEs, complementing existing uniform-in-time results for Euler-type schemes \cite{fang2020adaptive}.

\subsection{Background}
Generative diffusion models \cite{song2020score} consist of two parts. The so-called forward process connects an unknown data distribution $\mathcal{P}_0$ on $\R^d$ to a standard Gaussian distribution through a forward SDE with additive noise
\begin{equation}
    \label{eq:RefernceSDE}
    \d\Forward_t = a\left(t,\Forward_t\right)\d t + b\left(t\right)\d W_t,
    \quad \Forward_0 \sim \mathcal{P}_0.
\end{equation}
Here, the drift $a\colon[0,\infty)\times \R^d\to \R^d$ and the diffusion coefficient $b\colon[0,\infty)\to \R^{d \times d}$ are carefully chosen functions to ensure $x_t\stackrel{D}{\longrightarrow}\mathcal{N}(0,\sigma^2I_d)$ as $t\to\infty$, with a known constant $\sigma^2\in \R^+$ and the $d$-dimensional identity matrix $I_d$\footnote{ A forward SDE with multiplicative noise, i.e. a spatially dependent diffusion coefficient $b(t,x)$ is possible in principle, at the price of an additional term $\nabla\!\cdot\!(bb^\top)$ in \eqref{eq:ReferenceReverseSDE} and $-\tfrac12\nabla\!\cdot\!(bb^\top)$ in \eqref{eq:ProbabilityFlowODE}, cf.\ \cite{anderson1982reverse,song2020score}. As a further consequence, the reverse-time SDE would also lose its additive noise property. Because all schedules used in practice, in particular all of Table \ref{tab:PopularChoices}, have a spatially independent diffusion coefficient, we restrict to this case throughout.}. Furthermore, the coefficient functions should be chosen in such a way that the conditional forward process $x_t|x_0$ can be solved in closed form. In Table \ref{tab:PopularChoices} we list popular choices for the coefficient functions in the GDM literature. 
\begin{table}[!htbp]
    \centering
    \begin{tabular}{|c|c|c|c|}
    \hline
         Name & $a(t,x)$ & $b(t)$ & References \\
         \hline
         Variance Preserving (VP) & $-\tfrac{1}{2}\beta(t)\,x$ & $\sqrt{\beta(t)}\, I_d$ &
         \cite{song2020score, ho2020denoising} \\
         \hline
         Variance Exploding (VE) & $0$ & $\sqrt{\tfrac{\d}{\d t}\sigma^2(t)}\, I_d$ &
         \cite{song2020score, song2019generative} \\
         \hline
         Ornstein--Uhlenbeck (OU) & $-x$ & $\sqrt{2}\, I_d$ &
         \cite{bruno2025OU_convergence_Log_Concave, gentilonisilveri2025OU_Convergence} \\
         \hline
    \end{tabular}
    \caption{Popular choices of the coefficient functions in
    \eqref{eq:RefernceSDE} in the GDM literature. Here, $\beta\colon[0,\infty)
    \to (0,\infty)$ is a noise schedule and $\sigma^2(t)$ an increasing
    variance function.}
    \label{tab:PopularChoices}
\end{table}
The second part of generative diffusion models, the so-called backward process, is used to actually generate samples from the data distribution $\mathcal{P}_0$. 
It consists of simulating the reverse-time SDE of the forward process \eqref{eq:RefernceSDE}. Using classical time-reversal results \cite{anderson1982reverse,haussmann1986time}, the reverse-time dynamics are given by the SDE with additive noise
\begin{equation}
\label{eq:ReferenceReverseSDE}
    \d\Backward_{t}
    = \bigl[-a\left(\tau,\Backward_t\right)
      + b\left(\tau\right)b\left(\tau\right)^{\!\top}s(\tau, \Backward_t)\bigr]\d t
      + b\left(\tau\right)\d W_t,
    \quad \Backward_0 \sim \mathcal{L}\left(\Forward_T\right), 
\end{equation}
with $\tau=T-t$, the law $\mathcal{L}(x_T)$ of the random variable $x_T$ and the so-called \emph{score function} $s\colon[0,\infty)\times \R^d \to \R^d$ given by 
\[s\left(t,x\right) := \nabla_{x} \log p_t\left(x\right).\] 
Here and below, $p_t(x)$ denotes the density of $\mathcal{L}(x_t)$. Additionally, as observed in \cite{song2020score}, the reverse-time SDE \eqref{eq:ReferenceReverseSDE} shares its marginal distributions $p_t$ with the so-called \emph{probability flow ODE}
\begin{equation}
\label{eq:ProbabilityFlowODE}
    \d\Backward_{t}
    = \Bigl[-a\left(\tau,\Backward_t\right)
      + \tfrac12 b\left(\tau\right)b\left(\tau\right)^{\!\top}s(\tau, \Backward_t)\Bigr]\d t,
    \quad \Backward_0 \sim \mathcal{L}\left(\Forward_T\right),
\end{equation}
which gives us a second, deterministic route to simulate the backward process. To simulate the reverse-time SDE \eqref{eq:ReferenceReverseSDE} and probability flow ODE \eqref{eq:ProbabilityFlowODE}, we replace the typically intractable score function $s$ by the so-called \emph{score network} $s_\theta\colon[0,\infty)\times \R^d \to \R^d$, a (deep) neural network which is trained to approximate the score function using the denoising score matching (DSM) objective
\begin{equation}
    \label{eq:DSM}
    \mathbb{E}_{t \sim \mathcal{U}[0,T]}
    \Bigl[\lambda\left(t\right)\,\mathbb{E}_{\Forward_0,\,\Forward_t}
    \|s_\theta\left(t,\Forward_t\right) - \nabla\log p_{t|0}\left(\Forward_t|\Forward_0\right)\|^2\Bigr],
\end{equation}
where $\nabla\log p_{t|0}\left(\Forward_t|\Forward_0\right)$ is the score of the conditional forward process $x_t|x_0$, which is analytically tractable. Vincent \cite{vincent2011connection} showed that the DSM objective is equivalent to minimizing the expected squared error with respect to the true score
\begin{equation}
    \label{eq:TrueScoreMatchingObjective}
    \mathbb{E}_{t \sim \mathcal{U}[0,T]}
    \Bigl[\lambda\left(t\right)\,\mathbb{E}_{\Forward_t}
    \|s_\theta\left(t,\Forward_t\right) - s(t,x_t)\|^2\Bigr].
\end{equation}
Thus, $s_\theta$ trained with the tractable objective \eqref{eq:DSM} is indeed approximating the typically intractable $s$. With a score approximation at hand, we generate approximate samples from $\mathcal{P}_0$ by simulating any of the reverse-time dynamics \eqref{eq:ReferenceReverseSDE} or \eqref{eq:ProbabilityFlowODE}, via a discretization method, with $s$ replaced by $s_\theta$ and $\mathcal{L}(x_T)$ replaced by $\mathcal{N}(0,\sigma_T^2I_d)$. As a result, diffusion models inherit three central sources of error:
\begin{itemize}
    \item[(i)] an \emph{initialization error}, arising from starting the reverse process at a homogeneous Gaussian distribution $\mathcal{N}\left(0,\sigma^2I_d\right)$ instead of $\mathcal{L}\left(x_T\right)$,
    \item[(ii)] a \emph{score-matching error}, due to replacing the true score $s$ with a neural network approximation $s_\theta$, and
    \item[(iii)]  a \emph{discretization error}, introduced by the numerical scheme used to simulate the not explicitly solvable reverse-time dynamics.
\end{itemize}
In this work, we provide an error analysis for general discretization methods applied to the SDE \eqref{eq:ReferenceReverseSDE} used in generative diffusion models which includes all three error contributions. We complement the theoretical analysis with extensive simulations on both toy examples and standard image-generation benchmarks. Beyond providing an environment to easily prove the convergence order of existing and novel solvers, our results also lead to practical parameter-selection guidelines, helping to identify where computational resources are best invested when designing new diffusion models.

\subsection{Related Works}
\label{sec:RelatedWork}
The theoretical literature on convergence guarantees for diffusion models splits along which dynamics are used in generation: works analyzing the deterministic probability flow ODE, and works analyzing the reverse SDE itself.

\paragraph{ODE-based samplers.} The deterministic route, i.e. discretizing the probability flow ODE \eqref{eq:ProbabilityFlowODE}, has become the default choice for fast sampling outside of knowledge distillation methods: without a stochastic term, classical higher-order solvers such as Heun's method \cite{karras2022elucidating,Karras2024edm2} or dedicated exponential multistep schemes \cite{lu2022dpmsolver, zhang2023deis,lu2025dpm} apply out of the box, and typically produce acceptable samples with far fewer score evaluations than stochastic samplers. A corresponding convergence theory has developed quickly: polynomial-time guarantees for the probability flow ODE with stochastic corrector steps \cite{chen2023probabilityflowODE}, sharp guarantees for fully deterministic discrete-time samplers \cite{li2024sharpODE}, Wasserstein guarantees for general forward processes \cite{gao2024convergenceODE} and under weak log-concavity \cite{kremling2025weaklogconcave}, and, recently, convergence theory for genuinely higher-order ODE solvers \cite{wang2025highorderODE}. However, rather than ODE-based sampling being strictly better than SDE-based sampling, the two families exhibit complementary behavior in empirical studies: ODE-based samplers are fast but their sample quality plateaus as the computational budget grows, whereas SDE-based samplers continue to improve and attain better final quality \cite{xu2023restart, karras2022elucidating}. A common explanation is that the noise in \eqref{eq:ReferenceReverseSDE} continually contracts errors accumulated in earlier steps, while the deterministic dynamics \eqref{eq:ProbabilityFlowODE} transports them to the output \cite{xu2023restart}. The ODE's practical advantage at small computational budgets, in turn, is conventionally attributed to its solvers' smaller discretization error. An advantage that is usually taken for granted, because for the SDE out-of-the-box higher-order attempts have shown disappointing empirical results \cite{jolicoeur2021gotta, karras2022elucidating}, and the theory of higher-order stochastic samplers has remained essentially undeveloped. Notably, the theoretical comparison in \cite{xu2023restart,beyler2025convergence} models the SDE discretization error as $\O(\sqrt h)$ for Euler--Maruyama against the ODE's $\O(h)$ when using the Euler method. Our results, together with \cite{wang2026linearorderEM}, revise precisely this premise: in the 2-Wasserstein distance, the Euler--Maruyama discretization of the reverse SDE already attains order $\O(h)$ (cf.\ \cite{wang2026linearorderEM} and Section~\ref{sec:ExplicitSchemes}), and the methods we analyze surpass it, attaining order-$1.5$. This narrows the presumed discretization gap between the two families and, in our view, makes the stochastic sampler, with its error-contracting properties, an underexplored rather than inferior choice. Understanding whether, and by which schemes, the stochastic sampler can attain provably higher order is therefore a natural step towards combining the speed of ODE samplers with the favorable absolute sample quality of SDE samplers, and it is the central focus of this work.

\paragraph{SDE-based samplers.} Convergence guarantees for stochastic samplers have primarily been derived in TV or KL divergence via Girsanov-based arguments \cite{chen2023samplingeasylearningscore, chen2023improved, benton2023nearly, lee2022convergence}, and in discrete time via elementary arguments \cite{li2023towards}. These analyses accommodate very general data distributions, but are structurally tied to first-order accuracy: as observed in \cite[Appendix A.1]{SDE-DPM-2-paper}, the prevalent proof technique caps the achievable KL order independently of the order of the underlying scheme, and converting KL bounds into $\mathcal W_2$ via Talagrand-type inequalities, while not even being possible in every scenario, costs a further half order. More recently, several works have obtained guarantees directly in the 2-Wasserstein distance \cite{bruno2025OU_convergence_Log_Concave, gao2025wassersteinconvergenceguaranteesgeneral, strasman2025NonAsymptoticBounds, yu2026advancing, gentilonisilveri2025OU_Convergence, beyler2025convergence, wang2026linearorderEM}, which is known to correlate well with human perceptual similarity \cite{rubner20002WassersetinCorrespondsToHumanPerception} and is closely related to the widely used FID \cite{heusel2017gans} (cf. Section \ref{subsec:MeasureingW2}). Within this line, most works restrict attention to the Ornstein--Uhlenbeck forward process \cite{bruno2025OU_convergence_Log_Concave, yu2026advancing, gentilonisilveri2025OU_Convergence, wang2026linearorderEM}. Only \cite{gao2025wassersteinconvergenceguaranteesgeneral} treats the more general coefficient functions $f, g$ (cf. Section \ref{subsec:ProblemSetup}), which we consider essential given the demonstrated impact of the noise schedule on generation quality \cite{chen2023importanceOfNoiseschedule}, yet obtains only the rate $\O(\sqrt h)$ for its Euler-type scheme. Crucially, no existing work establishes a stochastic sampler with a provably \emph{strictly better} end-to-end convergence order than Euler--Maruyama. This is visible at a glance in Table~\ref{tab:IterationComplexity} in Section \ref{subsec:MainResult}. The stochastic Runge--Kutta method of \cite{wu2024stochastic} improves the dimension dependence in KL, but uses a single score evaluation per step and remains at first order. The second-order methods analyzed in \cite{SDE-DPM-2-paper} attain strong order one, and the same work explicitly lists the sampling complexity of strong order-$1.5$ schemes such as the SRA family \cite{rossler2010runge} as unexplored \cite[Appendix A.1]{SDE-DPM-2-paper}. On the empirical side, higher order stochastic schemes were tested for diffusion models in \cite{jolicoeur2021gotta} and \cite{blasingame2025rex} without accompanying guarantees. Our order-$1.5$ result closes precisely this gap.

Several practical stochastic samplers come without \emph{end-to-end} convergence guarantees, but with convergence-order proofs for the discretization error alone at fixed terminal time, e.g.\ SA-Solver \cite{xue2023sasolver} and SEEDS \cite{gonzalez2023seeds}. These results are complementary to ours: one-step members of such families can be brought under Theorem~\ref{thm:Main} by verifying the conditions of Definition~\ref{def:generalstrongscheme}, in the same way as we do for the exponential integrator in Section~\ref{sec:ExplicitSchemes}, whereas multistep schemes would require extending the definition to remainders that are centered with respect to earlier $\sigma$-algebras.

\paragraph{Score-matching assumptions.} A further point of divergence in the literature concerns the assumption on the score-matching error. Most works, e.g.\ \cite{gao2025wassersteinconvergenceguaranteesgeneral, yu2026advancing, strasman2025NonAsymptoticBounds}, assume the error to be bounded along the iterates of the discretization scheme. As argued in \cite{beyler2025convergence}, this assumption is disconnected from training: the denoising score-matching objective \cite{vincent2011connection} controls the expected squared error along the marginals of the \emph{forward} process, or equivalently, of the true reverse process, whereas the law of the sampler's iterates differs from these marginals by exactly the initialization and discretization errors the theorems seek to bound, so a small training loss does not imply the assumed bound. We therefore adopt the training-aligned assumption and quantify its price explicitly (Lemma~\ref{lem:ScoreErrorTransfer}). To our knowledge, the resulting interplay between initialization decay and score-error transfer has not been made explicit before.

\subsection{Our Contribution}
Our main contribution is an all-at-once analysis of generative diffusion models, i.e.\ an analysis that simultaneously covers all three error sources, initialization, score matching, and discretization, see Theorems~\ref{thm:Main} and~\ref{thm:MainDissipative}. The central ingredient of our analysis is the notion of a \emph{\SchemeClassName} (Definition~\ref{def:generalstrongscheme}), together with a single convergence theorem in the 2-Wasserstein distance that covers this entire class of discretization methods at once (cf. Theorem \ref{thm:Main}) and exploits dissipativity of the reverse-time drift whenever it is available (cf. Theorem \ref{thm:MainDissipative}). This turns the convergence analysis of any concrete sampler into the verification of a finite checklist: matching the scheme's terms against the coefficient functions of the underlying Itô--Taylor expansion up to explicit, index-wise tolerances. We carry out this verification for the Euler--Maruyama scheme, for the exponential integrator, and for the explicit SRA family of solvers introduced by R\"o\ss{}ler \cite[Section 6.5]{rossler2010runge}, a derivative-free class of stochastic Runge--Kutta schemes of strong order $1.5$ for SDEs with additive noise. To the best of our knowledge, this result is the first to offer full convergence guarantees for schemes exhibiting strictly higher order than Euler--Maruyama in the setting of generative diffusion models (cf. Table~\ref{tab:IterationComplexity}). Furthermore, in contrast to most existing results, our analysis covers a broad class of forward SDEs with general spatially linear drift with time-dependent coefficient $f$, and spatially independent diffusion with time-dependent coefficient $g$ (cf. Section \ref{subsec:ProblemSetup}), rather than a fixed variance-preserving, variance-exploding or Ornstein--Uhlenbeck schedule. Keeping the analysis general in this respect is important, as it has been shown that choosing an optimal \emph{noise schedule}, that is, an appropriate choice of the coefficient functions, can significantly improve the quality of the matched score function, and thus the performance of generative diffusion models \cite{chen2023importanceOfNoiseschedule}. Moreover, our analysis is modular in the initialization error which enters only through the abstract bound of Assumption \ref{ass:data}, so that existing and future results on the convergence of the forward process, e.g.\ \cite{gao2025wassersteinconvergenceguaranteesgeneral}, can be inserted directly.

While our proofs rely on standard results of Kloeden and Platen \cite{kloeden2010numerical} on strong solutions and moments of multiple stochastic integrals, our argument departs from the classical one in two related respects, which together we regard as a contribution in its own right. First, since $T$ must grow for the initialization error to vanish, we track the dependence on $T$ explicitly throughout the proofs. In standard references such as \cite{kloeden2010numerical,milstein2004stochastic}, the time horizon is fixed and this dependence is absorbed into generic constants. Second, our proof retains the sign of the cross term between the error and the drift difference in the discretization recursion (cf.\ \eqref{eq:DriftCrossTerm} in Appendix~\ref{sec:ProofOfMainTheorem}), which the classical argument usually discards. Without further structure, the resulting bounds still grow exponentially in $T$ (cf. Theorem \ref{thm:Main}). Under a dissipative drift, however, the retained sign reduces the dependence of the discretization and score-matching errors on $T$ from exponential to polynomial (cf. Theorem \ref{thm:MainDissipative}), and admits more general initialization errors (cf. Remark~\ref{rem:Dissipativity}). Uniform-in-time strong error bounds under related dissipativity assumptions are known for Euler-type schemes \cite{fang2020adaptive}, and our results complement these for higher-order samplers. We further discuss our results in the light of non-asymptotic bounds (Remark \ref{rem:NonAsymptotic}) and derive the resulting iteration complexity together with an optimal parameter selection rule for the terminal time, the score accuracy, and the step size (Corollary \ref{cor:IterationComplexityDissipative}), improving the complexity exponent by the factor $1/\gamma$ for schemes of order $\gamma$.

Finally, we present extensive numerical experiments demonstrating that the predicted convergence orders are directly observable in practice, both for toy problems where the data distribution is a Gaussian (mixture) and for the real-world datasets CIFAR-10 \cite{krizhevsky2009CIFAR10} and ImageNet \cite{deng2009imagenet}. In contrast to experiments reporting only downstream sample quality, our experiments verify the theoretical rates explicitly. Moreover, contrary to popular believe, they show that out-of-the-box higher-order stochastic samplers are highly profitable in practice: compared to Euler--Maruyama, SRA1 reduces the number of NFEs required to reach the error floor by a factor of $8$ to $16$ on the toy problems, and by at least a factor of four on ImageNet.

\subsection{Notation and Problem Setup}
\subsubsection{Problem Setup}
\label{subsec:ProblemSetup}
Throughout this work, we study forward processes whose drift is linear in the spatial
variable, $a(t,x)=-f(t)x$, and whose diffusion coefficient is a scalar
multiple of the identity, $b(t)=g(t)I_d$, with scalar-valued
$f,g\colon[0,\infty)\to\R$, i.e.\ processes of the form
\begin{equation}\label{eq:Forward}
    \d\Forward_t = -f\left(t\right)\Forward_t\,\d t + g\left(t\right)\d W_t,
    \qquad \Forward_0 \sim \mathcal{P}_0,
\end{equation}
which include many diffusion models used in practice, cf.\ \cite{song2020score,dockhorn2021score,gao2025wassersteinconvergenceguaranteesgeneral} and Table \ref{tab:PopularChoices}.
Using the time-reversal formula \eqref{eq:ReferenceReverseSDE}, the corresponding reverse-time SDE is
\begin{equation}\label{eq:Backward}
    \d\Backward_t
    = \left[f\left(\tau\right)\Backward_t + g^2\left(\tau\right)\,s\left(\tau,\Backward_t\right)\right]\d t
      + g\left(\tau\right)\d W_t,
    \qquad \Backward_0 \sim \mathcal{L}\left(\Forward_T\right).
\end{equation}
Furthermore, conditioned on $x_0$, the solution of the SDE \eqref{eq:Forward} can be derived explicitly as $x_t|x_0 \sim \mathcal{N}\left(\phi_f(t)x_0, \varphi_{f,g}(t)I_d\right)$, where
\begin{equation}
    \label{eq:StationarySolutionForward}
        \phi_f\left(t\right) = e^{\left(-\int_0^tf\left(s\right)\d s\right)} \quad \text{and } \quad \varphi_{f,g}\left(t\right) =\int_0^te^{-2\int_s^tf\left(v\right)\d v}g^2\left(s\right)\d s.
\end{equation}
This closed form is not only needed to make the DSM training realizable, but also allows us to choose $f$ and $g$ such that the terminal distribution is a uniform Gaussian with diagonal covariance and mean zero. In fact, we consider the backward process to be initialized from \(\mathcal{N}\!\left(0, \sigma_T^2 I_d\right)\), with \(\sigma_T^2 := \varphi_{f,g}\left(T\right)\). We proceed by introducing the following processes:
\begin{align}
    \label{eq:BackwardNormalstartingSDE}\d\BackwardNormalstarting_t
    &= \left[f(\tau)\BackwardNormalstarting_t + g^2(\tau)s\left(\tau,\BackwardNormalstarting_t\right)\right]\d t
      + g(\tau)\d W_t,
    &&\BackwardNormalstarting_0 \sim \mathcal{N}(0,\sigma^2_TI_d),\\
    \label{eq:BackwardNetworkScoreSDE}\d\BackwardNetwork_t
    &= \left[f(\tau)\BackwardNetwork_t + g^2(\tau)s_\theta\left(\tau,\BackwardNetwork_t\right)\right]\d t
      + g(\tau)\d W_t,
    &&\BackwardNetwork_0 \sim \mathcal{L}(\Forward_T),\\
    \label{eq:BackwardNormalstartingScoreSDE}\d\BackwardNormalstartingNetwork_t
    &= \left[f(\tau)\BackwardNormalstartingNetwork_t + g^2(\tau)s_\theta\left(\tau,\BackwardNormalstartingNetwork_t\right)\right]\d t
      + g(\tau)\d W_t,
    &&\BackwardNormalstartingNetwork_0 \sim \mathcal{N}(0,\sigma^2_TI_d),
\end{align}
where all processes \eqref{eq:Backward}-\eqref{eq:BackwardNormalstartingScoreSDE}
are driven by the same Wiener process $W_t$, and $\Backward_0=\BackwardNetwork_0$ as well as $\BackwardNormalstarting_0=\BackwardNormalstartingNetwork_0$ almost surely. We choose $\Backward_0$, $\BackwardNormalstarting_0$ and $W$ to be mutually independent.\footnote{Since the 2-Wasserstein distance $\mathcal W_2$ is defined as an infimum over all couplings of two laws, any explicit coupling of the processes yields a valid upper bound on the distances of interest. Coupling through a shared $W_t$ is what makes the pairwise differences of these processes tractable, since the diffusion terms cancel exactly, as exploited in the proof of Lemma~\ref{lem:ScoreMatchingError3}. The joint law of $\Backward_0$ and $\BackwardNormalstarting_0$, on the other hand, plays no role in our estimates:these two processes are only ever compared through their laws (Assumption~\ref{ass:data} and Lemma~\ref{lem:ScoreErrorTransfer}), so independence is merely a concrete choice.} One can interpret $\BackwardNormalstarting_t$ as being the process which introduces the initialization error, $\BackwardNetwork_t$ as being the process which introduces the score matching error, and $\BackwardNormalstartingNetwork_t$ as being the process which exhibits both error types. Furthermore, we denote the discretized version of any time continuous process $\ReferenceSDE_t$ by $\Discrete_t$ or $\hat{x}_t$, depending on the discretization used (cf. Definitions \ref{def:ItôTaylorScheme} and \ref{def:generalstrongscheme}). Thus, the processes $\BackwardNormalstartingNetworkDiscretizedGeneralStrong_k$ and $\BackwardNormalstartingNetworkDiscretizedItoTaylor_k$ exhibit all three error types. All discretizations are started from the initial value of the process they approximate, i.e.\ $\BackwardNormalstartingNetworkDiscretizedGeneralStrong_0=\BackwardNormalstartingNetworkDiscretizedItoTaylor_0=\BackwardNormalstartingNetwork_0$, and use the increments of the same Wiener process $W$. To simplify notation, we consider an equidistant discretization grid $0=t_0<t_1<...<t_K=T$ with uniform step size $h=\tfrac{T}{K}$, and write $\tau_k=T-t_k$ for the reverse time discretization. However, we want to stress that all proofs also work on shrinking non-equidistant discretizations. Moreover, by $\|\cdot\|$, we denote the Euclidean norm, and $\|\cdot\|_{L_2}= \sqrt{\mathbb{E}[\|\cdot\|^2]}$. Components of a vector $u\in\R^d$ are denoted by $u^{(i)}$. Furthermore, for a sufficiently smooth function $\varphi : \mathbb{R} \times \mathbb{R}^d \to \mathbb{R}^{d'}$, $d'\in\mathbb N$, we write $\partial_t\varphi:=\frac{\partial\varphi}{\partial t}$, and denote by $D\varphi(t,x)\in\R^{d'\times d}$ and $D^2\varphi(t,x)$ its first and second derivative with respect to $x$, that is, for $u,v\in\R^d$,
\begin{equation}
    \label{eq:DifferentialOperatorDefinition}
    D\varphi(t,x)u:=\frac{\mathrm d}{\mathrm d\varepsilon}\varphi(t,x+\varepsilon u)\Big|_{\varepsilon=0},
    \qquad
    D^2\varphi(t,x)[u,v]:=\frac{\mathrm d}{\mathrm d\varepsilon}D\varphi(t,x+\varepsilon v)u\Big|_{\varepsilon=0}.
\end{equation}
Thus, $D\varphi(t,x)$ is the Jacobian matrix of $\varphi(t,\cdot)$ at $x$, and $D^2\varphi(t,x)$ is a symmetric bilinear map whose values $D^2\varphi(t,x)[e_i,e_j]$ are the second partial derivatives with respect to the $i$-th and $j$-th spatial coordinate, where $e_i$ denotes the $i$-th unit vector. Furthermore, we write $\text{Poly}(T)$ for an arbitrary polynomial in $T$ that can vary from line to line, i.e. there exists a finite $p\in \mathbb N$ such that  $\text{Poly}(T) \leq C(1+T)^p $ with $C$ independent of $T$. Throughout, $(\Omega,\mathcal F,\mathbb P)$ carries the $d$-dimensional Wiener process $W$ and the initial values $\Backward_0$ and $\BackwardNormalstarting_0=\BackwardNormalstartingNetwork_0$, which are independent of $W$. We set
\[
    \mathcal F_t := \sigma\big(\Backward_0,\ \BackwardNormalstartingNetwork_0,\ W_s : 0\le s\le t\big)\vee\mathcal N, \qquad t\in[0,T],
\]
where $\mathcal N$ denotes the $\mathbb P$-null sets. Then $W$ is an $(\mathcal F_t)$-Wiener process and the solutions of \eqref{eq:Backward}-\eqref{eq:BackwardNormalstartingScoreSDE} are $(\mathcal F_t)$-adapted.

\subsubsection{Itô--Taylor Notation}
Since our analysis relies on the It\^o--Taylor expansion, we introduce the standard
notation of Kloeden and Platen \cite{kloeden2010numerical}. We write $\mathcal M$ for the set of all multi-indices
\[
    \alpha = (j_1, \ldots, j_\ell), \qquad j_k \in \{0, 1, \ldots, d\},
\]
of arbitrary length $\ell\ge1$, together with the empty multi-index $\nu$, and let $\mathcal A\subseteq\mathcal M$ denote a generic subset of multi-indices. Here $j_k = 0$ corresponds to integration
with respect to time, while $j_k \in \{1, \ldots, d\}$ corresponds to integration with
respect to the $j_k$-th component of the $d$-dimensional Wiener process $W_t$. For a
multi-index $\alpha = (j_1, \ldots, j_\ell)$, we define its length and number of non-stochastic components by
\[
    l(\alpha) := \ell,
    \qquad
    n(\alpha) := \sum_{k=1}^{\ell} \mathbf{1}_{\{j_k = 0\}},
\]
and write $-\alpha := (j_2, \ldots, j_{\ell})$ for the multi-index obtained by
removing the first component, and $\alpha- := (j_1, \ldots, j_{\ell-1})$ for the multi-index obtained by removing the last component, with the convention $-\alpha = \alpha- = \nu$ when $\ell = 1$.
For a given hierarchical set $\mathcal{A}$, i.e.\ a set satisfying $-\alpha\in\mathcal A$ for every $\alpha \in \mathcal A\setminus\{\nu\}$, its remainder set is defined by
\[
    \mathcal{B}(\mathcal{A})
    := \{\, \alpha \notin \mathcal{A} : -\alpha \in \mathcal{A} \,\}.
\]
The multiple It\^o integrals are defined recursively over an interval $[s, t]$ with
$0 \leq s < t \leq T$. For a sufficiently regular adapted process
$h : \Omega \times [0,T] \to \mathbb{R}$, set
\[
    I_\nu[h]_{s,t} := h(s),
\]
and for $\alpha = (j_1, \ldots, j_\ell)$ with $\ell \geq 1$, define recursively
\[
    I_\alpha[h]_{s,t}
    :=
    \begin{cases}
        \displaystyle
        \int_s^t I_{\alpha-}[h]_{s,u} \, \d u,
        & j_\ell = 0, \\[2ex]
        \displaystyle
        \int_s^t I_{\alpha-}[h]_{s,u} \, \d W_u^{(j_\ell)},
        & j_\ell \in \{1, \ldots, d\}.
    \end{cases}
\]
For a vector-valued process $h = (h^{(1)}, \ldots, h^{(d)})^\top$ we extend this definition
componentwise by setting
\[
    I_\alpha[h]_{s,t}
    := \bigl(I_\alpha[h^{(1)}]_{s,t},\,\ldots,\,I_\alpha[h^{(d)}]_{s,t}\bigr)^\top,
\]
so that $I_\alpha[f_\alpha]_{s,t} \in \mathbb{R}^d$ in the expansion below. We write $h\in\mathcal H_\alpha$ to indicate that $h$ satisfies the measurability and integrability conditions under which $I_\alpha[h]_{s,t}$ above is well-defined for every $0\le s<t\le T$, see \cite[Section 5]{kloeden2010numerical} for the precise conditions.

Consider the $d$-dimensional SDE \eqref{eq:BackwardNormalstartingScoreSDE}, i.e. the SDE exhibiting both the initialization error and score matching error, driven by
a $d$-dimensional Wiener process, with drift
\[
    \mu(t, x) = f(t)\, x + g^2(t)\, s_\theta(t, x),
\]
and scalar diffusion coefficient $g(t) \in \mathbb{R}$, so that the diffusion matrix
is $g(t) I_d$. We define the differential operators $L^0$ and $L^r$ acting on such functions by
\[
    L^0 \varphi(t, x)
    := \frac{\partial \varphi}{\partial t}(t, x)
    + D\varphi(t, x)\, \mu(t, x)
    + \frac{1}{2} g(t)^2 \sum_{i=1}^d D^2\varphi(t, x)[e_i,e_i],
\]
and for $r = 1, \ldots, d$,
\[
    L^r \varphi(t, x)
    := g(t)\, D\varphi(t, x)\,e_r.
\]
For a multi-index $\alpha = (j_1, \ldots, j_\ell)$, the coefficient functions
$f_\alpha : \mathbb{R} \times \mathbb{R}^d \to \mathbb{R}^d$ are defined
by the recursion
\begin{equation}
    \label{eq:CoefficientFunctions}
    f_\nu(t, x) := x,
    \qquad
    f_\alpha(t, x) := L^{j_1} f_{-\alpha}(t, x).
\end{equation}
With this notation in place, the It\^o--Taylor expansion of \eqref{eq:BackwardNormalstartingScoreSDE}
over one time step $[t_k, t_{k+1}]$ reads a.s.
\[
    \BackwardNormalstartingNetwork_{t_{k+1}}
    = \BackwardNormalstartingNetwork_{t_k}
    + \sum_{\alpha \in \mathcal{A}_\gamma \setminus \{\nu\}}
        I_\alpha\bigl[f_\alpha(t_k, \BackwardNormalstartingNetwork_{t_k})\bigr]_{t_k,\, t_{k+1}}
    + \sum_{\alpha \in \mathcal{B}(\mathcal{A}_\gamma)}
        I_\alpha\bigl[f_\alpha(\,\cdot\,, \BackwardNormalstartingNetwork_{\,\cdot\,})\bigr]_{t_k,\, t_{k+1}},
\]
where the first sum collects the frozen-coefficient terms of the scheme and the second
sum is the hierarchical remainder. Strong It\^o--Taylor schemes are obtained by
discarding the remainder.

\begin{defi}
    \label{def:ItôTaylorScheme}
    We define a strong Itô--Taylor scheme for \eqref{eq:BackwardNormalstartingScoreSDE} of order $\gamma$ as
    \[\BackwardNormalstartingNetworkDiscretizedItoTaylor_{k+1}=\BackwardNormalstartingNetworkDiscretizedItoTaylor_k+\sum_{\alpha \in \mathcal A_\gamma\backslash\{\nu\}}I_\alpha[f_\alpha(t_k,\BackwardNormalstartingNetworkDiscretizedItoTaylor_k)]_{t_k,t_{k+1}},\]
    where
    \[\mathcal A_\gamma = \{\alpha \in \mathcal M: l(\alpha)+n(\alpha)\leq 2 \gamma \text{  or  } l(\alpha)=n(\alpha)=\gamma +0.5\}.\]
\end{defi}
\begin{defi}
    \label{def:generalstrongscheme}
    We define an \emph{\SchemeClassName} (\SchemeClassShortHand) for \eqref{eq:BackwardNormalstartingScoreSDE} of order $\gamma$ as
    \[\BackwardNormalstartingNetworkDiscretizedGeneralStrong_{k+1}=\BackwardNormalstartingNetworkDiscretizedGeneralStrong_k+\sum_{\alpha \in \mathcal A_\gamma\backslash\{\nu\}}I_\alpha[g_{\alpha,k}]_{t_k,t_{k+1}} + R_k,\]
    with $\mathcal{A}_\gamma$ from Definition \ref{def:ItôTaylorScheme} and
    \[\mathbb{E}\Big[\max_{0\leq k \leq K-1}\|g_{\alpha,k}-f_\alpha(t_k,\BackwardNormalstartingNetworkDiscretizedGeneralStrong_k)\|^2\Big]\leq C h^{2\gamma -\phi(\alpha)} \]
    with
    \[\phi(\alpha)=\begin{cases}
        2(l(\alpha)-1) & \text{if}\quad  l(\alpha)=n(\alpha)\\
        l(\alpha)+n(\alpha)-1&\text{if}\quad l(\alpha)\neq n(\alpha)
    \end{cases}.\]
    We further require that the coefficients $g_{\alpha , k}$ are $\mathcal{F}_{t_k}$ measurable, and the remainder $R_k$ to admit a decomposition
    \[
        R_k = R_k^M + R_k^D
    \]
    into a (martingale) term $R_k^M$ that is centered conditionally on $\mathcal F_{t_k}$,
    \[
        \mathbb E\big[R_k^M \mid \mathcal F_{t_k}\big] = 0 \quad \text{a.s.}, \qquad
        \max_{0\leq k \leq K-1}\mathbb{E}\left[\big\|R_k^M\big\|^2\right] \leq C_R\,h^{2\gamma+1},
    \]
    and a higher order term $R_k^D$, in the sense that
    \[
        \max_{0\leq k\leq K-1}\mathbb{E}\left[\big\|R_k^D\big\|^2\right]\leq C_R\,h^{2\gamma+2}.
    \]
    Furthermore, both remainder components $R_k^M$ and $R_k^D$ need to be $\mathcal{F}_{t_{k+1}}$ measurable.
\end{defi}

\begin{rem}
    Since we are interested in the behavior of \SchemeClassShortHand with varying terminal time $T$, the constant in Definition \ref{def:generalstrongscheme} need to follow a specific $T$ dependence themselves. In particular, we assume them to be $C(T),C_R(T)=\O(e^{\text{Poly}(T)})$ in general, and in the dissipative case (cf. Assumption \ref{ass:Dissipativity}), we assume them to grow at most polynomial in $T$. Furthermore, we note that Definitions \ref{def:ItôTaylorScheme} and \ref{def:generalstrongscheme} are specifically formulated for the SDE \eqref{eq:BackwardNormalstartingScoreSDE}, because we defined the coefficient functions $f_\alpha$ for this SDE. The definitions can of course be applied to SDEs which do not include an initialization or score-matching error. 
\end{rem}
 The only difference of our definition to Kloeden and Platen's \cite{kloeden2010numerical} is in the condition on the remainder term $R_k$, which they bound at once via
\[
    \mathbb{E}\Bigg[\sup_{0\leq k \leq K-1}\Big\|\sum_{n=0}^k R_n\Big\|^2\Bigg]
    \leq C h^{2\gamma}.
\]
This difference is a direct consequence of the fact that we choose to bound the error of the scheme iterate by iterate rather than on the whole interval at once. Our definition is in fact the stronger of the two (cf. Theorem \ref{rem:GeneralStrongSchemeDefi}), yet it considerably simplifies verifying the conditions for a specific scheme, because we remove one step in the usual proof of the remainder condition. This fact is exploited in Appendix \ref{sec:MainOrderProofs}, where we also prove that the explicit SRA family is an \SchemeClassShortHand, whereas, to the best of our knowledge, it was previously not known whether this family belongs to the general strong schemes.
\begin{thm}
   
\label{rem:GeneralStrongSchemeDefi} 
Definition \ref{def:generalstrongscheme} it is stronger than the notion of a general strong scheme of Kloeden and Platen \cite[Section 11.5]{kloeden2010numerical}.
\end{thm}
Furthermore, we will see later that our definition of an \SchemeClassShortHand relates naturally to Milstein's fundamental mean-squared theorem \cite[Theorem 1.1]{milstein2004stochastic}, cf. Remark \ref{rem:MilsteinRelation}.

\section{Assumptions and Main Result}
\label{sec:main}
In this section, we present our main result. We begin by defining and discussing the assumptions used throughout this work.

\subsection{Assumptions and Merit Discussion}
\label{sec:assumptions}
To analyze the reverse-time SDE and its discretization, we impose several structural assumptions.
\begin{assumption} (existence of solution)
    \label{ass:Solution}
    We assume $\|x_0\|_{L_2}<\infty$ for $x_0 \sim \mathcal{P}_0$ and all SDEs \eqref{eq:Forward}, \eqref{eq:Backward}, \eqref{eq:BackwardNormalstartingSDE}, \eqref{eq:BackwardNetworkScoreSDE}, and \eqref{eq:BackwardNormalstartingScoreSDE} admit strong solutions.
\end{assumption}

\begin{assumption} (bounded initialization error)
    \label{ass:data} We assume  that the initialization error is bounded by 
    \[\sup_{0\leq t\leq T}\mathcal W_2\big(\mathcal L(\Backward_t),\mathcal L(\BackwardNormalstarting_t)\big)\leq B_{\text{init}}(T), \qquad \lim_{T\to\infty}B_\text{init}(T)=0.\]
\end{assumption}

\begin{assumption}(score matching error bound) 
\label{ass:ScoreMatchingErrorBound} We assume the existence of some $\varepsilon>0$ such that
    \[\sup_{t \in [0,T]} \mathbb{E}\big[\bigl\|s\left(T-t,\Backward_t\right)-s_\theta\left(T-t,\Backward_t\right)\bigr\|^2\big] \leq \varepsilon^2.\]

\end{assumption}

\begin{assumption} (regularity of coefficient functions)
    \label{ass:CoefficientFunctions}
    For a given index set $\mathcal A$, we assume the coefficient functions $f_\alpha$ defined in \eqref{eq:CoefficientFunctions} to exist, to be Lipschitz continuous and of linear growth, i.e. there exists a constant $C_{\mathcal{A}}(T)>0$ growing at most polynomial in $T$, such that
    \begin{align*}
        \|f_\alpha(t,x)-f_\alpha(t,y)\|^2\leq C_{\mathcal{A}}(T)^2\|x-y\|^2
    \end{align*}
    and
    \begin{align*}
        \|f_\alpha(t,x)\|^2\leq C_{\mathcal{A}}(T)^2(1+\|x\|^2),
    \end{align*}
    for all $\alpha \in \mathcal A$, $x,y\in \R^d$, $t\in [0,T]$.
\end{assumption}

\begin{assumption} (regularity of score function)
    \label{ass:ScoreFunctionLipschitz}
    We assume the score network to be pointwise Lipschitz continuous, i.e. there exists a constant $L_{s_\theta}>0$, such that
    \[\|s_{\theta}(t,x) - s_{\theta}(t,y)\|\leq L_{s_\theta}(T)\|x-y\|\]
    for all $t\in[0,T]$, $x,y\in\R^d$, and $L_{s_\theta}(T)$ growing at most polynomial in $T$.
\end{assumption}

\begin{assumption}(regularity of the true score)
\label{ass:TrueScoreLipschitz}
We assume the true score function to be pointwise Lipschitz continuous, i.e. there exists $L_s>0$ such that \[\|s(t,x)-s(t,y)\|\le L_s(T)\|x-y\|\]
for all $t\in[0,T]$, $x,y\in\R^d$, and $L_{s}(T)$ growing at most polynomial in $T$.
\end{assumption}

\begin{assumption} (model parameters)
\label{ass:DriftAndDiffusion} In general, we assume $f,g:[0, \infty) \to \R$ to be deterministic, to be chosen such that the terminal law of the forward process approaches $\mathcal{N}(0,\sigma_T^2I_d)$ with growing $T$ and with $\sigma_T^2$ growing at most polynomial in $T$, and $g$ to be bounded via $M_g(T):=\sup_{t\in [0,T]}|g(t)|$.
\end{assumption}

\begin{assumption} (dissipative drift)
\label{ass:Dissipativity}
There exists $\lambda>0$ such that the drift $\mu(t,x)=f(T-t)x+g^2(T-t)s_\theta(T-t,x)$ of \eqref{eq:BackwardNormalstartingScoreSDE} satisfies
\[
    \langle \mu(t,x)-\mu(t,y),\,x-y\rangle \le -\lambda\|x-y\|^2
    \qquad\text{for all } t\in[0,T],\; x,y\in\R^d.
\]
This assumption is imposed only for the dissipative case, cf. Theorem~\ref{thm:MainDissipative}.
\end{assumption}

We proceed by discussing the practicality of the assumptions stated above.
    \begin{description}

    \item[On Assumption \ref{ass:Solution}:] This is standard in the SDE discretization literature. We note that for the score-based reverse processes \eqref{eq:BackwardNetworkScoreSDE} and \eqref{eq:BackwardNormalstartingScoreSDE}, Lipschitz continuity and linear growth of the drift and diffusion coefficients, implied by Assumptions~\ref{ass:CoefficientFunctions}, and the finite second moment of their initial condition, already guarantee the existence of a unique strong solution. Also, Assumption~\ref{ass:DriftAndDiffusion}, together with $\|x_0\|_{L_2}<\infty$ assures a unique strong solution for the forward SDE \eqref{eq:Forward}. The assumption is thus only genuinely restrictive for the reverse processes \eqref{eq:Backward} and \eqref{eq:BackwardNormalstartingSDE}, whose drifts involve the true score, which we do not explicitly assume to be linearly growing.

    \item[On Assumption \ref{ass:data}:] This assumption is kept abstract on purpose: our error decomposition is modular, and any bound on the initialization error that holds uniformly in $t\in[0,T]$ can be inserted into Theorem~\ref{thm:Main} or Theorem~\ref{thm:MainDissipative} directly. Uniformity is required since the score-matching error is transferred along the whole trajectory (Lemma~\ref{lem:ScoreErrorTransfer}). The most common route is strong log-concavity of the data distribution, e.g.\ \cite{bruno2025OU_convergence_Log_Concave,strasman2025NonAsymptoticBounds,yu2026advancing}, and for general $f,g$, \cite[Eq.~(3.12)]{gao2025wassersteinconvergenceguaranteesgeneral} provide exactly such a uniform bound, which decays for VP and OU schedules. Bounds at the terminal time have further been obtained for weakly log-concave data \cite{kremling2025weaklogconcave, gentilonisilveri2025OU_Convergence} and data with bounded support \cite{beyler2025convergence, lee2023ExponentialIntegrator}. We expect uniform analogues of these by the same arguments, but do not pursue this here. These bounds typically carry further technical assumptions, which we then inherit alongside them.

    We note that VE schedules, i.e.\ $f\equiv 0$, generally do not satisfy Assumption~\ref{ass:data}: since $\mathbb E[\Forward_T]=\mathbb E[\Forward_0]$, already at $t=0$ we have $B_\text{init}(T)\geq\mathcal W_2\big(\mathcal L(\Forward_T),\mathcal N(0,\sigma_T^2I_d)\big)\geq\|\mathbb E[\Forward_0]\|$ for every $T$. This is a consequence of requiring the bound uniformly in $t$, as bounds at the terminal time do decay for VE schedules \cite{gao2025wassersteinconvergenceguaranteesgeneral}. Our full error bound therefore does not cover VE schedules, in particular the EDM2 setting of Problems~\ref{prob:ImageNet64} and~\ref{prob:ImageNet512}. We nonetheless include these experiments to test our discretization error bound, which is independent of Assumption~\ref{ass:data}.

    \item[On Assumption \ref{ass:ScoreMatchingErrorBound}:] While a lot of the literature assumes the score error to be bounded along the path induced by the discretization of the reverse process $\BackwardNormalstartingNetworkDiscretized_k$, we take the bound over the path of true backward process $\Backward_t$. We believe it to be the more natural assumption as discussed in \cite{beyler2025convergence} and Section \ref{sec:RelatedWork}, as, contrary to the prevalent assumption, it is formulated along the same marginals as the learning objective.

    \item[On Assumption \ref{ass:CoefficientFunctions}:] This assumption is standard in a lot of numerical SDE literature, e.g. in Kloeden and Platen \cite{kloeden2010numerical}. It is true if all derivatives of $f,g$ and the score network $s_\theta$ are Lipschitz continuous up to a certain order. In some works, e.g. by Chang \cite{chang1987numerical}, this assumption is sometimes replaced by the same derivatives being bounded, which would be applicable here as well. One could relax this condition to sufficiently differentiable coefficient functions and existence of non exploding solutions, cf. \cite{jentzen2009pathwise}. We also note that this assumption implies the boundedness of $g$ in Assumption \ref{ass:DriftAndDiffusion}, as $M_g(T)\leq C_{\mathcal A}(T)$, because $f_{(j)}=ge_j$.

    \item[On Assumption \ref{ass:ScoreFunctionLipschitz}:] While it seems redundant in light of Assumption \ref{ass:CoefficientFunctions}, we need it in order to explicitly transfer the score matching error from the true reverse process $\Backward_t$ to the backward process including the initialization error $\BackwardNormalstarting_t$. We note that this assumption is naturally fulfilled by a network with bounded weights.
    
    \item[On Assumption \ref{ass:TrueScoreLipschitz}:] This assumption is standard in the existing literature, either through explicitly assuming it on the whole time interval, e.g. \cite{beyler2025convergence,strasman2025NonAsymptoticBounds}, only at $t=0$ and extending it to the whole time interval \cite{kremling2025weaklogconcave}, or inferring it from conditions on the data distribution, e.g. \cite{gao2025wassersteinconvergenceguaranteesgeneral,yu2026advancing}. We note, however, that requiring a Lipschitz constant uniform on $[0,T]$ is an implicit regularity assumption on the data distribution: it excludes, in particular, data supported on a lower-dimensional manifold, for which the true score is Lipschitz for every $t<T$ but $L_s$ diverges as the reverse process approaches the data. This restriction can be lifted by early stopping, i.e. terminating the reverse process at time $T-\delta$, at the cost of an additional bias of order $\sqrt{d\,\delta}$. This can be inserted into our error decomposition by standard arguments, cf. e.g. \cite{beyler2025convergence}.
    
    \item[On Assumption \ref{ass:DriftAndDiffusion}:] The terminal distribution condition can directly be fulfilled by looking at \eqref{eq:StationarySolutionForward}, whereas the boundedness condition is used in several places, e.g. in Lemma \ref{lem:ScoreMatchingError3}, and is used to ensure that the coefficients are not increasing too rapidly.

    \item[On Assumption \ref{ass:Dissipativity}:] This is an assumption on the
    drift, that is, on $f$, $g$, and the \emph{learned} score $s_\theta$,
    and imposes no direct condition on the data distribution. Its structure
    is revealed by the identity
    \begin{equation}
        \label{eq:DissipativityIdentity}
        \langle\mu(t,x)-\mu(t,y),\,x-y\rangle
        = f(\tau)\|x-y\|^2
        + g^2(\tau)\,\langle s_\theta(\tau,x)-s_\theta(\tau,y),\,x-y\rangle,
    \end{equation}
    which shows that the drift is dissipative with rate $\lambda$ if and only if $s_\theta(\tau,\cdot)$ is one-sidedly Lipschitz with constant $-\big(\lambda + f(\tau)\big)/g^2(\tau)$. For the \emph{true} score, this one-sided bound is precisely strong log-concavity of the forward marginal $p_\tau$. The assumption thus couples the noise schedule and the score through the threshold $f(\tau)/g^2(\tau)$. Whether it holds for a \emph{learned} score is then a question of inheritance: if $s_\theta$ approximates a true score satisfying the bound sufficiently well, it satisfies a comparable bound. Alternatively, the one-sided Lipschitz property can be enforced architecturally. For further discussion confer Remark \ref{rem:Dissipativity}.

    For schedules with $g(0)=0$, such as the VP-SDE of Section~\ref{sec:experiments}, identity \eqref{eq:DissipativityIdentity} shows that the attainable dissipativity rate is proportional to $g^2(\tau)$, and thus vanishes towards the end of the backward process. The uniform rate $\lambda$ in Assumption~\ref{ass:Dissipativity} is then determined by the early-stopping parameter $\delta$, and the constants of Theorem~\ref{thm:MainDissipative} deteriorate as $\delta\to0$. Since the drift is strongly contractive during most of the backward process, this is likely pessimistic. A natural refinement is a time-dependent rate $\lambda(t)$ in Assumption~\ref{ass:Dissipativity}, which we leave to future work, cf.\ Remark~\ref{rem:StepwiseConstants}. For schedules with $g$ bounded away from zero, such as OU, this issue does not arise.

    \end{description}

\subsection{Main Result}
\label{subsec:MainResult}
We state two main theorems. The difference between the theorems is that one is proven in a more general setting, without the dissipative drift assumption, while the other includes this assumption. The proofs of the theorems can be found in Appendix \ref{sec:ProofOfMainTheorem}.
\begin{thm}(General case)
    \label{thm:Main}
   Let $\BackwardNormalstartingNetworkDiscretizedGeneralStrong_0,...,\BackwardNormalstartingNetworkDiscretizedGeneralStrong_K$ be the result of an \SchemeClassShortHand of order $\gamma$. Under our imposed assumptions, we get the error bound
   \[\mathcal{W}_2\left(\mathcal{L}\left(\BackwardNormalstartingNetworkDiscretizedGeneralStrong_K\right), p_0\right)=\O\left(e^{\text{Poly}(T)}(B_\text{init}(T)+h^\gamma +\varepsilon)\right).\]
\end{thm}

\begin{thm}(Dissipative case)
    \label{thm:MainDissipative}
    Let $\BackwardNormalstartingNetworkDiscretizedGeneralStrong_0,\dots,\BackwardNormalstartingNetworkDiscretizedGeneralStrong_K$ be the result of an \SchemeClassShortHand of order $\gamma$. Impose, in addition to the assumptions of Theorem~\ref{thm:Main}, Assumption~\ref{ass:Dissipativity}, and let $h\le h_0$, where $h_0=h_0(\lambda,T)$ is the step-size threshold of Lemma~\ref{lem:GeneraslStrongSchmeDissipative}, with $1/h_0\le\text{Poly}(T)$. Then
    \[
        \mathcal{W}_2\big(\mathcal{L}(\BackwardNormalstartingNetworkDiscretizedGeneralStrong_K),p_0\big)
        =\O\big(\text{Poly}(T)\,(B_\text{init}(T)+h^\gamma+\varepsilon)\big).
    \]
\end{thm}

\begin{rem}[On dissipativity]
\label{rem:Dissipativity}
Note that dissipativity is not treated in standard works such as \cite{kloeden2010numerical}. Its most striking consequence is that the $T$ dependence can be relaxed from $e^{\text{Poly}(T)}$ to a polynomial one. In particular, this allows the initialization error to be shrinking slower, effectively accommodating a broader class of data distributions.

Moreover, we want to note that Assumption~\ref{ass:Dissipativity} demands strong monotonicity uniformly on $[0,T]$. This amounts to strong log-concavity of $p_t$ for all times $t$ and therefore excludes multimodal targets. This restriction is, however, confined to a window of \emph{fixed} length: if $p_0$ is merely weakly log-concave the marginals $p_t$ become strongly log-concave for all $t \ge t_0$ with a finite, explicitly computable $t_0$, a phenomenon termed \emph{regime shift} \cite[Props.~3-4]{kremling2025weaklogconcave}, see also \cite{gentilonisilveri2025OU_Convergence}. Although \cite{kremling2025weaklogconcave} study probability flow ODEs, these two propositions concern only the forward marginals $p_t$ and thus apply unchanged in our setting. Consequently, increasing $T$ only extends the dissipative regime, on which the mechanism of Theorem~\ref{thm:MainDissipative} applies, while the non-dissipative window retains fixed length $t_0$ independent of $T$. This is consistent with our experiments in Section \ref{sec:experiments}, where the sampling error of a multimodal Problem~\ref{prob:highdim_gmm} does not deteriorate with growing $T$ despite the failure of Assumption~\ref{ass:Dissipativity}.
\end{rem}

\begin{rem}[On the polynomial $T$-dependence]
\label{rem:StepwiseConstants}
The polynomial $T$-dependence in Theorem~\ref{thm:MainDissipative} originates from bounding all coefficients uniformly on $[0,T]$. This is likely pessimistic: the coefficients are largest at the beginning of the backward process, and the errors committed there are subsequently contracted by the dissipative drift, cf.\ the factor $(1-\lambda h)$ in the recursions of Lemmas~\ref{lem:ItôTaylorDissipative} and~\ref{lem:GeneraslStrongSchmeDissipative}. A finer, step-wise tracking of the constants might therefore reduce the polynomial factor, which would be in line with the observations in Figure~\ref{fig:Dissipativity}. We leave such a refinement to future work.
\end{rem}

The asymptotic bound of Theorem~\ref{thm:MainDissipative} translates directly into an iteration-complexity statement of the kind common in the diffusion-model literature. We note that the notation $\tilde \O$ is used to suppress any logarithmic factors.
\begin{cor}[Iteration complexity, dissipative case]
\label{cor:IterationComplexityDissipative}
Under the assumptions of Theorem~\ref{thm:MainDissipative} and $B_\text{init}(T)\le C_0e^{-cT}$ for some $c>0$, the choices $T=\tfrac2c\log(1/\zeta)$, $\varepsilon\le\zeta/\text{Poly}(T)$, and $h\le\min\{h_0(T),(\zeta/\text{Poly}(T))^{1/\gamma}\}$ yield $\mathcal W_2=\O(\zeta)$ with the number of steps being
\[
    K=\tilde\O\big(\zeta^{-1/\gamma}\big).
\]
For $\gamma=1.5$ this gives $K=\tilde\O(\zeta^{-2/3})$.
\end{cor}
Table~\ref{tab:IterationComplexity} places our results next to the existing literature. Among the $\mathcal W_2$ results, the SRA family is the first stochastic sampler with a provable rate exceeding that of Euler--Maruyama.
We want to note furthermore, that in Section~\ref{sec:ExplicitSchemes}, we show that the Euler--Maruyama and exponential integrator samplers analyzed in several of these works fit into our framework.

\begin{table}[!htp]
\centering
\begin{threeparttable}
\setlength{\tabcolsep}{4pt}
\begin{tabular}{lllllll}
\toprule
Reference & Metric & Process & Sampler & $T$ & $\varepsilon_{\text{score}}$ & $K$ (OU) \\
\midrule
Chen et al. \cite{chen2023samplingeasylearningscore} & TV & OU & EI & $\log(1/\zeta)$ & $\tilde\O(\zeta)$ & $\tilde\O(\zeta^{-2})$ \\
Li et al. \cite{li2023towards} & TV/KL & VP\tnote{a}  & DDPM & - & $\tilde\O(\zeta)$ & $\tilde\O(\zeta^{-2})$\\
Li et al. \cite{li2023towards} & TV/KL & VP\tnote{a}  & accel.\ DDPM & - & - & $\tilde\O(\zeta^{-1})$ \\
Chen, Lee, Lu \cite{chen2023improved} & KL & OU  & EM & $\log(1/\zeta)$ & $\tilde\O(\zeta)$ & $\tilde\O(\zeta^{-2})$ \\
Chen, Lee, Lu \cite{chen2023improved} & KL & OU  & EI & $\log(1/\zeta)$ & $\tilde\O(\zeta)$ & $\tilde\O(\zeta^{-2})$ \\
Benton et al. \cite{benton2023nearly} & KL & OU & EI & $\log(1/\zeta)$ & $\O(\zeta)$ & $\tilde\O(\zeta^{-2})$ \\
Wu, Chen, Wei \cite{wu2024stochastic} & KL & OU & SRK (1 eval.)\tnote{b} & $\log(1/\zeta)$ & $\tilde\O(\zeta^{6})$ & $\tilde\O(\zeta^{-1})$ \\
He et al. \cite{SDE-DPM-2-paper} & KL & OU & SDE-DPM-2 & $\log(1/\zeta)$ & $\tilde \O(\zeta)$ & $\tilde\O(\zeta^{-1})$ \\
He et al. \cite{SDE-DPM-2-paper} & KL & OU & EI & $\log(1/\zeta)$ & $\tilde \O(\zeta)$ & $\tilde\O(\zeta^{-2})$ \\
He et al. \cite{SDE-DPM-2-paper} & KL & OU & RK-2 & $\log(1/\zeta)$ & $\tilde \O(\zeta)$ & $\tilde \O(\zeta^{-2})$ \\
Yu et al. \cite{yu2026advancing} & $\mathcal W_2$ & OU & EM/EI &$\log(1/\zeta)$ & $\O(\zeta)$ & $\tilde \O(\zeta^{-2})$ \\
Yu et al. \cite{yu2026advancing} & $\mathcal W_2$ & OU & REM/REI &$\log(1/\zeta)$ & $\O(\zeta)$ & $\tilde \O(\zeta^{-2})$ \\
Yu et al.\cite{yu2026advancing} & $\mathcal W_2$ & OU & 2nd order &$\log(1/\zeta)$ & $\O(\zeta)$ & $\tilde \O(\zeta^{-1})$ \\
Beyler, Bach \cite{beyler2025convergence} & $\mathcal W_2$ & VE & EM & $(\zeta^{-2/3}$)\tnote{c} & $\O(\zeta^{5/3})$ & $\O(\zeta^{-8/3})$\tnote{d} \\ 
Gao, Nguyen, Zhu \cite{gao2025wassersteinconvergenceguaranteesgeneral} & $\mathcal W_2$ & $f,g$\tnote{e} & EM & $\log(1/\zeta)$ & $\O(\zeta)$ & $\tilde \O(\zeta^{-2}) $ \\
Bruno et al. \cite{bruno2025OU_convergence_Log_Concave}\tnote{f} & $\mathcal W_2$ & OU & EM & $\log(1/\zeta)$ & $\tilde \O(\zeta)$ & $\tilde\O(\zeta^{-2})$ \\
Strasman et al. \cite{strasman2025NonAsymptoticBounds} & $\mathcal W_2$ & VP & EI & $\log(1/\zeta)$ & $\tilde \O(\zeta) $ & $\tilde \O (\zeta^{-2})$ \\
Gentiloni-Silveri, Ocello \cite{gentilonisilveri2025OU_Convergence} & $\mathcal W_2$  & OU & EM & $\log(1/\zeta)$ & $\tilde \O(\zeta) $ & $\tilde \O (\zeta^{-2})$ \\
Wang, Wang \cite{wang2026linearorderEM} &$\mathcal{W}_2$&  OU & EM & $\log(1/\zeta)$ & $\tilde \O(\zeta)$ & $\tilde O(\zeta^{-1})$   \\
\midrule
Ours\tnote{g}\; $\gamma=1$ & $\mathcal W_2$ & $f,g$\tnote{e} & EM/EI & $\log(1/\zeta)$ & $\tilde\O(\zeta)$ & $\tilde \O(\zeta^{-1})$ \\
Ours\tnote{g}\; $\gamma=1.5$ & $\mathcal W_2$ & $f,g$\tnote{e} & explicit SRA family & $\log(1/\zeta)$ & $\tilde \O(\zeta)$ & $\tilde \O(\zeta^{-2/3})$ \\
Ours\tnote{g}\; $\gamma=$ general & $\mathcal W_2$ & $f,g$\tnote{e} & general & $\log(1/\zeta)$ & $\tilde \O(\zeta)$ & $\tilde \O(\zeta^{-1/\gamma})$ \\
\bottomrule
\end{tabular}
\begin{tablenotes}[flushleft]
\footnotesize
\item[a] Time-discrete case only. 
\item[b] Single score evaluation per step. Requires bounded support and early stopping (the stated $K$ suppresses a factor $\delta^{-1}$). The score-error dependency is acknowledged as suboptimal by the authors.
\item[c] Polynomial rather than logarithmic $T$: bounded support yields only $T^{-3/2}$-decay of the propagated initialization error.
\item[d] Derived from their proposition (variant without early stopping) by parameter balancing. Not stated in this form by the authors. 
\item[e] $f$ and $g$ are meant to be fully general, as defined in Section \ref{subsec:ProblemSetup}. 
\item[f] Also convergence for the Milstein scheme and $K=\tilde \O(\zeta^{-1})$. However, in the setting considered, the Milstein scheme coincides with the EM scheme. 
\item[g] Under Assumption~\ref{ass:Dissipativity}, which for the true score amounts to strong log-concavity of every $p_\tau$ with constant above $f(\tau)/g^2(\tau)$, cf.\ Corollary~\ref{cor:IterationComplexityDissipative}. Theorem~\ref{thm:Main} alone does not yield a polynomial iteration complexity.
\end{tablenotes}

\end{threeparttable}
\caption{Iteration complexity of \emph{stochastic} (SDE-based) samplers to reach error $\zeta$ in the stated metric, specialized to the OU forward process where the original result is more general. $h$ is determined by $h=T/K$. Results are stated for $\mathrm{KL}\le\zeta^2$, TV $\le\zeta$, and $\mathcal W_2\le\zeta$. The metrics are not mutually comparable without further assumptions, and the score matching error depends on how it is defined in each work, in most cases along the trajectory of the iterates. The results further rely on different assumptions on the data distribution, ranging from finite second moments for most KL/TV results to log-concavity or bounded support for most $\mathcal W_2$ results. Note the $\tilde \O$-notation above Corollary \ref{cor:IterationComplexityDissipative}.} 

\label{tab:IterationComplexity}
\end{table}

\begin{rem}[Towards non-asymptotic bounds]
\label{rem:NonAsymptotic}
Although Theorems~\ref{thm:Main} and \ref{thm:MainDissipative} are stated asymptotically, essentially all of the underlying machinery is already fully explicit. It is therefore possible by our techniques to derive explicit non-asymptotic bounds, cf. the proof in Appendix \ref{sec:proof} below.
\end{rem}

\section{Explicit Schemes}
\label{sec:ExplicitSchemes}
In this section, we determine the order of convergence of several schemes by applying Theorem~\ref{thm:Main}. In view of the theorem, this essentially reduces to verifying the conditions of Definition~\ref{def:generalstrongscheme} for the schemes under consideration. To this end, we compute the coefficient functions associated with SDEs of the form \eqref{eq:BackwardNormalstartingScoreSDE}, that is, SDEs with additive noise, in particular with a general vector-valued drift coefficient $\mu(t,x_t)\in\R^d$ and a spatially independent, scalar diffusion coefficient $g(t)$:
\[
\mathrm dx_t = \mu(t,x_t)\,\mathrm dt + g(t)\,\mathrm dW_t, \qquad x_t\in\R^d,
\]
driven by a $d$-dimensional Wiener process $W_t$, so that the diffusion matrix is $g(t)I_d$. The proofs corresponding to the schemes can be found in Appendix \ref{sec:ProofsExplicitSchemes}. For the introduction of our schemes we need to introduce three stochastic integrals:
\begin{subequations}
\label{eq:StochasticIntegrals}
\begin{align}
    h &:= I_{(0)}[1]_{t_k,t_{k+1}} = \int_{t_k}^{t_{k+1}}\d t\quad \Rightarrow \quad I_{(0,0)}[1]_{t_k, t_{k+1}}=\int _{t_k}^{t_{k+1}}\int_{t_k}^t\d s\d t = \tfrac{h^2}{2},\\
    \Delta W_{t_k} &:= \sum_{j=1}^d I_{(j)}[1]_{t_k,t_{k+1}}\,e_j = W_{t_{k+1}}-W_{t_k} \;\in\R^d,
    \qquad \Delta W_{t_k}^{(j)} := I_{(j)}[1]_{t_k,t_{k+1}},\\
    \Delta Z_{t_k} &:= \sum_{j=1}^d I_{(j,0)}[1]_{t_k,t_{k+1}}\,e_j = \int_{t_k}^{t_{k+1}}\bigl(W_t-W_{t_k}\bigr)\,\d t \;\in\R^d,
    \qquad \Delta Z_{t_k}^{(j)} := I_{(j,0)}[1]_{t_k,t_{k+1}}.
\end{align}    
\end{subequations}

An immediate consequence of their definition is that for every $k$ the pair $(\Delta W_{t_k},\Delta Z_{t_k})$ is $\mathcal F_{t_{k+1}}$-measurable and independent of $\mathcal F_{t_k}$. In Appendix \ref{sec:ProofsExplicitSchemes}, we discuss how to simulate those stochastic integrals in practice. We want to note that simulating higher stochastic integrals is possible, but requires significantly higher computational cost, especially in high dimensions, cf.\ \cite[Chapter 5.8 or Exercise 10.5.1]{kloeden2010numerical}, which is one reason why we restrict ourselves to strong order $\gamma=1.5$.
\subsection{Euler--Maruyama}
The first scheme we consider is the Euler--Maruyama (EM) scheme
\begin{equation}
        \label{eq:EMDiscretization}
        x_{k+1}=x_k + h\,\mu(t_k,x_k) + g(t_k)\Delta W_{t_k}.
\end{equation}
Even though results on this scheme exist already (cf. Table \ref{tab:IterationComplexity}), we show its order of convergence to reiterate the generality of our main theorem.   
\begin{thm}
    \label{lem:EMScheme}
    In our setting, EM is a strong Itô--Taylor scheme, and thus an \SchemeClassShortHand with zero remainder $R_k$ of order $\gamma=1$.
\end{thm}
\subsection{The Exponential Integrator Scheme}
Unlike the Euler--Maruyama scheme above, which applies to any SDE of the general form introduced at the start of this section, the exponential integrator we investigate next exploits the specific linear-plus-score structure $\mu(t,x)=f(t)x+g^2(t)s_\theta(t,x)$ of the drift in \eqref{eq:BackwardNormalstartingScoreSDE}, together with the reverse-time reparametrization $\tau=T-t$ from Section~\ref{sec:main}. We therefore return to the notation $\BackwardNormalstartingNetworkDiscretized_k$ for the discretized process, matching its use throughout the rest of this paper, rather than the generic $x_k$ used above. The difference to the Euler--Maruyama method is that the exponential integrator only discretizes the unknown score-function part. 

Specifically, the exponential integrator is the solution of the frozen SDE
\[
    \d\BackwardNormalstartingNetworkDiscretized_t = \big[f(\tau)\BackwardNormalstartingNetworkDiscretized_t + g^2(\tau)s_\theta(\tau_k,\BackwardNormalstartingNetworkDiscretized_{t_k})\big]\d t + g(\tau)\,\d W_t, \qquad t\in[t_k,t_{k+1}].
\]
The solution to this SDE can be explicitly computed as
\begin{align*}
    \BackwardNormalstartingNetworkDiscretized_{k+1}
    &= \Phi_{t_{k+1},t_k}\,\BackwardNormalstartingNetworkDiscretized_k
      + s_\theta(\tau_k,\BackwardNormalstartingNetworkDiscretized_k)\int_{t_k}^{t_{k+1}}\Phi_{t_{k+1},s}\,g^2(T-s)\,\d s
      + \int_{t_k}^{t_{k+1}}\Phi_{t_{k+1},s}\,g(T-s)\,\d W_s,
\end{align*}
where we define the integrating factor in $t\in[t_k,t_{k+1}]$ as
\[
    \Phi_{t,t_k} := \exp\!\left(\int_{t_k}^t f(T-u)\,\d u\right), \qquad \Phi_{t_k,t_k}=1.
\]
We abbreviate $\Phi_k:=\Phi_{t_{k+1},t_k}$ and write the scheme as

\begin{subequations}\label{eq:EI_compact}
\begin{align}
&a_k := \int_{t_k}^{t_{k+1}}\Phi_{t_{k+1},s}\,g^2(T-s)\d s, \\
&b_k := \int_{t_k}^{t_{k+1}}\Phi_{t_{k+1},s}\,g(T-s)\,\d W_s  \\
&\BackwardNormalstartingNetworkDiscretized_{k+1} = \Phi_k\BackwardNormalstartingNetworkDiscretized_k + a_k\,s_\theta(\tau_k,\BackwardNormalstartingNetworkDiscretized_k) + b_k.
\end{align}
\end{subequations}
We note that $ b_k \sim \mathcal N(0,\sigma_k^2I_d)$ with
\[
    \sigma_k^2 := \int_{t_k}^{t_{k+1}}\Phi_{t_{k+1},s}^2\,g^2(T-s)\,\d s.
\]
However, for its analysis in Section \ref{sec:proof}, we need additional structural assumptions not yet stated in Section \ref{sec:assumptions}, as they are only needed in this specific case.
\begin{assumption}(regularity for the exponential integrator)
\label{ass:EI}
The function $f$ is bounded by $F(T)\in \R$ and Lipschitz continuous with constant $L_f(T)$ on $[0,T]$, the function $g$ is Lipschitz continuous with constant $L_g(T)$ on $[0,T]$, and the score network grows at most linearly, i.e.\ there exists $C_s(T)>0$ such that
\[
    \|s_\theta(\tau,x)\|\le C_s(T)(1+\|x\|)\qquad\text{for all }\tau\in[0,T],\ x\in\R^d.
\]
All constants, i.e. $L_f(T),F(T),L_g(T),C_s(T)$, depend at most polynomial on $T$.
\end{assumption}
Together with Assumption~\ref{ass:ScoreFunctionLipschitz}, the growth condition is equivalent to $\sup_{\tau\in[0,T]}\|s_\theta(\tau,0)\|<\infty$, which holds, e.g., for networks with bounded weights and a bounded time embedding. Note that Assumption~\ref{ass:CoefficientFunctions} only yields linear growth of $g^2s_\theta$, and thus of $s_\theta$ only where $g$ is bounded away from zero. 
\begin{thm}
    \label{lem:EIScheme}
    Under Assumptions \ref{ass:Solution} to \ref{ass:DriftAndDiffusion} and Assumption \ref{ass:EI}, the exponential integrator \eqref{eq:EI_compact} is an \SchemeClassShortHand of order $\gamma=1$ for all $h\le h_0^\text{EI}$, with $h_0^\text{EI}$ from Lemma~\ref{lem:MomentsEISRA}\ref{lem:MomentsEISRA_EI}.
\end{thm}

\begin{rem}[Versions of the exponential integrator]
\label{rem:EIParameterization}
In practice, exponential integrators for diffusion models are typically applied to a noise- or data-prediction parameterization of the network, whose frozen output varies more slowly than the score \cite{lu2022dpmsolver,zhang2023deis,lu2025dpm}. In the data-prediction parameterization, the reverse SDE gains an additional term linear in $x$, which can be included in the exactly integrated linear part. Hence, modern implementations of the exponential integrator differ from ours, which is a likely reason for its relatively poor performance compared to EM on Problems~\ref{prob:highdim_gmm} and~\ref{prob:highdim_sg}. We note however, that this reparameterization does not change the order: the coefficient functions $f_\alpha$ depend only on the drift $\mu$, which is unchanged, and the proof of Theorem~\ref{lem:EIScheme} carries over with $f$ replaced by an $\tilde f$ that absorbs the extra linear term.
\end{rem}

\subsection{The explicit SRA family}
As with the Euler--Maruyama scheme, the explicit SRA family (cf. \cite[Section 6.5]{rossler2010runge}) we consider next applies to the general SDE with additive noise $\mathrm dx_t=\mu(t,x_t)\,\mathrm dt+g(t)\,\mathrm dW_t$ introduced at the start of this section. We therefore return to the generic notation $x_k$, $\mu$, $g$. This family of stochastic Runge--Kutta methods is a natural candidate to consider in our case, as it is specifically tailored for SDEs with additive noise, and includes members with only two stages, considerably reducing the computational cost compared to methods achieving the same order for fully general SDEs. Using the notation we introduced in this section, an $s$-stage member of the SRA family takes the form
\begin{subequations}\label{eq:SRA}
\begin{align}
&x_{k+1}=x_k + \sum_{i=1}^s\alpha_i\mu(t_k+c_i^{(0)}h, H_i)h + \sum_{i=1}^s\left(\beta_i^{(1)}\Delta W_{t_k}+\frac{\beta_i^{(2)}}{h}\Delta Z_{t_k}\right)g(t_k+c_i^{(1)}h)\label{eq:SRA_a}\\
&H_i=x_k+\sum_{j=1}^sA_{ij}^{(0)}\mu(t_k+c_j^{(0)}h,H_j)h + \sum_{j=1}^s B_{ij}^{(0)}g(t_k+c_j^{(1)}h)\frac{1}{h}\Delta Z_{t_k} \label{eq:SRA_b},
\end{align}
\end{subequations}
with parameters $c^{(0)},c^{(1)},\beta^{(1)}, \beta^{(2)}, \alpha \in \R^s$ and $A^{(0)},B^{(0)}\in \R^{s\times s}$. Furthermore, an SRA scheme is \emph{explicit} when $A^{(0)}$ and $B^{(0)}$ are strictly lower triangular. Following \cite{rossler2010runge}, the parenthesized superscripts $(0)$, $(1)$, $(2)$ label the coefficient sets and stage values of the scheme and do not denote components. Entries of the coefficients are indexed by subscripts. The associated Butcher tableau takes the form
\begin{equation}
\label{eq:ButcherTableau}
    \centering
    \renewcommand{\arraystretch}{1.3}
    \begin{tabular}{c}
        $\begin{array}[t]{c|c|c|c}
            c^{(0)} & A^{(0)} & B^{(0)} &  c^{(1)}\\[0.5ex]
            \hline\rule{0pt}{2.8ex}
            & \alpha^T & {\beta^{(1)}}^T & {\beta^{(2)}}^T
        \end{array}.$
    \end{tabular}
\end{equation}
Similarly to the EI scheme, the SRA-family needs an additional set of assumptions.
\begin{assumption}(regularity for the explicit SRA family)
\label{ass:SRA}
The drift $\mu$ is three times continuously differentiable in $(t,x)$, and all its partial derivatives up to total order three grow at most linearly in $x$, i.e.\ there exists $C_{\mathrm{SRA}}>0$ such that
\[
    \big\|\partial_t^aD^b\mu(t,x)\big\|\le C_{\mathrm{SRA}}(T)(1+\|x\|)\qquad\text{for all }a+b\le3,\ t\in[0,T],\ x\in\R^d,
\]
where $\|\cdot\|$ denotes the operator norm of the $b$-linear map $D^b$, and $C_{\mathrm{SRA}}(T)$ depends at most polynomial on $T$.
\end{assumption}
Assumption~\ref{ass:CoefficientFunctions} does not imply this, since it only controls the coefficient functions $f_\alpha$, which are sums of derivatives of $\mu$ in which second and higher spatial derivatives appear either multiplied by powers of $g$ or contracted with $\mu$. Runge--Kutta schemes, however, evaluate $\mu$ at perturbed points, and the remainders of their expansion contain the individual derivatives. For $\mu(t,x)=f(T-t)x+g^2(T-t)s_\theta(T-t,x)$, Assumption~\ref{ass:SRA} holds if $f,g\in C^3$ with derivatives bounded polynomial in $T$ and all partial derivatives of $s_\theta$ up to total order three grow at most linearly. The regularity of $g$ needed for the SRA family, $g\in C^2$ with bounded $g'$ and $g''$, already follows from Assumption~\ref{ass:CoefficientFunctions}, since $f_{(0,j)}=g'e_j$ and $f_{(0,0,j)}=g''e_j$ with $(0,j)\in\mathcal A_{3/2}$ and $(0,0,j)\in\mathcal B(\mathcal A_{3/2})$.

\begin{thm}
\label{lem:SRA}
    Under Assumptions \ref{ass:Solution} to \ref{ass:DriftAndDiffusion} and Assumption \ref{ass:SRA}, the explicit SRA family \eqref{eq:SRA} is an \SchemeClassShortHand of order $\gamma=1.5$ for all $h\le h_0^\text{SRA}$, with $h_0^\text{SRA}$ from Lemma~\ref{lem:MomentsEISRA}\ref{lem:MomentsEISRA_SRA}, provided that its coefficients satisfy the order conditions
    \begin{equation}
    \tag{OC}\label{eq:OrderConditions}
    \begin{array}{rcl@{\qquad}rcl@{\qquad}rcl}
    \alpha^Te &=& 1 & {\beta^{(1)}}^Te &=& 1 & {\beta^{(2)}}^Te &=& 0\\[1ex]
    \alpha^TB^{(0)}e &=& 1 & \alpha^T(B^{(0)}e)^2 &=& \frac{3}{2} & \alpha^TA^{(0)}e &=& \frac12\\[1ex]
    {\beta^{(1)}}^Tc^{(1)} &=& 1 & {\beta^{(2)}}^Tc^{(1)} &=& -1 & c^{(0)} &=& A^{(0)}e,
    \end{array}
    \end{equation}
    where $e=(1,\ldots,1)^T$, together with the node condition
    \begin{equation}
    \tag{N}\label{eq:SRANodes}
        c^{(0)},\,c^{(1)}\in[0,1]^s,
    \end{equation}
    and, under Assumption~\ref{ass:Dissipativity}, additionally the positivity condition
    \begin{equation}
    \tag{P}\label{eq:SRAPositivity}
        \alpha_i\ge0,\qquad i=1,\ldots,s.
    \end{equation}
\end{thm}
Conditions \eqref{eq:SRANodes} and \eqref{eq:SRAPositivity} are not part of R\"o\ss{}ler's order conditions \cite{rossler2010runge}. We need them only for the moment bound of Lemma~\ref{lem:MomentsEISRA}. Condition \eqref{eq:SRANodes} keeps all stage times inside $[t_k,t_{k+1}]$, and \eqref{eq:SRAPositivity} lets the dissipativity of the drift pass through the convex combination of the stages.

This theorem is particularly powerful because it shows that every SRA-type method is an \SchemeClassShortHand, so that our main theorem applies to all of them at once. Within the \SchemeClassShortHand, the SRA-type methods form a much more specific subclass, as each of them follows the fixed template of a Butcher tableau of the form \eqref{eq:ButcherTableau}. Nevertheless, this subclass remains broad: it contains not only the three SRA schemes introduced by R\"o\ss{}ler, but every method whose tableau satisfies the order conditions \eqref{eq:OrderConditions}, \eqref{eq:SRANodes} and \eqref{eq:SRAPositivity}. This yields a large design space in which new methods can be constructed simply by choosing suitable coefficients, each of them automatically covered by our theory. The subclass also includes existing methods beyond R\"o\ss{}ler's schemes, such as a method adapted from Chang~\cite{chang1987numerical} (see Section~\ref{sec:experiments}). Even the EM method is of SRA type, with the tableau
\begin{equation*}
\centering
\renewcommand{\arraystretch}{1.3}
\begin{tabular}{c}
    $\begin{array}[t]{c|c|c|c}
        0 & & &  0\\
        \hline
        & 1 &1 & 0
    \end{array}$
\end{tabular}
\end{equation*}
which satisfies the order conditions for strong order $\gamma = 1$ (cf.~\cite[Thm.~6.6]{rossler2010runge}). For EM, however, we retain the direct proof, as it better highlights the method's structural simplicity.

\begin{cor}
    \label{cor:ExplicitSchemesSpecialization}
    Specializing $\mu(t,x):=f(t)x+g^2(t)s_\theta(t,x)$ and $x_k:=\BackwardNormalstartingNetworkDiscretized_k$ in Theorems~\ref{lem:EMScheme} and~\ref{lem:SRA}, the Euler--Maruyama and Runge--Kutta-type schemes, applied to \eqref{eq:BackwardNormalstartingScoreSDE}, are \SchemeClassShortHand in the sense of Definition~\ref{def:generalstrongscheme} of order $\gamma=1$ and $\gamma=1.5$ respectively. Together with Theorem~\ref{lem:EIScheme}, which is already stated directly for \eqref{eq:BackwardNormalstartingScoreSDE}, Theorem~\ref{thm:Main} applies to all schemes, yielding the stated convergence guarantees.
\end{cor}

\section{Experimental Results}
\label{sec:experiments}
In this section, we compare the practical performance of the discretization methods discussed in Section \ref{sec:ExplicitSchemes} in numerical examples. Specifically, we investigate the Euler--Maruyama (EM) method \eqref{eq:EMDiscretization}, the exponential integrator (EI) method \eqref{eq:EI_compact}, and the following two members of the explicit SRA family:
\begin{equation*}
    \centering
    \renewcommand{\arraystretch}{1.3}
    \begin{tabular}{c@{\hspace{3em}}c@{\hspace{3em}}c}
        SRA1 & Chang \\[1ex]
        $\begin{array}[t]{c|cc|cc|cc}
            0 &          &          &               &               & \multicolumn{2}{l}{1} \\[0.5ex]
            \frac34 & \frac34  &          & \frac32       &               & \multicolumn{2}{l}{0} \\[0.5ex]
            \hline\rule{0pt}{2.8ex}
                      & \frac13 & \frac23& 1 & 0 & -1 & 1
        \end{array}$
        &
        $\begin{array}[t]{c|ccc|ccc|ccc}
            0 &          &          &          &               &               &               & \multicolumn{3}{l}{1} \\[0.5ex]
            \frac12 &  \frac12  &          &          &         &               &               & \multicolumn{3}{l}{0} \\[0.5ex]
            \frac12 & \frac12   &    &          &         &  \frac32       &               & \multicolumn{3}{l}{0} \\[0.5ex]
            \hline\rule{0pt}{2.8ex}
                      & 0 & \frac13 & \frac23 & 1 & 0 & 0 & -1 & 1 & 0
        \end{array}$
    \end{tabular}
\end{equation*}
The first method has two stages and was introduced by R\"o\ss{}ler in \cite[Section 6.5]{rossler2010runge}, whereas the second method has three stages and is an adjustment of a method due to Chang \cite{chang1987numerical} to SDEs having time dependent diffusion coefficients, rather than constant ones. Note that the order conditions \eqref{eq:OrderConditions}, \eqref{eq:SRANodes} and \eqref{eq:SRAPositivity} are satisfied in both cases.

\subsection{Considered Problems}
\label{subsec:Problems}
We consider six distributions as examples:
\begin{prob}
    \label{prob:highdim_gmm} $\mathcal{P}_0$ is a 3072-dimensional Gaussian
    Mixture with four equally weighted components with pairwise distinct means $\neq 0_{3072}$ and identical covariance $\Sigma = 2 I_{3072}$.
\end{prob}
\begin{prob}
\label{prob:highdim_sg}
    $\mathcal{P}_0$ is a 3072-dimensional Gaussian distribution with mean $\neq 0_{3072}$ and covariance $\Sigma = \tfrac12I_{3072}$.
\end{prob}
\begin{prob}
\label{prob:CIFAR10}
    $\mathcal{P}_0$ is the distribution of the CIFAR-10 dataset
    \cite{krizhevsky2009CIFAR10}.
\end{prob}
\begin{prob}
\label{prob:LatentCIFAR10}
    $\mathcal{P}_0$ is the distribution of the CIFAR-10 dataset in a latent
    space.
\end{prob}
\begin{prob}
\label{prob:ImageNet64}
    $\mathcal{P}_0$ is the distribution of the ImageNet dataset \cite{deng2009imagenet} scaled to a resolution of $64\times 64$.
\end{prob}

\begin{prob}
\label{prob:ImageNet512}
    $\mathcal{P}_0$ is the distribution of the ImageNet dataset \cite{deng2009imagenet} at resolution $512\times 512$, in the $64\times 64\times 4$-dimensional latent space of the Stable Diffusion autoencoder \cite{rombach2022high}, i.e.\ the distribution of the encoded images.
\end{prob}
These problems are chosen to cover complementary aspects of our theory. Problem \ref{prob:CIFAR10} serves as a classical real-world benchmark in image generation \cite{krizhevsky2009CIFAR10}, and we choose the dimensionality of Problems \ref{prob:highdim_gmm} and \ref{prob:highdim_sg} to match it. While Problem \ref{prob:highdim_sg} has a strongly log-concave distribution at all times, Problem \ref{prob:highdim_gmm} does not have this property, and thus allows us to verify the predictions made in Remark \ref{rem:Dissipativity} for non-strongly log-concave distributions. We calculate the score function and prove the claims about the two problems in the following subsection. Problem \ref{prob:LatentCIFAR10} is meant to illustrate the effect of different discretization schemes in the context of latent GDMs, such as Stable Diffusion \cite{rombach2022high}. As encoder, we use the $2048$-dimensional latent representations obtained from the penultimate layer of the Inception v3 model provided by the \texttt{pytorch\_fid} repository. We emphasize that Problem \ref{prob:LatentCIFAR10} serves as a proof of concept: we do not generate images, only features in a latent space, since no public decoder exists. Moreover, this choice of latent space allows for an interesting experiment: employing the FID metric on images generated in pixel space, like we do in Problem \ref{prob:CIFAR10}, measures an approximated squared $\mathcal{W}_2$ distance in a latent space (cf. Section \ref{subsec:MeasureingW2}), thus making our predictions about convergence not applicable. By generating directly in this latent space, we eliminate the generation and measure space mismatch, which means our predictions should again hold. Moreover, unlike for Problems \ref{prob:CIFAR10} and \ref{prob:LatentCIFAR10}, the score functions of Problems \ref{prob:highdim_gmm} and \ref{prob:highdim_sg} are available in closed form, which allows us to conduct experiments using the true score function in generation, i.e. experiments in the absence of a score matching error. To complete the set of problems, we consider the schemes in a further real-world benchmark, where we generate images of resolution $64\times 64$ in Problem \ref{prob:ImageNet64}, and latent representations of images of resolution $512\times 512$ in Problem \ref{prob:ImageNet512}. For Problem \ref{prob:ImageNet512}, we measure $\mathcal{W}_2$ directly between these latents, i.e.\ in the space in which the samples are generated, whereas the FID is computed on the decoded images.

For Problems \ref{prob:highdim_gmm} to \ref{prob:LatentCIFAR10}, we use the original score\_sde codebase of \cite{song2020score} (cf. Section "Code Availability"), and consider the forward SDE
\[
    \d \Forward_t = -\tfrac{1}{2}t \Forward_t \d t + \sqrt{t}\,\d W_t,
    \quad \text{i.e.} \quad f(t)=\tfrac{1}{2}t \quad \text{and} \quad
    g(t)=\sqrt{t},
\]
which is the standard VP-SDE introduced in \cite{song2020score}, cf. Table \ref{tab:PopularChoices}. The parameter functions of this SDE fulfill Assumption \ref{ass:DriftAndDiffusion}, grow polynomial, and all derivatives of $g$ are bounded under early stopping. Unless stated otherwise, we choose $T=4$ as the terminal time and $\delta=10^{-3}$ as the early stopping parameter. The score network is trained following the procedure described in \cite{song2020score}, employing the \emph{ncsn++} network architecture introduced therein. We choose this architecture to stay close to how GDMs are used in the real world, even though it is unclear whether this network fulfills any of Assumptions \ref{ass:CoefficientFunctions} or \ref{ass:ScoreFunctionLipschitz}. To ensure comparability, we use the same model architecture and training procedure across Problems \ref{prob:highdim_gmm} to \ref{prob:LatentCIFAR10}, even though the \emph{ncsn++} architecture, a U-Net-like convolutional neural network, is specifically designed for images in pixel space, i.e., for Problem \ref{prob:CIFAR10}. As a result, performance on the other considered problems could likely be improved by employing network architectures tailored to their respective data distributions.

For Problems \ref{prob:ImageNet64} and \ref{prob:ImageNet512} we use the EDM2 codebase as it is, without further adapting it to our setting. EDM2 identifies time with the noise level, $t=\sigma$, and uses the forward SDE
\begin{equation}\label{eq:fwd-edm}
    \d \Forward_t = \sqrt{2t}\,\d W_t, \qquad t\in[0,T],
    \quad \text{i.e.} \quad f(t)=0 \quad \text{and} \quad
    g(t)=\sqrt{2t},
\end{equation}
so that $\Forward_t = \Forward_0 + t Z$ with $Z\sim\mathcal N(0,I)$, with $T=80$ and early stopping at $\delta=2\times 10^{-3}$. Note, that the choices of $f\equiv 0$ and $g=\sqrt{2t}$ also fulfill Assumption \ref{ass:DriftAndDiffusion}, grow polynomial, and all derivatives of $g$ are bounded under early stopping. As discussed below Assumption~\ref{ass:data}, this VE setting does, however, not satisfy Assumption~\ref{ass:data}, so these experiments test the convergence order of the discretization error only. Generation starts from $\mathcal N(0,T^2 I)$ and runs through the grid
\[
    t_k = \big(T^{1/7} + \tfrac{k}{K-1}(\delta^{1/7}-T^{1/7})\big)^{7},
    \qquad k=0,\dots,K-1,
\]
the default schedule of EDM2, followed by a final denoising step from $\delta$ to $0$. Each scheme of Section~\ref{sec:ExplicitSchemes} is applied to the reverse SDE of \eqref{eq:fwd-edm} on this grid, i.e.\ with the non-uniform step sizes $h_k = t_k-t_{k+1}$. Our convergence results remain valid on such non-uniform grids when the step size $h$ is replaced by the maximal step size, so we report convergence with respect to
\[
    h_{\max} = \max_{k} h_k = t_0-t_1 .
\]
The iteration complexity of Corollary~\ref{cor:IterationComplexityDissipative}, however, assumes an equidistant grid and does not apply here. We note that the grid is equidistant in $t^{1/7}$, with step size $h = (T^{1/7}-\delta^{1/7})/(K-1)$, and that $h_{\max}\approx 7\,T^{6/7}h\approx 299\,h$ as $K$ grows. Convergence orders in $h_{\max}$ and in $h$ therefore coincide. As networks we use the pretrained class-conditional EDM2 networks \texttt{edm2-img64-s} and \texttt{edm2-img512-s} \cite{Karras2024edm2}, without guidance.

\subsection{Verifying properties of Problems \ref{prob:highdim_gmm} and \ref{prob:highdim_sg}}
\label{subsec:VerifyingAssumptions}
We begin by computing the score function, which is possible for arbitrary coefficients $f$ and $g$ satisfying Assumption~\ref{ass:DriftAndDiffusion}: Let $\mathcal{P}_0 = \sum_{k=1}^m \xi_k \mathcal{N}\left(\mu_k, \Sigma_k\right)$ be a Gaussian mixture with a probability vector $\xi \in [0,1]^m$, and define $\Sigma_{k,t} = \phi\left(t\right)^2 \Sigma_k + \varphi\left(t\right) I$, with $\phi$ and $\varphi$ from \eqref{eq:StationarySolutionForward}, and $\pi_{k,t}$ as the density of $\mathcal{N}\left(\phi\left(t\right) \mu_k, \Sigma_{k,t} \right)$. Then, we have $p_t\left(x\right) = \sum_{k=1}^m \xi_k \pi_{k,t}\left(x\right)$, and the score function is
\begin{align}
    \label{eq:TrueScoreFunction}
    s(t,x)=\nabla \log p_t\left(x\right) = - \sum_{k=1}^m w_k\left(x\right)
    \Sigma_{k,t}^{-1} \left(x - \phi\left(t\right) \mu_k\right),
\end{align}
with weights $w_k\left(x\right) = \tfrac{\xi_k \pi_{k,t}\left(x\right)}{p_t\left(x\right)}$, where $t$ denotes the forward time. In the experiments in which we generate samples with this exact score, we have $s_\theta = s$, so that Assumptions~\ref{ass:ScoreFunctionLipschitz} and~\ref{ass:Dissipativity} become verifiable properties of \eqref{eq:TrueScoreFunction}, cf.\ the discussion of Assumption~\ref{ass:Dissipativity}. Concerning Assumption~\ref{ass:ScoreFunctionLipschitz}, the score \eqref{eq:TrueScoreFunction} is Lipschitz uniformly on $[0,T]$.

For Assumption~\ref{ass:Dissipativity}, we now specialize to our choice of $f$ and $g$, for which the explicit forms of $\phi$ and $\varphi$ are
\begin{equation*}
         \phi(t)=e^{-t^2/4} \quad \text{and}\quad \varphi(t)=1-e^{-t^2/2}.
\end{equation*}
Recall from identity \eqref{eq:DissipativityIdentity} that the drift of \eqref{eq:BackwardNormalstartingScoreSDE} is dissipative with rate $\lambda$ at backward time $t$ if and only if the forward marginal $p_\tau$ at $\tau = T-t$ is strongly log-concave with constant at least $\big(\lambda + f(\tau)\big)/g^2(\tau)$. For our VP-SDE, the threshold $f(\tau)/g^2(\tau)$ equals $\tfrac12$. For Problem \ref{prob:highdim_sg}, where $\Sigma_1 = \tfrac12 I_d$, we have $p_\tau = \mathcal{N}\big(\phi(\tau)\mu_1,\, \Sigma_{1,\tau}\big)$ with
\[
    \Sigma_{1,\tau} = v(\tau)\, I_d, \qquad
    v(\tau) = \tfrac12\phi(\tau)^2 + \varphi(\tau)
    = 1-\tfrac12 e^{-\tau^2/2} \le 1.
\]
Its log-concavity constant is thus $v(\tau)^{-1} \ge 1 > \tfrac12$, and Assumption~\ref{ass:Dissipativity} holds with $\lambda = \delta/2$ on $[0, T-\delta]$ for any $\delta > 0$. For Problem \ref{prob:highdim_gmm}, where $\Sigma_k = 2 I_d$ for all $k$, Assumption~\ref{ass:Dissipativity} fails near $t = T$. Gaussian mixtures are the prototypical example of weakly but not strongly log-concave distributions \cite[Prop.~4.1]{gentilonisilveri2025OU_Convergence}\cite[Section 3]{kremling2025weaklogconcave}. For our instance, non-log-concavity is immediate: with mode separation $\|\mu_1 - \mu_2\| = \sqrt{89}$ at component covariance $2I_d$, the density $p_\tau$ at the midpoint $\tfrac{\phi(\tau)}{2}(\mu_1+\mu_2)$ of the two scaled means is smaller than at either mean for all small $\tau$, violating quasi-concavity. Nevertheless, $\mathcal{P}_0$ is weakly log-concave, so the regime shift discussed in Remark~\ref{rem:Dissipativity} applies, and Assumption~\ref{ass:Dissipativity} is satisfied on $[0, T-\tau_1]$ for some finite $\tau_1 > 0$ independent of $T$.

\subsection{Measuring the 2-Wasserstein Distance}
\label{subsec:MeasureingW2}
To quantify the 2-Wasserstein distance between generated samples and the data distribution, we adopt a strategy similar to the widely used Fréchet Inception Distance (FID) \cite{heusel2017gans}: We approximate both the data distribution $\mathcal{P}_0$ and the generated samples by Gaussians $\mathcal{N}(\mu_d,\Sigma_d)$ and $\mathcal{N}(\mu_g,\Sigma_g)$, where $\mu_d$ and $\Sigma_d$ are the empirical mean and covariance of the data, and $\mu_g$ and $\Sigma_g$ are the empirical mean and covariance of our generated samples, respectively, and compute
\[\mathcal{W}_2\left(\mathcal{P}_0, \mathcal{L}\left(\BackwardNormalstartingNetworkDiscretized_K\right)\right)\approx\mathcal{W}_2\left(\mathcal{N}(\mu_d,\Sigma_d), \mathcal{N}(\mu_g,\Sigma_g)\right).\]
However, unlike FID, which approximates the Gaussians in the Inception v3 latent space used for Problem \ref{prob:LatentCIFAR10} and reports the squared distance $\mathcal W_2^2$, we compute the mean and variance directly in the space in which the samples are generated, i.e.\ in the respective latent space for Problems \ref{prob:LatentCIFAR10} and \ref{prob:ImageNet512}, and in pixel space otherwise, and report $\mathcal W_2$ itself. A further difference to FID is our treatment of the covariance matrix: while FID employs the full covariance matrix, we often use a diagonal approximation. This choice is motivated by a practical issue from which FID also suffers, namely instability when only a small number of samples is available. This instability is more pronounced when using the full covariance matrix and further exacerbated as the dimensionality increases. As a result, using the full covariance matrix becomes ill-suited for measuring small $\mathcal{W}_2$ distances in high dimensions. This effect is clearly illustrated in Table \ref{tab:NumberOfSamplesVsW2SG}. 

\begin{table}[!htp]
    \centering
    \begin{tabular}{|c|c|c|c|c|c|c|}
    \hline
         \# of samples & $5\cdot 10^4$ &$10^5$ & $2\cdot 10^5$ &$4\cdot 10^5$   \\
         \hline
         $\mathcal{W}_2$ using diagonal $\Sigma_g$& $0.303$ & $0.214$ & $0.149$ & $0.110$ \\
         \hline
         $\mathcal{W}_2$ using full $\Sigma_g$& $6.89$ & $4.871$ & $3.430$& $2.418$   \\
         \hline
    \end{tabular}
    \caption{Our $\mathcal{W}_2$ measure evaluated for different numbers of $\mathcal{N}(0,I_{3072})$ samples, sampled via \texttt{torch.randn}, and the analytic distribution $\mathcal{N}(0,I_{3072})$, i.e. $\mathcal{W}_2 \left(\mathcal{N}(0,I_{3072}), \mathcal{N}( \mu_g, \Sigma_g)\right)$, with $\mu_g$ and $\Sigma_g$ being the empirical mean and (diagonal) covariance of the samples.}
    \label{tab:NumberOfSamplesVsW2SG}
\end{table}
In this experiment, we sample from a $3072$-dimensional standard Gaussian using the \texttt{torch.randn} function and evaluate our $\mathcal{W}_2$ measure. Although the full covariance matrix offers better approximation properties in principle, these results motivate the use of a diagonal approximation for high dimensional problems. Accordingly, we use a diagonal covariance matrix for the high-dimensional settings in Problems \ref{prob:highdim_gmm} to \ref{prob:CIFAR10}, and Problems \ref{prob:ImageNet64} and \ref{prob:ImageNet512}. In contrast, for the relatively low-dimensional setting in Problem \ref{prob:LatentCIFAR10}, we employ the full covariance matrix in order to retain its superior approximation properties. Furthermore, keeping the covariance full in this case allows us to effectively recover the square root of FID values, cf. Section \ref{sec:EvaluatingCifar}. Finally, whenever the true mean and covariance of the data distribution are known, we use these quantities directly in our measure rather than relying on empirical approximations of $\mu_d$ and $\Sigma_d$. This applies to the toy problems in Problems \ref{prob:highdim_gmm} and \ref{prob:highdim_sg} as well as to the experiment reported in Table \ref{tab:NumberOfSamplesVsW2SG}.

\subsection{Evaluating Toy Problems}
In Figure \ref{fig:StepsizeVsW2TrueScore}, we consider simulations using the true score function $s$, i.e. we discretize $\BackwardNormalstarting_t$, and plot the step size $h$ against $\mathcal{W}_2$. We consider Problem \ref{prob:highdim_gmm} in the top row and Problem \ref{prob:highdim_sg} on the bottom row. Each row consists of the EM method, the EI method, Chang's method, and the SRA1 method. Note that we have access to the true mean and variance of the data in both problems, which means we do not have to use the empirical mean and variance for our $\mathcal{W}_2$ evaluation. We simulate $300\,000$ samples for each sampler with each stepsize. The green regions along the $x$ and $y$ axes indicate stability wherever possible: for large stepsizes, numerical methods for SDEs can become unstable even for simple test functions \cite[Chapter 9.8]{kloeden2010numerical}, while for small $\mathcal{W}_2$ values, instability arises from the limited number of samples affecting our $\mathcal{W}_2$ estimation. Note that these stability regions therefore naturally depend on both the specific problem and the choice of sampler. In particular, the SRA1 sampler often converges before we can identify a stability region. Additionally, we plot reference slopes for comparison and display the empirical convergence exponent $\omega_{LS}$ within the stability region, computed via least squares, i.e. 
\begin{equation*}
    \begin{pmatrix}
        b_{LS}\\
        \omega_{LS} 
    \end{pmatrix} = \begin{pmatrix}
            1 & \log(h_1)\\
            \vdots & \vdots\\
            1 & \log(h_n)
        \end{pmatrix}^\dagger \begin{pmatrix}
            \log(\mathcal{W}_2(h_1))\\
            \vdots \\
            \log(\mathcal{W}_2(h_n))
        \end{pmatrix},
\end{equation*} with $A^\dagger$ being the Moore--Penrose inverse of a matrix $A\in \R^{m\times n}$, $m,n\in \N$.

\begin{figure}[]
\centering
\begin{tabular}{cccc}
\hspace{1.2cm}{\scriptsize  EM}& \hspace{0.7cm} {\scriptsize EI}& \hspace{7mm} {\scriptsize Chang} &\hspace{9mm} {\scriptsize SRA1} \vspace{1mm}\\ 
\begin{tikzpicture}

\definecolor{darkgray176}{RGB}{176,176,176}
\definecolor{lightgray204}{RGB}{204,204,204}
\definecolor{steelblue31119180}{RGB}{31,119,180}
\definecolor{mygreen}{RGB}{31,150,31}

\begin{axis}[
ticklabel style={font=\scriptsize},
legend cell align={left},
legend style={
  fill opacity=1,
  draw opacity=1,
  text opacity=1,
  font=\scriptsize,
  at={(0.03,0.97)},
  anchor=north west,
  draw=lightgray204
},
width=0.24\textwidth,    
height=4cm, 
log basis x={10},
log basis y={10},
tick align=outside,
tick pos=left,
x grid style={darkgray176},
xmajorgrids,
xmin=0.003705672245534736, xmax=2.698565695347128,
xminorgrids,
xmode=log,
xtick style={color=black},
xtick={1e-05,0.0001,0.001,0.01,0.1,1,10},
xticklabels={
  \(\displaystyle {10^{-5}}\),
  \(\displaystyle {10^{-4}}\),
  \(\displaystyle {10^{-3}}\),
  \(\displaystyle {10^{-2}}\),
  \(\displaystyle {10^{-1}}\),
  \(\displaystyle {10^{0}}\),
  \(\displaystyle {10^{1}}\)
},
y grid style={darkgray176},
ymajorgrids,
ymin=0.05, ymax=315.072658458867,
yminorgrids,
ylabel={{\scriptsize Problem \ref{prob:highdim_gmm}}},
ymode=log,
ytick style={color=black},
ytick={0.01,0.1,1,10,100,1000,10000},
yticklabels={
  \(\displaystyle {10^{-2}}\),
  \(\displaystyle {10^{-1}}\),
  \(\displaystyle {10^{0}}\),
  \(\displaystyle {10^{1}}\),
  \(\displaystyle {10^{2}}\),
  \(\displaystyle {10^{3}}\),
  \(\displaystyle {10^{4}}\)
}
]

\addlegendimage{empty legend}
\addlegendentry{$\omega_{\text{LS}} = 1.00$}

\addplot [semithick, blue, mark=x, mark size=3, mark options={solid}]
table {%
2 84.4626879211461
0.8 9.59465224155552
0.666666666666667 7.71286668259183
0.571428571428572 6.50065147605524
0.5 5.63300259998227
0.4 4.4602470258662
0.2666666666666667 2.94893360875306
0.16 1.77048073547816
0.08 0.899954294174925
0.04 0.470831919060095
0.02 0.284302607709736
0.01 0.211296139804571
0.005 0.185378946874542
};
\label{graph:empirical}
\def\xmin{0.003705672245534736}
\def\ymin{0.05}
\def\xminshade{0.036}
\def\xmaxshade{1}
\def\yminshade{0.3}
\def\ymaxshade{12}

\def\xclipmin{0.05}
\def\xclipmax{1}
\def\yclipmin{0.1}
\def\yclipmax{100}


\draw[mygreen, line width=3pt] 
  (axis cs:\xminshade,\ymin) -- (axis cs:\xmaxshade,\ymin);

\draw[mygreen, line width=3pt] 
  (axis cs:\xmin,\yminshade) -- (axis cs:\xmin,\ymaxshade);

\begin{scope}
  \clip (axis cs:\xclipmin,\yclipmin) rectangle (axis cs:\xclipmax,\yclipmax);

  \def\yoffset{0.7}

  \addplot[
      thick,
      black,
      domain=0.001:\xclipmax,
      samples=2
  ] { \yoffset*0.20831919060095*((x)/0.04) };
  \label{graph:h1}
\end{scope}

\end{axis}
\end{tikzpicture} & 
\begin{tikzpicture}

\definecolor{darkgray176}{RGB}{176,176,176}
\definecolor{lightgray204}{RGB}{204,204,204}
\definecolor{steelblue31119180}{RGB}{31,119,180}
\definecolor{mygreen}{RGB}{31,150,31}

\begin{axis}[
ticklabel style={font=\scriptsize},
legend cell align={left},
legend style={
  fill opacity=1,
  draw opacity=1,
  text opacity=1,
  font=\scriptsize,
  at={(0.03,0.97)},
  anchor=north west,
  draw=lightgray204
},
    width=0.24\textwidth,    
    height=4cm, 
    log basis x={10},
    log basis y={10},
    tick align=outside,
    tick pos=left,
    x grid style={darkgray176},
    xmajorgrids,
    xmin=0.003705672245534736, xmax=2.698565695347128,
    xminorgrids,
    xmode=log,
    xtick style={color=black},
    xtick={1e-05,0.0001,0.001,0.01,0.1,1,10},
    xticklabels={
      \(\displaystyle {10^{-5}}\),
      \(\displaystyle {10^{-4}}\),
      \(\displaystyle {10^{-3}}\),
      \(\displaystyle {10^{-2}}\),
      \(\displaystyle {10^{-1}}\),
      \(\displaystyle {10^{0}}\),
      \(\displaystyle {10^{1}}\)
    },
    y grid style={darkgray176},
    ymajorgrids,
    ymin=0.05, ymax=315.072658458867,
    yminorgrids,
    ymode=log,
    ytick style={color=black},
    ytick={0.01,0.1,1,10,100,1000,10000},
    yticklabels={
      \(\displaystyle {10^{-2}}\),
      \(\displaystyle {10^{-1}}\),
      \(\displaystyle {10^{0}}\),
      \(\displaystyle {10^{1}}\),
      \(\displaystyle {10^{2}}\),
      \(\displaystyle {10^{3}}\),
      \(\displaystyle {10^{4}}\)
    }
]

\addlegendimage{empty legend}
\addlegendentry{$\omega_{\text{LS}} = 1.05$}

\addplot [semithick, blue, mark=x, mark size=3, mark options={solid}]
table {%
2 404.132029293101
0.8 26.8833375734652
0.666666666666667 18.4820871977506
0.571428571428571 14.4209829414057
0.5 12.0245100381644
0.4 9.08002184565481
0.266666666666667 5.6447744798412
0.16 3.21783596800583
0.08 1.55547473954518
0.04 0.778012349965509
0.02 0.416893693630361
0.01 0.256359986750001
0.005 0.201886537318003
};

\def\xmin{0.003705672245534736}
\def\ymin{0.05}
\def\xminbox{0.008}
\def\xmaxbox{0.92}
\def\yminbox{0.22}
\def\ymaxbox{28}

\def\xclipmin{0.05}
\def\xclipmax{1}
\def\yclipmin{0.01}
\def\yclipmax{100}


\draw[mygreen, line width=3pt] 
  (axis cs:\xminbox,\ymin) -- (axis cs:\xmaxbox,\ymin);

\draw[mygreen, line width=3pt] 
  (axis cs:\xmin,\yminbox) -- (axis cs:\xmin,\ymaxbox);


\begin{scope}
  \clip (axis cs:\xclipmin,\yclipmin) rectangle (axis cs:\xclipmax,\yclipmax);

  \def\yoffset{0.5}

  \addplot[
      thick,
      black,
      domain=0.001:\xclipmax,
      samples=2
  ]{ \yoffset*0.20831919060095*(x/0.04) };
\end{scope}

\end{axis}
\end{tikzpicture} & 
\begin{tikzpicture}

\definecolor{darkgray176}{RGB}{176,176,176}
\definecolor{lightgray204}{RGB}{204,204,204}
\definecolor{steelblue31119180}{RGB}{31,119,180}
\definecolor{mygreen}{RGB}{31,150,31}

\begin{axis}[
ticklabel style={font=\scriptsize},
legend cell align={left},
legend style={
  fill opacity=1,
  draw opacity=1,
  text opacity=1,
  font=\scriptsize,
  at={(0.03,0.97)},
  anchor=north west,
  draw=lightgray204
},
    width=0.24\textwidth,    
    height=4cm, 
    log basis x={10},
    log basis y={10},
    tick align=outside,
    tick pos=left,
    x grid style={darkgray176},
    xmajorgrids,
    xmin=0.003705672245534736, xmax=2.698565695347128,
    xminorgrids,
    xmode=log,
    xtick style={color=black},
    xtick={1e-05,0.0001,0.001,0.01,0.1,1,10},
    xticklabels={
      \(\displaystyle {10^{-5}}\),
      \(\displaystyle {10^{-4}}\),
      \(\displaystyle {10^{-3}}\),
      \(\displaystyle {10^{-2}}\),
      \(\displaystyle {10^{-1}}\),
      \(\displaystyle {10^{0}}\),
      \(\displaystyle {10^{1}}\)
    },
    y grid style={darkgray176},
    ymajorgrids,
    ymin=0.05, ymax=315.072658458867,
    yminorgrids,
    ymode=log,
    ytick style={color=black},
    ytick={0.01,0.1,1,10,100,1000,10000},
    yticklabels={
      \(\displaystyle {10^{-2}}\),
      \(\displaystyle {10^{-1}}\),
      \(\displaystyle {10^{0}}\),
      \(\displaystyle {10^{1}}\),
      \(\displaystyle {10^{2}}\),
      \(\displaystyle {10^{3}}\),
      \(\displaystyle {10^{4}}\)
    }
]

\addlegendimage{empty legend}
\addlegendentry{$\omega_{\text{LS}} = 1.72$}

\addplot [semithick, blue, mark=x, mark size=3, mark options={solid}]
table {%
2 176.767222391825
0.8 3.67981515597317
0.666666666666667 3.17689274511864
0.571428571428572 2.45879839852719
0.5 1.91716267860555
0.4 1.23965384160362
0.2666666666666667 0.571808869131425
0.16 0.26121047106121
0.08 0.182794878256953
0.04 0.173664384514865
0.02 0.173324625995465
0.01 0.179082412412592
0.005 0.175458543245563
};
\label{graph:empirical}

\def\xmin{0.003705672245534736}
\def\ymin{0.05}
\def\xminbox{0.12}
\def\xmaxbox{0.92}
\def\yminbox{0.2}
\def\ymaxbox{3.8}

\def\xclipmin{0.25}
\def\xclipmax{2}
\def\yclipmin{0.01}
\def\yclipmax{100}


\draw[mygreen, line width=3pt] 
  (axis cs:\xminbox,\ymin) -- (axis cs:\xmaxbox,\ymin);

\draw[mygreen, line width=3pt] 
  (axis cs:\xmin,\yminbox) -- (axis cs:\xmin,\ymaxbox);

\begin{scope}
  \clip (axis cs:\xclipmin,\yclipmin) rectangle (axis cs:\xclipmax,\yclipmax);

  \def\yoffset{1.5}

  \addplot[
      thick,
      orange,
      domain=0.001:\xclipmax,
      samples=2
  ]{ \yoffset*0.20831919060095*(x/0.4)^(1.5) };
  \label{graph:h1.5}
\end{scope}

\end{axis}
\end{tikzpicture}&
\begin{tikzpicture}

\definecolor{darkgray176}{RGB}{176,176,176}
\definecolor{lightgray204}{RGB}{204,204,204}
\definecolor{steelblue31119180}{RGB}{31,119,180}
\definecolor{mygreen}{RGB}{31,150,31}

\begin{axis}[
ticklabel style={font=\scriptsize},
legend cell align={left},
legend style={
  fill opacity=1,
  draw opacity=1,
  text opacity=1,
  font=\scriptsize,
  at={(0.03,0.97)},
  anchor=north west,
  draw=lightgray204
},
    width=0.24\textwidth,    
    height=4cm, 
    log basis x={10},
    log basis y={10},
    tick align=outside,
    tick pos=left,
    x grid style={darkgray176},
    xmajorgrids,
    xmin=0.003705672245534736, xmax=2.698565695347128,
    xminorgrids,
    xmode=log,
    xtick style={color=black},
    xtick={1e-05,0.0001,0.001,0.01,0.1,1,10},
    xticklabels={
      \(\displaystyle {10^{-5}}\),
      \(\displaystyle {10^{-4}}\),
      \(\displaystyle {10^{-3}}\),
      \(\displaystyle {10^{-2}}\),
      \(\displaystyle {10^{-1}}\),
      \(\displaystyle {10^{0}}\),
      \(\displaystyle {10^{1}}\)
    },
    y grid style={darkgray176},
    ymajorgrids,
    ymin=0.05, ymax=315.072658458867,
    yminorgrids,
    ymode=log,
    ytick style={color=black},
    ytick={0.01,0.1,1,10,100,1000,10000},
    yticklabels={
      \(\displaystyle {10^{-2}}\),
      \(\displaystyle {10^{-1}}\),
      \(\displaystyle {10^{0}}\),
      \(\displaystyle {10^{1}}\),
      \(\displaystyle {10^{2}}\),
      \(\displaystyle {10^{3}}\),
      \(\displaystyle {10^{4}}\)
    }
]

\addlegendimage{empty legend}
\addlegendentry{$\omega_{\text{LS}} = 1.50$}

\addplot [semithick, blue, mark=x, mark size=3, mark options={solid}]
table {%
2.00000   24.33362118
0.80000   0.66549002
0.66667   0.65573696
0.57143   0.54515731
0.50000   0.43272988
0.40000   0.30850727
0.26667   0.22517450
0.20000   0.19978660
0.16000   0.18410699
0.08000   0.17688785
0.04000   0.17192217
0.02000   0.17520467
0.01000   0.17571139
0.00500   0.17642277
};
\label{graph:empirical}

\def\xmin{0.003705672245534736}
\def\ymin{0.05}
\def\xminbox{0.33}
\def\xmaxbox{0.7}
\def\yminbox{0.25}
\def\ymaxbox{0.7}

\def\xclipmin{0.4}
\def\xclipmax{2}
\def\yclipmin{0.01}
\def\yclipmax{100}


\draw[mygreen, line width=3pt] 
  (axis cs:\xminbox,\ymin) -- (axis cs:\xmaxbox,\ymin);

\draw[mygreen, line width=3pt] 
  (axis cs:\xmin,\yminbox) -- (axis cs:\xmin,\ymaxbox);

\begin{scope}
  \clip (axis cs:\xclipmin,\yclipmin) rectangle (axis cs:\xclipmax,\yclipmax);

  \def\yoffset{0.5}

  \addplot[
      thick,
      orange,
      domain=0.001:\xclipmax,
      samples=2
  ]{ \yoffset*0.20831919060095*(x/0.4)^(1.5) };
  \label{graph:h1.5}
\end{scope}

\end{axis}
\end{tikzpicture} \\
\begin{tikzpicture}

\definecolor{darkgray176}{RGB}{176,176,176}
\definecolor{lightgray204}{RGB}{204,204,204}
\definecolor{steelblue31119180}{RGB}{31,119,180}
\definecolor{mygreen}{RGB}{31,150,31}

\begin{axis}[
ticklabel style={font=\scriptsize},
legend cell align={left},
legend style={
  fill opacity=1,
  draw opacity=1,
  text opacity=1,
  font=\scriptsize,
  at={(0.03,0.97)},
  anchor=north west,
  draw=lightgray204
},
width=0.24\textwidth,
height=4cm,
log basis x={10},
log basis y={10},
tick align=outside,
tick pos=left,
x grid style={darkgray176},
xlabel={{ \scriptsize h vs $\mathcal{W}_2$} },
xmajorgrids,
xmin=0.003705672245534736, xmax=2.698565695347128,
xminorgrids,
xmode=log,
xtick style={color=black},
xtick={1e-05,0.0001,0.001,0.01,0.1,1,10},
xticklabels={
  \(\displaystyle {10^{-5}}\),
  \(\displaystyle {10^{-4}}\),
  \(\displaystyle {10^{-3}}\),
  \(\displaystyle {10^{-2}}\),
  \(\displaystyle {10^{-1}}\),
  \(\displaystyle {10^{0}}\),
  \(\displaystyle {10^{1}}\)
},
y grid style={darkgray176},
ymajorgrids,
ymin=0.05, ymax=315.072658458867,
yminorgrids,
ymode=log,
ylabel={{\scriptsize Problem \ref{prob:highdim_sg}}},
ytick style={color=black},
ytick={0.001,0.01,0.1,1,10,100,1000,10000},
yticklabels={
  \(\displaystyle {10^{-3}}\),
  \(\displaystyle {10^{-2}}\),
  \(\displaystyle {10^{-1}}\),
  \(\displaystyle {10^{0}}\),
  \(\displaystyle {10^{1}}\),
  \(\displaystyle {10^{2}}\),
  \(\displaystyle {10^{3}}\),
  \(\displaystyle {10^{4}}\)
}
]
\addlegendimage{empty legend}
\addlegendentry{$\omega_{\text{LS}} = 1.03$}
\addplot [semithick, blue, mark=x, mark size=3, mark options={solid}]
table {%
2 276.704168850992
0.8 11.7932576579938
0.666666666666667 8.78594833728367
0.571428571428571 7.10200633576531
0.5 5.99656714177888
0.4 4.60745164182299
0.266666666666667 2.94694593153325
0.16 1.72323073079564
0.08 0.849824546374611
0.04 0.424826615586377
0.02 0.221436373840846
0.01 0.130888580983391
0.005 0.0983976097226812
};

\def\xmin{0.003705672245534736}
\def\ymin{0.05}
\def\xminbox{0.008}
\def\xmaxbox{1}
\def\yminbox{0.12}
\def\ymaxbox{30}

\def\xclipmin{0.05}
\def\xclipmax{1}
\def\yclipmin{0.01}
\def\yclipmax{100}


\draw[mygreen, line width=3pt] 
  (axis cs:\xminbox,\ymin) -- (axis cs:\xmaxbox,\ymin);

\draw[mygreen, line width=3pt] 
  (axis cs:\xmin,\yminbox) -- (axis cs:\xmin,\ymaxbox);

\begin{scope}
  \clip (axis cs:\xclipmin,\yclipmin) rectangle (axis cs:\xclipmax,\yclipmax);

  \def\yoffset{0.5}

  \addplot[
      thick,
      black,
      domain=0.001:\xclipmax,
      samples=2
  ]{ \yoffset*0.20831919060095*(x/0.04) };
\end{scope}

\end{axis}
\end{tikzpicture} & \begin{tikzpicture}

\definecolor{darkgray176}{RGB}{176,176,176}
\definecolor{lightgray204}{RGB}{204,204,204}
\definecolor{steelblue31119180}{RGB}{31,119,180}
\definecolor{mygreen}{RGB}{31,150,31}

\begin{axis}[
ticklabel style={font=\scriptsize},
legend cell align={left},
legend style={
  fill opacity=1,
  draw opacity=1,
  text opacity=1,
  font=\scriptsize,
  at={(0.03,0.97)},
  anchor=north west,
  draw=lightgray204
},
width=0.24\textwidth,
height=4cm,
log basis x={10},
log basis y={10},
tick align=outside,
tick pos=left,
x grid style={darkgray176},
xlabel={{ \scriptsize h vs $\mathcal{W}_2$} },
xmajorgrids,
xmin=0.003705672245534736, xmax=2.698565695347128,
xminorgrids,
xmode=log,
xtick style={color=black},
xtick={1e-05,0.0001,0.001,0.01,0.1,1,10},
xticklabels={
  \(\displaystyle {10^{-5}}\),
  \(\displaystyle {10^{-4}}\),
  \(\displaystyle {10^{-3}}\),
  \(\displaystyle {10^{-2}}\),
  \(\displaystyle {10^{-1}}\),
  \(\displaystyle {10^{0}}\),
  \(\displaystyle {10^{1}}\)
},
y grid style={darkgray176},
ymajorgrids,
ymin=0.05, ymax=315.072658458867,
yminorgrids,
ymode=log,
ytick style={color=black},
ytick={0.001,0.01,0.1,1,10,100,1000,10000},
yticklabels={
  \(\displaystyle {10^{-3}}\),
  \(\displaystyle {10^{-2}}\),
  \(\displaystyle {10^{-1}}\),
  \(\displaystyle {10^{0}}\),
  \(\displaystyle {10^{1}}\),
  \(\displaystyle {10^{2}}\),
  \(\displaystyle {10^{3}}\),
  \(\displaystyle {10^{4}}\)
}
]
\addlegendimage{empty legend}
\addlegendentry{$\omega_{\text{LS}} = 1.00$}
\addplot [semithick, blue, mark=x, mark size=3, mark options={solid}]
table {%
2 1416.47652772322
0.8 19.1393252385686
0.666666666666667 14.1748026797588
0.571428571428571 11.6484692322358
0.5 9.9052608036142
0.4 7.63251139570707
0.266666666666667 4.85013171595703
0.16 2.80349020595261
0.08 1.36533150624193
0.04 0.676456107528273
0.02 0.345026169993193
0.01 0.187632621258156
0.005 0.123179952885052
};

\def\xmin{0.003705672245534736}
\def\ymin{0.05}
\def\xminbox{0.003}
\def\xmaxbox{1}
\def\yminbox{0.1}
\def\ymaxbox{30}

\def\xclipmin{0.05}
\def\xclipmax{1}
\def\yclipmin{0.01}
\def\yclipmax{100}


\draw[mygreen, line width=3pt] 
  (axis cs:\xminbox,\ymin) -- (axis cs:\xmaxbox,\ymin);

\draw[mygreen, line width=3pt] 
  (axis cs:\xmin,\yminbox) -- (axis cs:\xmin,\ymaxbox);

\begin{scope}
  \clip (axis cs:\xclipmin,\yclipmin) rectangle (axis cs:\xclipmax,\yclipmax);

  \def\yoffset{0.5}

  \addplot[
      thick,
      black,
      domain=0.001:\xclipmax,
      samples=2
  ]{ \yoffset*0.20831919060095*(x/0.04) };
\end{scope}

\end{axis}
\end{tikzpicture}  & 
\begin{tikzpicture}

\definecolor{darkgray176}{RGB}{176,176,176}
\definecolor{lightgray204}{RGB}{204,204,204}
\definecolor{steelblue31119180}{RGB}{31,119,180}
\definecolor{mygreen}{RGB}{31,150,31}

\begin{axis}[
ticklabel style={font=\scriptsize},
legend cell align={left},
legend style={
  fill opacity=1,
  draw opacity=1,
  text opacity=1,
  font=\scriptsize,
  at={(0.03,0.97)},
  anchor=north west,
  draw=lightgray204
},
width=0.24\textwidth,
height=4cm,
log basis x={10},
log basis y={10},
tick align=outside,
tick pos=left,
x grid style={darkgray176},
xlabel={{\scriptsize h vs $\mathcal{W}_2$}},
xmajorgrids,
xmin=0.003705672245534736, xmax=2.698565695347128,
xminorgrids,
xmode=log,
xtick style={color=black},
xtick={1e-05,0.0001,0.001,0.01,0.1,1,10},
xticklabels={
  \(\displaystyle {10^{-5}}\),
  \(\displaystyle {10^{-4}}\),
  \(\displaystyle {10^{-3}}\),
  \(\displaystyle {10^{-2}}\),
  \(\displaystyle {10^{-1}}\),
  \(\displaystyle {10^{0}}\),
  \(\displaystyle {10^{1}}\)
},
y grid style={darkgray176},
ymajorgrids,
ymin=0.05, ymax=315.072658458867,
yminorgrids,
ymode=log,
ytick style={color=black},
ytick={0.001,0.01,0.1,1,10,100,1000,10000},
yticklabels={
  \(\displaystyle {10^{-3}}\),
  \(\displaystyle {10^{-2}}\),
  \(\displaystyle {10^{-1}}\),
  \(\displaystyle {10^{0}}\),
  \(\displaystyle {10^{1}}\),
  \(\displaystyle {10^{2}}\),
  \(\displaystyle {10^{3}}\),
  \(\displaystyle {10^{4}}\)
}
]

\addlegendimage{empty legend}
\addlegendentry{$\omega_{\text{LS}} = 1.78$}

\addplot [semithick, blue, mark=x, mark size=3, mark options={solid}]
table {%
2 316.605405406334
0.8 7.13847720592525
0.666666666666667 5.45060447194491
0.571428571428571 4.14055782164146
0.5 3.22469908330352
0.4 2.10899758305732
0.266666666666667 0.968589819188535
0.16 0.37205969335644
0.08 0.12989111763414
0.04 0.0918285372143328
0.02 0.0887445561394428
0.01 0.0878282896427145
0.005 0.0880968794495154
};

\def\xmin{0.003705672245534736}
\def\ymin{0.05}
\def\xminbox{0.06}
\def\xmaxbox{1.2}
\def\yminbox{0.11}
\def\ymaxbox{10.3}

\def\xclipmin{0.25}
\def\xclipmax{1.5}
\def\yclipmin{0.01}
\def\yclipmax{100}


\draw[mygreen, line width=3pt] 
  (axis cs:\xminbox,\ymin) -- (axis cs:\xmaxbox,\ymin);

\draw[mygreen, line width=3pt] 
  (axis cs:\xmin,\yminbox) -- (axis cs:\xmin,\ymaxbox);

\begin{scope}
  \clip (axis cs:\xclipmin,\yclipmin) rectangle (axis cs:\xclipmax,\yclipmax);

  \def\yoffset{1.5}

  \addplot[
      thick,
      orange,
      domain=0.001:\xclipmax,
      samples=2
  ]{ \yoffset*0.20831919060095*(x/0.4)^(1.5) };
\end{scope}

\end{axis}

\end{tikzpicture}&
\begin{tikzpicture}

\definecolor{darkgray176}{RGB}{176,176,176}
\definecolor{lightgray204}{RGB}{204,204,204}
\definecolor{steelblue31119180}{RGB}{31,119,180}
\definecolor{mygreen}{RGB}{31,150,31}

\begin{axis}[
ticklabel style={font=\scriptsize},
legend cell align={left},
legend style={
  fill opacity=1,
  draw opacity=1,
  text opacity=1,
  font=\scriptsize,
  at={(0.03,0.97)},
  anchor=north west,
  draw=lightgray204
},
width=0.24\textwidth,
height=4cm,
log basis x={10},
log basis y={10},
tick align=outside,
tick pos=left,
x grid style={darkgray176},
xlabel={{\scriptsize h vs $\mathcal{W}_2$}},
xmajorgrids,
xmin=0.003705672245534736, xmax=2.698565695347128,
xminorgrids,
xmode=log,
xtick style={color=black},
xtick={1e-05,0.0001,0.001,0.01,0.1,1,10},
xticklabels={
  \(\displaystyle {10^{-5}}\),
  \(\displaystyle {10^{-4}}\),
  \(\displaystyle {10^{-3}}\),
  \(\displaystyle {10^{-2}}\),
  \(\displaystyle {10^{-1}}\),
  \(\displaystyle {10^{0}}\),
  \(\displaystyle {10^{1}}\)
},
y grid style={darkgray176},
ymajorgrids,
ymin=0.05, ymax=315.072658458867,
yminorgrids,
ymode=log,
ytick style={color=black},
ytick={0.001,0.01,0.1,1,10,100,1000,10000},
yticklabels={
  \(\displaystyle {10^{-3}}\),
  \(\displaystyle {10^{-2}}\),
  \(\displaystyle {10^{-1}}\),
  \(\displaystyle {10^{0}}\),
  \(\displaystyle {10^{1}}\),
  \(\displaystyle {10^{2}}\),
  \(\displaystyle {10^{3}}\),
  \(\displaystyle {10^{4}}\)
}
]

\addplot [semithick, blue, mark=x, mark size=3, mark options={solid}]
table {%
2.00000   134.20702003
0.80000   2.20037332
0.66667   1.04777387
0.57143   0.60509251
0.50000   0.39840088
0.40000   0.23510195
0.26667   0.14150733
0.20000   0.11516020
0.16000   0.10273079
0.08000   0.09097478
0.04000   0.09040910
0.02000   0.08740811
0.01000   0.08763982
0.00500   0.08592142
};

\def\xmin{0.003705672245534736}
\def\ymin{0.05}
\def\xminbox{0.24}
\def\xmaxbox{1.2}
\def\yminbox{0.13}
\def\ymaxbox{5}

\def\xclipmin{0.1}
\def\xclipmax{0.5}
\def\yclipmin{0.01}
\def\yclipmax{100}


\begin{scope}
  \clip (axis cs:\xclipmin,\yclipmin) rectangle (axis cs:\xclipmax,\yclipmax);

  \def\yoffset{10}

  \addplot[
      thick,
      orange,
      domain=0.001:\xclipmax,
      samples=2
  ]{ \yoffset*0.20831919060095*(x/0.4)^(1.5) };
\end{scope}

\end{axis}

\end{tikzpicture}\\
\end{tabular}
\caption{Simulations for Problems \ref{prob:highdim_gmm} and \ref{prob:highdim_sg} using the true score function $s$, cf. \eqref{eq:TrueScoreFunction}. The top row considers Problem \ref{prob:highdim_gmm}, while the bottom row considers Problem \ref{prob:highdim_sg}. $\mathcal{W}_2$ of the considered methods vs the stepsize $h$. Depicted is the empirical curve (\ref{graph:empirical}), a reference line $\mathcal{O}(h)$ (\ref{graph:h1}), and a reference line $\mathcal{O}(h^{3/2})$ (\ref{graph:h1.5}). The green regions on the axes denote the stability region limited by the stepsize on the x-axis, and our $\mathcal{W}_2$ measure on the y-axis.}
\label{fig:StepsizeVsW2TrueScore}
\end{figure}

In Figure~\ref{fig:Dissipativity}, we investigate the influence of the terminal time $T$. For this purpose, we simulate Problems~\ref{prob:highdim_gmm} and~\ref{prob:highdim_sg} with the true score $s$ and a fixed step size, and let $T$ grow. We use $h=0.1$ for EM and EI, and $h=0.5$ for SRA1 and Chang's method, since with $h=0.1$ the higher-order methods already reach the floor of our $\mathcal W_2$ measure. These choices eliminate the score matching error and fix the step size, thus isolating the influence of $T$. Note that in order to keep $h=T/K$ fixed, we need to use more discretization points for higher $T$. Since we use the true score function, we again simulate $300\,000$ samples.
\begin{figure}[]
\centering
\begin{tabular}{cc}
\hspace{1cm}{\scriptsize  Problem \ref{prob:highdim_gmm}}& \hspace{1cm}{\scriptsize  Problem \ref{prob:highdim_sg}}\vspace{1mm}\\ 
\begin{tikzpicture}

\definecolor{darkgray176}{RGB}{176,176,176}
\definecolor{lightgray204}{RGB}{204,204,204}
\definecolor{steelblue31119180}{RGB}{31,119,180}

\begin{axis}[
ticklabel style={font=\scriptsize},
legend style={fill opacity=0.8, draw opacity=1, text opacity=1, draw=lightgray204},
width=0.32\textwidth,    
height=4cm, 
log basis x={10},
log basis y={10},
tick align=outside,
tick pos=left,
x grid style={darkgray176},
xlabel={\scriptsize T},
xmajorgrids,
xmin=0.328936063770756, xmax=24.3208358131731,
xminorgrids,
xmode=log,
xtick style={color=black},
xtick={0.01,0.1,1,10,100,1000},
xticklabels={
  \(\displaystyle {10^{-2}}\),
  \(\displaystyle {10^{-1}}\),
  \(\displaystyle {10^{0}}\),
  \(\displaystyle {10^{1}}\),
  \(\displaystyle {10^{2}}\),
  \(\displaystyle {10^{3}}\)
},
y grid style={darkgray176},
ylabel={\scriptsize $\mathcal{W}_2$},
ymajorgrids,
ymin=0.4, ymax=67.283793985104,
yminorgrids,
ymode=log,
ytick style={color=black},
ytick={0.01,0.1,1,10,100,1000},
yticklabels={
  \(\displaystyle {10^{-2}}\),
  \(\displaystyle {10^{-1}}\),
  \(\displaystyle {10^{0}}\),
  \(\displaystyle {10^{1}}\),
  \(\displaystyle {10^{2}}\),
  \(\displaystyle {10^{3}}\)
}
]
\addplot [semithick, EM_color, mark=o, mark size=2, mark options={solid}]
table {%
0.4 55.2810499252835
0.5 49.9657005621611
0.6 44.7785140823059
0.7 39.7507094420435
1 25.9897430133726
1.5 9.59420219810311
2 2.68136375029735
2.5 1.23113429245946
3 1.08615671276024
3.5 1.10211271526729
4 1.10280295392445
5 1.10978927249325
5.5 1.10796228158296
6 1.10988203763266
6.5 1.11064971676445
10 1.10556508931611
20 1.11192040743222
};
\label{graph:EM}

\addplot [semithick, EI_color, mark=triangle, mark size=2, mark options={solid}]
table {%
 0.4                     5.655969e+01
         0.5                     5.121318e+01
         0.6                     4.598404e+01
         0.7                     4.089791e+01
           1                    2.688866e+01
         1.5                    9.886090e+00
           2                    2.544121e+00
         2.5                    1.673873e+00
           3                    1.867134e+00
         3.5                    1.945814e+00
           4                   1.963064e+00
         4.5                    1.963576e+00
           5                    1.964146e+00
           6                    1.966329e+00
         6.5                    1.961297e+00
          10                   1.964074e+00
          20                   1.961049e+00
};
\label{graph:EI}

\addplot[semithick,SRA1_color, mark=square, mark size=2, mark options={solid}]
table {%
  0.5                     5.377942e+01
           1                     2.871342e+01
           2                     3.745211e+00
           3                     5.713681e-01
           5                    4.496663e-01
          10                    4.474294e-01
          20                    1.278628e+12

};
\label{graph:SRA1}

\addplot[semithick, Chang_color, mark=x, mark size=2, mark options={solid}]
table {%
 0.5                     5.376428e+01
           1                     2.874180e+01
           2                     4.275965e+00
           3                     1.834198e+00
           5                    1.932498e+00
          10                    1.919732e+00
          20                   2.886282e+12

};
\label{graph:HO}

\end{axis}

\end{tikzpicture} & 
\begin{tikzpicture}

\definecolor{darkgray176}{RGB}{176,176,176}
\definecolor{lightgray204}{RGB}{204,204,204}
\definecolor{steelblue31119180}{RGB}{31,119,180}

\begin{axis}[
ticklabel style={font=\scriptsize},
legend style={fill opacity=0.8, draw opacity=1, text opacity=1, draw=lightgray204},
width=0.32\textwidth,    
height=4cm, 
log basis x={10},
log basis y={10},
tick align=outside,
tick pos=left,
x grid style={darkgray176},
xlabel={{\scriptsize T}},
xmajorgrids,
xmin=0.322334601603383, xmax=27.2283954012651,
xminorgrids,
xmode=log,
xtick style={color=black},
xtick={0.01,0.1,1,10,100,1000},
xticklabels={
  \(\displaystyle {10^{-2}}\),
  \(\displaystyle {10^{-1}}\),
  \(\displaystyle {10^{0}}\),
  \(\displaystyle {10^{1}}\),
  \(\displaystyle {10^{2}}\),
  \(\displaystyle {10^{3}}\)
},
y grid style={darkgray176},
ymajorgrids,
ymin=0.3, ymax=67.283793985104,
yminorgrids,
ymode=log,
ytick style={color=black},
ytick={0.01,0.1,1,10,100,1000},
yticklabels={
  \(\displaystyle {10^{-2}}\),
  \(\displaystyle {10^{-1}}\),
  \(\displaystyle {10^{0}}\),
  \(\displaystyle {10^{1}}\),
  \(\displaystyle {10^{2}}\),
  \(\displaystyle {10^{3}}\)
}
]
\addplot [semithick, EM_color, mark=o, mark size=2, mark options={solid}]
table {%
0.4 49.0817780866283
0.5 43.9889305284509
0.6 38.9932357297595
0.7 34.212022669149
1 21.8786460889166
1.5 8.76532270932968
2 2.6368415623065
2.5 0.638033339154379
3 0.890500012354069
3.5 1.03610655930466
4 1.0642348967845
5 1.06915044146703
5.5 1.06955094094578
6 1.06951373914444
6.5 1.0706020257311
10 1.07167681926063
20 1.07033655086197
};

\addplot [semithick, EI_color, mark=triangle, mark size=2, mark options={solid}]
table {%
0.4                     5.114039e+01
         0.5                     4.621091e+01
         0.6                     4.130358e+01
         0.7                    3.654657e+01
           1                    2.401857e+01
         1.5                    1.034556e+01
           2                    3.968969e+00
         2.5                    1.969177e+00
           3                    1.719297e+00
         3.5                    1.715980e+00
           4                    1.717027e+00
         4.5                    1.718070e+00
           5                    1.715236e+00
           6                    1.718448e+00
         6.5                    1.717497e+00
          10                   1.716424e+00
          20                  1.716644e+00
};

\addplot[semithick, Chang_color, mark=x, mark size=2, mark options={solid}]
table {%
 0.5                     4.805804e+01
           1                     2.611547e+01
           2                     5.535651e+00
           3                     3.211790e+00
           5                    3.228417e+00
          10                    3.226247e+00
          20                   8.226936e+11

};

\addplot[semithick,SRA1_color, mark=square, mark size=2, mark options={solid}]
table {%
 0.5                     4.726413e+01
           1                     2.570303e+01
           2                     4.909533e+00
           3                     7.602886e-01
           5                    3.913492e-01
          10                    3.874062e-01
          20                    3.732625e+11

};

\end{axis}
\end{tikzpicture}  
\end{tabular}
\caption{Simulations to test the impact of dissipativity for Problems \ref{prob:highdim_gmm} and \ref{prob:highdim_sg} using the true score function $s$, cf. \eqref{eq:TrueScoreFunction}. We plot the achieved $\mathcal{W}_2$ of the EM method (\ref{graph:EM}), the EI method (\ref{graph:EI}), Chang's method (\ref{graph:HO}), and the SRA1 method (\ref{graph:SRA1}) versus the terminal time $T$.}
\label{fig:Dissipativity}
\end{figure}

In Figure \ref{fig:StepsizeVsW2Trainend}, we run simulations using a score network trained for 100 epochs, that is, we discretize $\BackwardNormalstartingNetwork_t$. Because of the high computational cost of simulating using the score model, we reduce the number of samples generated to $50\,000$. The layout of the figure is identical to Figure \ref{fig:StepsizeVsW2TrueScore}.

\begin{figure}[]
    \centering
    \begin{tabular}{cccc}
    \hspace{1.2cm}{\scriptsize  EM}& \hspace{0.7cm} {\scriptsize EI}& \hspace{7mm} {\scriptsize Chang} &\hspace{9mm} {\scriptsize SRA1} \vspace{1mm}\\ 
\begin{tikzpicture}
\definecolor{crimson2143940}{RGB}{214,39,40}
\definecolor{darkgray176}{RGB}{176,176,176}
\definecolor{darkorange25512714}{RGB}{255,127,14}
\definecolor{forestgreen4416044}{RGB}{44,160,44}
\definecolor{lightgray204}{RGB}{204,204,204}
\definecolor{mediumpurple148103189}{RGB}{148,103,189}
\definecolor{steelblue31119180}{RGB}{31,119,180}
\definecolor{mygreen}{RGB}{31,150,31}

\begin{axis}[
ticklabel style={font=\scriptsize},
legend cell align={left},
legend style={
  fill opacity=0.8,
  draw opacity=1,
  text opacity=1,
  font=\scriptsize,
  at={(0.03,0.97)},
  anchor=north west,
  draw=lightgray204
},
width=0.24\textwidth,    
height=4cm, 
log basis x={10},
log basis y={10},
tick align=outside,
tick pos=left,
x grid style={darkgray176},
xmajorgrids,
xmin=0.00767270499010924, xmax=2.698565695347128,
xminorgrids,
xmode=log,
xtick style={color=black},
xtick={0.0001,0.001,0.01,0.1,1,10},
xticklabels={
  \(\displaystyle {10^{-4}}\),
  \(\displaystyle {10^{-3}}\),
  \(\displaystyle {10^{-2}}\),
  \(\displaystyle {10^{-1}}\),
  \(\displaystyle {10^{0}}\),
  \(\displaystyle {10^{1}}\)
},
y grid style={darkgray176},
ymajorgrids,
ymin=0.15, ymax=551.341520850039,
ylabel={{\scriptsize Problem \ref{prob:highdim_gmm}}},
yminorgrids,
ymode=log,
ytick style={color=black},
ytick={1e-06,1e-05,0.0001,0.001,0.01,0.1,1,10,100,1000,10000},
yticklabels={
  \(\displaystyle {10^{-6}}\),
  \(\displaystyle {10^{-5}}\),
  \(\displaystyle {10^{-4}}\),
  \(\displaystyle {10^{-3}}\),
  \(\displaystyle {10^{-2}}\),
  \(\displaystyle {10^{-1}}\),
  \(\displaystyle {10^{0}}\),
  \(\displaystyle {10^{1}}\),
  \(\displaystyle {10^{2}}\),
  \(\displaystyle {10^{3}}\),
  \(\displaystyle {10^{4}}\)
}
]
\addlegendimage{empty legend}
\addlegendentry{$\omega_{\text{LS}} = 1.06$}

\addplot [semithick, blue, mark=x, mark size=3, mark options={solid}]
table {%
2 233.639067073599
0.8 12.0081570054919
0.666666666666667 9.12516948097949
0.571428571428572 7.55184621638932
0.5 6.38179802727837
0.4 4.92550619746985
0.2666666666666667 3.15517633326306
0.16 1.899422274416
0.08 1.01938340581972
0.04 0.632855095472424
0.02 0.496073377717046
0.01 0.449946867642707
};
\label{graph:empirical}

\def\xmin{0.00767270499010924}
\def\ymin{0.15}
\def\xminbox{0.062}
\def\xmaxbox{1.4}
\def\yminbox{1}
\def\ymaxbox{15}

\def\xclipmin{0.07}
\def\xclipmax{1}
\def\yclipmin{0.01}
\def\yclipmax{100}


\draw[mygreen, line width=3pt] 
  (axis cs:\xminbox,\ymin) -- (axis cs:\xmaxbox,\ymin);

\draw[mygreen, line width=3pt] 
  (axis cs:\xmin,\yminbox) -- (axis cs:\xmin,\ymaxbox);

\begin{scope}
  \clip (axis cs:\xclipmin,\yclipmin) rectangle (axis cs:\xclipmax,\yclipmax);

  \def\yoffset{1.5}

  \addplot[
      thick,
      black,
      domain=0.001:\xclipmax,
      samples=2
  ]{ \yoffset*0.20831919060095*(x/0.06)^(1) };
  \label{graph:h1}
\end{scope}

\end{axis}

\end{tikzpicture} & 
\begin{tikzpicture}

\definecolor{darkgray176}{RGB}{176,176,176}
\definecolor{lightgray204}{RGB}{204,204,204}
\definecolor{steelblue31119180}{RGB}{31,119,180}
\definecolor{mygreen}{RGB}{31,150,31}

\begin{axis}[
ticklabel style={font=\scriptsize},
legend cell align={left},
legend style={
  fill opacity=1,
  draw opacity=1,
  text opacity=1,
  font=\scriptsize,
  at={(0.03,0.97)},
  anchor=north west,
  draw=lightgray204
},
    width=0.24\textwidth,    
    height=4cm, 
    log basis x={10},
    log basis y={10},
    tick align=outside,
    tick pos=left,
    x grid style={darkgray176},
    xmajorgrids,
    xmin=0.00767270499010924, xmax=2.698565695347128,
    xminorgrids,
    xmode=log,
    xtick style={color=black},
    xtick={1e-05,0.0001,0.001,0.01,0.1,1,10},
    xticklabels={
      \(\displaystyle {10^{-5}}\),
      \(\displaystyle {10^{-4}}\),
      \(\displaystyle {10^{-3}}\),
      \(\displaystyle {10^{-2}}\),
      \(\displaystyle {10^{-1}}\),
      \(\displaystyle {10^{0}}\),
      \(\displaystyle {10^{1}}\)
    },
    y grid style={darkgray176},
    ymajorgrids,
    ymin=0.15, ymax=551.341520850039,
    yminorgrids,
    ymode=log,
    ytick style={color=black},
    ytick={0.01,0.1,1,10,100,1000,10000},
    yticklabels={
      \(\displaystyle {10^{-2}}\),
      \(\displaystyle {10^{-1}}\),
      \(\displaystyle {10^{0}}\),
      \(\displaystyle {10^{1}}\),
      \(\displaystyle {10^{2}}\),
      \(\displaystyle {10^{3}}\),
      \(\displaystyle {10^{4}}\)
    }
]

\addlegendimage{empty legend}
\addlegendentry{$\omega_{\text{LS}} = 1.02$}

\addplot [semithick, blue, mark=x, mark size=3, mark options={solid}]
table {%
2 3652.77905131947
0.8 1148.56191871065
0.666666666666667 524.996044377792
0.571428571428571 147.747549005302
0.5 35.6250696601614
0.4 13.0519554483082
0.266666666666667 6.70671470112185
0.16 3.62850941759641
0.08 1.73743119475338
0.04 0.920418254028223
0.02 0.59611577884161
0.01 0.483094615177154
};
\label{graph:empirical}

\def\xmin{0.00767270499010924}
\def\ymin{0.15}
\def\xminbox{0.017}
\def\xmaxbox{0.42}
\def\yminbox{0.55}
\def\ymaxbox{15}

\def\xclipmin{0.07}
\def\xclipmax{1}
\def\yclipmin{0.01}
\def\yclipmax{100}


\draw[mygreen, line width=3pt] 
  (axis cs:\xminbox,\ymin) -- (axis cs:\xmaxbox,\ymin);

\draw[mygreen, line width=3pt] 
  (axis cs:\xmin,\yminbox) -- (axis cs:\xmin,\ymaxbox);

\begin{scope}
  \clip (axis cs:\xclipmin,\yclipmin) rectangle (axis cs:\xclipmax,\yclipmax);

  \def\yoffset{1.5}

  \addplot[
      thick,
      black,
      domain=0.001:\xclipmax,
      samples=2
  ]{ \yoffset*0.20831919060095*(x/0.08)^(1) };
\end{scope}

\end{axis}
\end{tikzpicture} & 
\begin{tikzpicture}

\definecolor{crimson2143940}{RGB}{214,39,40}
\definecolor{darkgray176}{RGB}{176,176,176}
\definecolor{darkorange25512714}{RGB}{255,127,14}
\definecolor{forestgreen4416044}{RGB}{44,160,44}
\definecolor{lightgray204}{RGB}{204,204,204}
\definecolor{mediumpurple148103189}{RGB}{148,103,189}
\definecolor{steelblue31119180}{RGB}{31,119,180}
\definecolor{mygreen}{RGB}{31,150,31}

\begin{axis}[
ticklabel style={font=\scriptsize},
legend cell align={left},
legend style={
  fill opacity=0.8,
  draw opacity=1,
  text opacity=1,
  font=\scriptsize,
  at={(0.03,0.97)},
  anchor=north west,
  draw=lightgray204
},
width=0.24\textwidth,    
height=4cm, 
log basis x={10},
log basis y={10},
tick align=outside,
tick pos=left,
x grid style={darkgray176},
xmajorgrids,
xmin=0.00767270499010924, xmax=2.698565695347128,
xminorgrids,
xmode=log,
xtick style={color=black},
xtick={0.0001,0.001,0.01,0.1,1,10},
xticklabels={
  \(\displaystyle {10^{-4}}\),
  \(\displaystyle {10^{-3}}\),
  \(\displaystyle {10^{-2}}\),
  \(\displaystyle {10^{-1}}\),
  \(\displaystyle {10^{0}}\),
  \(\displaystyle {10^{1}}\)
},
y grid style={darkgray176},
ymajorgrids,
ymin=0.15, ymax=551.341520850039,
yminorgrids,
ymode=log,
ytick style={color=black},
ytick={1e-06,1e-05,0.0001,0.001,0.01,0.1,1,10,100,1000,10000},
yticklabels={
  \(\displaystyle {10^{-6}}\),
  \(\displaystyle {10^{-5}}\),
  \(\displaystyle {10^{-4}}\),
  \(\displaystyle {10^{-3}}\),
  \(\displaystyle {10^{-2}}\),
  \(\displaystyle {10^{-1}}\),
  \(\displaystyle {10^{0}}\),
  \(\displaystyle {10^{1}}\),
  \(\displaystyle {10^{2}}\),
  \(\displaystyle {10^{3}}\),
  \(\displaystyle {10^{4}}\)
}
]

\addlegendimage{empty legend}
\addlegendentry{$\omega_{\text{LS}} = 1.25$}

\addplot [semithick, blue, mark=x, mark size=3, mark options={solid}]
table {%
2 73.9803977005531
0.8 1.65226655977868
0.666666666666667 1.91909280355518
0.571428571428572 1.69933333202421
0.5 1.47839610004226
0.4 1.07799075831595
0.2666666666666667 0.629216532461755
0.16 0.467746871844141
0.08 0.450328097970018
0.04 0.456352956988253
0.02 0.449443197632836
0.01 0.448084839913565
};
\label{graph:empirical}

\def\xmin{0.00767270499010924}
\def\ymin{0.15}
\def\xminbox{0.24}
\def\xmaxbox{0.72}
\def\yminbox{0.6}
\def\ymaxbox{3}

\def\xclipmin{0.3}
\def\xclipmax{1.5}
\def\yclipmin{0.01}
\def\yclipmax{100}


\draw[mygreen, line width=3pt] 
  (axis cs:\xminbox,\ymin) -- (axis cs:\xmaxbox,\ymin);

\draw[mygreen, line width=3pt] 
  (axis cs:\xmin,\yminbox) -- (axis cs:\xmin,\ymaxbox);

\begin{scope}
  \clip (axis cs:\xclipmin,\yclipmin) rectangle (axis cs:\xclipmax,\yclipmax);

  \def\yoffset{1.0}

  \addplot[
      thick,
      orange,
      domain=0.001:\xclipmax,
      samples=2
  ]{ \yoffset*0.20831919060095*(x/0.3)^(1.5) };
\end{scope}
\end{axis}

\end{tikzpicture}&
\begin{tikzpicture}

\definecolor{crimson2143940}{RGB}{214,39,40}
\definecolor{darkgray176}{RGB}{176,176,176}
\definecolor{darkorange25512714}{RGB}{255,127,14}
\definecolor{forestgreen4416044}{RGB}{44,160,44}
\definecolor{lightgray204}{RGB}{204,204,204}
\definecolor{mediumpurple148103189}{RGB}{148,103,189}
\definecolor{steelblue31119180}{RGB}{31,119,180}
\definecolor{mygreen}{RGB}{31,150,31}

\begin{axis}[
ticklabel style={font=\scriptsize},
legend cell align={left},
legend style={
  fill opacity=0.8,
  draw opacity=1,
  text opacity=1,
  font=\scriptsize,
  at={(0.03,0.97)},
  anchor=north west,
  draw=lightgray204
},
width=0.24\textwidth,    
height=4cm, 
log basis x={10},
log basis y={10},
tick align=outside,
tick pos=left,
x grid style={darkgray176},
xmajorgrids,
xmin=0.00767270499010924, xmax=2.698565695347128,
xminorgrids,
xmode=log,
xtick style={color=black},
xtick={0.0001,0.001,0.01,0.1,1,10},
xticklabels={
  \(\displaystyle {10^{-4}}\),
  \(\displaystyle {10^{-3}}\),
  \(\displaystyle {10^{-2}}\),
  \(\displaystyle {10^{-1}}\),
  \(\displaystyle {10^{0}}\),
  \(\displaystyle {10^{1}}\)
},
y grid style={darkgray176},
ymajorgrids,
ymin=0.15, ymax=551.341520850039,
yminorgrids,
ymode=log,
ytick style={color=black},
ytick={1e-06,1e-05,0.0001,0.001,0.01,0.1,1,10,100,1000,10000},
yticklabels={
  \(\displaystyle {10^{-6}}\),
  \(\displaystyle {10^{-5}}\),
  \(\displaystyle {10^{-4}}\),
  \(\displaystyle {10^{-3}}\),
  \(\displaystyle {10^{-2}}\),
  \(\displaystyle {10^{-1}}\),
  \(\displaystyle {10^{0}}\),
  \(\displaystyle {10^{1}}\),
  \(\displaystyle {10^{2}}\),
  \(\displaystyle {10^{3}}\),
  \(\displaystyle {10^{4}}\)
}
]

\addplot [semithick, blue, mark=x, mark size=3, mark options={solid}]
table {%
           2                 1.204488e+01
         0.8                 7.571783e-01
    0.666667                 7.812070e-01
    0.571429                 6.798718e-01
         0.5                 6.749044e-01
         0.4                 5.669191e-01
    0.266667                 4.783430e-01
        0.16                 4.530715e-01
        0.08                 4.443003e-01
        0.04                 4.461193e-01
        0.02                 4.415414e-01
        0.01                 4.374016e-01
};
\label{graph:empirical}

\def\xmin{0.00767270499010924}
\def\ymin{0.15}
\def\xminbox{0.24}
\def\xmaxbox{0.72}
\def\yminbox{0.6}
\def\ymaxbox{3}

\def\xclipmin{0.3}
\def\xclipmax{1.3}
\def\yclipmin{0.01}
\def\yclipmax{100}


\begin{scope}
  \clip (axis cs:\xclipmin,\yclipmin) rectangle (axis cs:\xclipmax,\yclipmax);

  \def\yoffset{5.0}

  \addplot[
      thick,
      orange,
      domain=0.001:\xclipmax,
      samples=2
  ]{ \yoffset*0.20831919060095*(x/0.3)^(1.5) };
\end{scope}
\end{axis}

\end{tikzpicture}\\
\begin{tikzpicture}

\definecolor{crimson2143940}{RGB}{214,39,40}
\definecolor{darkgray176}{RGB}{176,176,176}
\definecolor{darkorange25512714}{RGB}{255,127,14}
\definecolor{forestgreen4416044}{RGB}{44,160,44}
\definecolor{lightgray204}{RGB}{204,204,204}
\definecolor{mediumpurple148103189}{RGB}{148,103,189}
\definecolor{steelblue31119180}{RGB}{31,119,180}
\definecolor{mygreen}{RGB}{31,150,31}

\begin{axis}[
ticklabel style={font=\scriptsize},
legend cell align={left},
legend style={
  fill opacity=0.8,
  draw opacity=1,
  text opacity=1,
  font=\scriptsize,
  at={(0.03,0.97)},
  anchor=north west,
  draw=lightgray204
},
width=0.24\textwidth,    
height=4cm, 
log basis x={10},
log basis y={10},
tick align=outside,
tick pos=left,
x grid style={darkgray176},
xlabel={{ \scriptsize h vs $\mathcal{W}_2$}},
xmajorgrids,
xmin=0.007, xmax=2.698565695347128,
xminorgrids,
xmode=log,
xtick style={color=black},
xtick={1e-05,0.0001,0.001,0.01,0.1,1,10},
xticklabels={
  \(\displaystyle {10^{-5}}\),
  \(\displaystyle {10^{-4}}\),
  \(\displaystyle {10^{-3}}\),
  \(\displaystyle {10^{-2}}\),
  \(\displaystyle {10^{-1}}\),
  \(\displaystyle {10^{0}}\),
  \(\displaystyle {10^{1}}\)
},
y grid style={darkgray176},
ymajorgrids,
ymin=0.15, ymax=551.341520850039,
yminorgrids,
ymode=log,
ylabel={{\scriptsize Problem \ref{prob:highdim_sg}}},
ytick style={color=black},
ytick={1e-06,1e-05,0.0001,0.001,0.01,0.1,1,10,100,1000,10000},
yticklabels={
  \(\displaystyle {10^{-6}}\),
  \(\displaystyle {10^{-5}}\),
  \(\displaystyle {10^{-4}}\),
  \(\displaystyle {10^{-3}}\),
  \(\displaystyle {10^{-2}}\),
  \(\displaystyle {10^{-1}}\),
  \(\displaystyle {10^{0}}\),
  \(\displaystyle {10^{1}}\),
  \(\displaystyle {10^{2}}\),
  \(\displaystyle {10^{3}}\),
  \(\displaystyle {10^{4}}\)
}
]

\addlegendimage{empty legend}
\addlegendentry{$\omega_{\text{LS}} = 1.00$}

\addplot [semithick, blue, mark=x, mark size=3, mark options={solid}]
table {%
2 255.052971163759
0.8 11.8696811196628
0.666666666666667 8.75695668302162
0.571428571428571 7.16300666818533
0.5 5.97558730294886
0.4 4.71727244854051
0.266666666666667 3.00386759413228
0.16 1.75548882071553
0.08 0.884479178365387
0.04 0.468223459523619
0.02 0.297382680313878
0.01 0.239094413502857
};

\def\xmin{0.007}
\def\ymin{0.15}
\def\xminbox{0.018}
\def\xmaxbox{1.2}
\def\yminbox{0.26}
\def\ymaxbox{15}

\def\xclipmin{0.07}
\def\xclipmax{1}
\def\yclipmin{0.01}
\def\yclipmax{100}


\draw[mygreen, line width=3pt] 
  (axis cs:\xminbox,\ymin) -- (axis cs:\xmaxbox,\ymin);

\draw[mygreen, line width=3pt] 
  (axis cs:\xmin,\yminbox) -- (axis cs:\xmin,\ymaxbox);

\begin{scope}
  \clip (axis cs:\xclipmin,\yclipmin) rectangle (axis cs:\xclipmax,\yclipmax);

  \def\yoffset{1.5}

  \addplot[
      thick,
      black,
      domain=0.001:\xclipmax,
      samples=2
  ]{ \yoffset*0.20831919060095*(x/0.08)^(1) };
\end{scope}

\end{axis}

\end{tikzpicture} & 
\begin{tikzpicture}

\definecolor{crimson2143940}{RGB}{214,39,40}
\definecolor{darkgray176}{RGB}{176,176,176}
\definecolor{darkorange25512714}{RGB}{255,127,14}
\definecolor{forestgreen4416044}{RGB}{44,160,44}
\definecolor{lightgray204}{RGB}{204,204,204}
\definecolor{mediumpurple148103189}{RGB}{148,103,189}
\definecolor{steelblue31119180}{RGB}{31,119,180}
\definecolor{mygreen}{RGB}{31,150,31}

\begin{axis}[
ticklabel style={font=\scriptsize},
legend cell align={left},
legend style={
  fill opacity=0.8,
  draw opacity=1,
  text opacity=1,
  font=\scriptsize,
  at={(0.03,0.97)},
  anchor=north west,
  draw=lightgray204
},
width=0.24\textwidth,    
height=4cm, 
log basis x={10},
log basis y={10},
tick align=outside,
tick pos=left,
x grid style={darkgray176},
xlabel={{ \scriptsize h vs $\mathcal{W}_2$}},
xmajorgrids,
xmin=0.007, xmax=2.698565695347128,
xminorgrids,
xmode=log,
xtick style={color=black},
xtick={1e-05,0.0001,0.001,0.01,0.1,1,10},
xticklabels={
  \(\displaystyle {10^{-5}}\),
  \(\displaystyle {10^{-4}}\),
  \(\displaystyle {10^{-3}}\),
  \(\displaystyle {10^{-2}}\),
  \(\displaystyle {10^{-1}}\),
  \(\displaystyle {10^{0}}\),
  \(\displaystyle {10^{1}}\)
},
y grid style={darkgray176},
ymajorgrids,
ymin=0.15, ymax=551.341520850039,
yminorgrids,
ymode=log,
ytick style={color=black},
ytick={1e-06,1e-05,0.0001,0.001,0.01,0.1,1,10,100,1000,10000},
yticklabels={
  \(\displaystyle {10^{-6}}\),
  \(\displaystyle {10^{-5}}\),
  \(\displaystyle {10^{-4}}\),
  \(\displaystyle {10^{-3}}\),
  \(\displaystyle {10^{-2}}\),
  \(\displaystyle {10^{-1}}\),
  \(\displaystyle {10^{0}}\),
  \(\displaystyle {10^{1}}\),
  \(\displaystyle {10^{2}}\),
  \(\displaystyle {10^{3}}\),
  \(\displaystyle {10^{4}}\)
}
]

\addlegendimage{empty legend}
\addlegendentry{$\omega_{\text{LS}} = 1.07$}

\def\xmin{0.007}
\def\ymin{0.15}
\def\xminbox{0.015}
\def\xmaxbox{0.6}
\def\yminbox{0.3}
\def\ymaxbox{10}

\def\xclipmin{0.07}
\def\xclipmax{1}
\def\yclipmin{0.01}
\def\yclipmax{100}


\draw[mygreen, line width=3pt] 
  (axis cs:\xminbox,\ymin) -- (axis cs:\xmaxbox,\ymin);

\draw[mygreen, line width=3pt] 
  (axis cs:\xmin,\yminbox) -- (axis cs:\xmin,\ymaxbox);

\begin{scope}
  \clip (axis cs:\xclipmin,\yclipmin) rectangle (axis cs:\xclipmax,\yclipmax);

  \def\yoffset{1.5}

  \addplot[
      thick,
      black,
      domain=0.001:\xclipmax,
      samples=2
  ]{ \yoffset*0.20831919060095*(x/0.08)^(1) };
\end{scope}

\addplot [semithick, blue , mark=x, mark size=3, mark options={solid}]
table {%
2 3679.12079851464
0.8 1106.04612174016
0.666666666666667 482.96788894345
0.571428571428571 92.5797660395651
0.5 14.4112213050661
0.4 8.35795200108208
0.266666666666667 5.11552087402865
0.16 2.89021635275794
0.08 1.40352717130183
0.04 0.709957321368005
0.02 0.401059038902651
0.01 0.276907569079649
};
\end{axis}

\end{tikzpicture} &
\begin{tikzpicture}

\definecolor{crimson2143940}{RGB}{214,39,40}
\definecolor{darkgray176}{RGB}{176,176,176}
\definecolor{darkorange25512714}{RGB}{255,127,14}
\definecolor{forestgreen4416044}{RGB}{44,160,44}
\definecolor{lightgray204}{RGB}{204,204,204}
\definecolor{mediumpurple148103189}{RGB}{148,103,189}
\definecolor{steelblue31119180}{RGB}{31,119,180}
\definecolor{mygreen}{RGB}{31,150,31}

\begin{axis}[
ticklabel style={font=\scriptsize},
legend cell align={left},
legend style={
  fill opacity=0.8,
  draw opacity=1,
  text opacity=1,
  font=\scriptsize,
  at={(0.03,0.97)},
  anchor=north west,
  draw=lightgray204
},
width=0.24\textwidth,    
height=4cm, 
log basis x={10},
log basis y={10},
tick align=outside,
tick pos=left,
x grid style={darkgray176},
xlabel={{\scriptsize h vs $\mathcal{W}_2$}},
xmajorgrids,
xmin=0.007, xmax=2.698565695347128,
xminorgrids,
xmode=log,
xtick style={color=black},
xtick={1e-05,0.0001,0.001,0.01,0.1,1,10},
xticklabels={
  \(\displaystyle {10^{-5}}\),
  \(\displaystyle {10^{-4}}\),
  \(\displaystyle {10^{-3}}\),
  \(\displaystyle {10^{-2}}\),
  \(\displaystyle {10^{-1}}\),
  \(\displaystyle {10^{0}}\),
  \(\displaystyle {10^{1}}\)
},
y grid style={darkgray176},
ymajorgrids,
ymin=0.15, ymax=551.341520850039,
yminorgrids,
ymode=log,
ytick style={color=black},
ytick={1e-06,1e-05,0.0001,0.001,0.01,0.1,1,10,100,1000,10000},
yticklabels={
  \(\displaystyle {10^{-6}}\),
  \(\displaystyle {10^{-5}}\),
  \(\displaystyle {10^{-4}}\),
  \(\displaystyle {10^{-3}}\),
  \(\displaystyle {10^{-2}}\),
  \(\displaystyle {10^{-1}}\),
  \(\displaystyle {10^{0}}\),
  \(\displaystyle {10^{1}}\),
  \(\displaystyle {10^{2}}\),
  \(\displaystyle {10^{3}}\),
  \(\displaystyle {10^{4}}\)
}
]

\addlegendimage{empty legend}
\addlegendentry{$\omega_{\text{LS}} = 1.48$}

\addplot [semithick, blue, mark=x, mark size=3, mark options={solid}]
table {%
2 10.8277845596165
0.8 5.84932950458424
0.666666666666667 4.93927057409476
0.571428571428571 3.8901049008072
0.5 3.0149648874185
0.4 2.02241601804163
0.266666666666667 0.913697723600764
0.16 0.425219977845183
0.08 0.244359806997349
0.04 0.222500316213053
0.02 0.220643601833427
};

\def\xmin{0.007}
\def\ymin{0.15}
\def\xminbox{0.06}
\def\xmaxbox{1.2}
\def\yminbox{0.24}
\def\ymaxbox{6}

\def\xclipmin{0.2}
\def\xclipmax{1}
\def\yclipmin{0.01}
\def\yclipmax{100}


\draw[mygreen, line width=3pt] 
  (axis cs:\xminbox,\ymin) -- (axis cs:\xmaxbox,\ymin);

\draw[mygreen, line width=3pt] 
  (axis cs:\xmin,\yminbox) -- (axis cs:\xmin,\ymaxbox);

\begin{scope}
  \clip (axis cs:\xclipmin,\yclipmin) rectangle (axis cs:\xclipmax,\yclipmax);

  \def\yoffset{1.5}

  \addplot[
      thick,
      orange,
      domain=0.001:\xclipmax,
      samples=2
  ]{ \yoffset*0.20831919060095*(x/0.25)^(1.5) };
\end{scope}

\end{axis}

\end{tikzpicture}&
\begin{tikzpicture}

\definecolor{crimson2143940}{RGB}{214,39,40}
\definecolor{darkgray176}{RGB}{176,176,176}
\definecolor{darkorange25512714}{RGB}{255,127,14}
\definecolor{forestgreen4416044}{RGB}{44,160,44}
\definecolor{lightgray204}{RGB}{204,204,204}
\definecolor{mediumpurple148103189}{RGB}{148,103,189}
\definecolor{steelblue31119180}{RGB}{31,119,180}
\definecolor{mygreen}{RGB}{31,150,31}

\begin{axis}[
ticklabel style={font=\scriptsize},
legend cell align={left},
legend style={
  fill opacity=0.8,
  draw opacity=1,
  text opacity=1,
  font=\scriptsize,
  at={(0.03,0.97)},
  anchor=north west,
  draw=lightgray204
},
width=0.24\textwidth,    
height=4cm, 
log basis x={10},
log basis y={10},
tick align=outside,
tick pos=left,
x grid style={darkgray176},
xlabel={{\scriptsize h vs $\mathcal{W}_2$}},
xmajorgrids,
xmin=0.007, xmax=2.698565695347128,
xminorgrids,
xmode=log,
xtick style={color=black},
xtick={1e-05,0.0001,0.001,0.01,0.1,1,10},
xticklabels={
  \(\displaystyle {10^{-5}}\),
  \(\displaystyle {10^{-4}}\),
  \(\displaystyle {10^{-3}}\),
  \(\displaystyle {10^{-2}}\),
  \(\displaystyle {10^{-1}}\),
  \(\displaystyle {10^{0}}\),
  \(\displaystyle {10^{1}}\)
},
y grid style={darkgray176},
ymajorgrids,
ymin=0.15, ymax=551.341520850039,
yminorgrids,
ymode=log,
ytick style={color=black},
ytick={1e-06,1e-05,0.0001,0.001,0.01,0.1,1,10,100,1000,10000},
yticklabels={
  \(\displaystyle {10^{-6}}\),
  \(\displaystyle {10^{-5}}\),
  \(\displaystyle {10^{-4}}\),
  \(\displaystyle {10^{-3}}\),
  \(\displaystyle {10^{-2}}\),
  \(\displaystyle {10^{-1}}\),
  \(\displaystyle {10^{0}}\),
  \(\displaystyle {10^{1}}\),
  \(\displaystyle {10^{2}}\),
  \(\displaystyle {10^{3}}\),
  \(\displaystyle {10^{4}}\)
}
]

\addplot [semithick, blue, mark=x, mark size=3, mark options={solid}]
table {%
 2                5.021403e+01
         0.8                 2.387941e+00
    0.666667                 1.203993e+00
    0.571429                 6.756195e-01
         0.5                 4.991644e-01
         0.4                3.452927e-01
    0.266667                2.588920e-01
        0.16                2.251788e-01
        0.08                2.250596e-01
        0.04               2.197721e-01
        0.02               2.201222e-01
};

\def\xmin{0.007}
\def\ymin{0.15}
\def\xminbox{0.06}
\def\xmaxbox{1.2}
\def\yminbox{0.24}
\def\ymaxbox{6}

\def\xclipmin{0.15}
\def\xclipmax{0.8}
\def\yclipmin{0.01}
\def\yclipmax{100}


\begin{scope}
  \clip (axis cs:\xclipmin,\yclipmin) rectangle (axis cs:\xclipmax,\yclipmax);

  \def\yoffset{5.5}

  \addplot[
      thick,
      orange,
      domain=0.001:\xclipmax,
      samples=2
  ]{ \yoffset*0.20831919060095*(x/0.25)^(1.5) };
\end{scope}

\end{axis}

\end{tikzpicture}
\end{tabular}
    \caption{Simulations for Problems \ref{prob:highdim_gmm} and \ref{prob:highdim_sg} using a learned score network $s_\theta$ which was trained for $100$ epochs. The top row considers Problem \ref{prob:highdim_gmm}, while the bottom row considers Problem \ref{prob:highdim_sg}. $\mathcal{W}_2$ of the considered methods vs the stepsize $h$. Depicted is the empirical curve (\ref{graph:empirical}), a reference line $\mathcal{O}(h)$ (\ref{graph:h1}), and a reference line $\mathcal{O}(h^{3/2})$ (\ref{graph:h1.5}). The green regions on the axes denote the stability region limited by the stepsize on the x-axis, and our $\mathcal{W}_2$ measure on the y-axis.}
    \label{fig:StepsizeVsW2Trainend}
\end{figure}

Figures \ref{fig:StepsizeVsW2TrueScore} and \ref{fig:StepsizeVsW2Trainend} show that the convergence bounds established in Section \ref{sec:ExplicitSchemes} agree with our empirical observations. The higher-order methods also offer a clear advantage in computation time, as shown in Table \ref{tab:computetime}. We measure compute time in neural function evaluations (NFEs), i.e., the number of times the score must be evaluated. This is a standard measure for comparing the compute cost of discretization methods in diffusion models, since evaluating the score network is by far the most expensive step in image generation, cf. e.g. \cite[Figure 3]{song2020score} and \cite[Figure 2]{karras2022elucidating}. Note that per step, Chang's method and SRA1 need three and two times as many NFEs as EM or EI.
\begin{table}[!htp]
    \centering
    \begin{tabular}{|c|c|c|c|c||c|}
    \hline
          &EM &  EI & Chang& SRA1&$\mathcal{W}_2$ floor\\ \hline
         Problem \ref{prob:highdim_gmm}, TS& $800$ & $>800$ & $150$& $50$ & $0.18$\\
         Problem \ref{prob:highdim_sg}, TS& $800$ & $>800$ & $300$& $100$ &$0.09$\\
         Problem \ref{prob:highdim_gmm}, SN& $400$ & $>400$ & $150$& $50$&$0.45$\\
         Problem \ref{prob:highdim_sg}, SN& $400$ & $>400$ & $150$& $50$&$0.24$\\
         \hline
    \end{tabular}
    \caption{Number of neural function evaluations (NFEs) for a discretization method to reach the $\mathcal{W}_2$ floor when using the true score function (TS) to generate samples, and when using the trained score network (SN) to generate samples. Entries $>n$ indicate that the floor was not reached within the tested range of $n$ NFEs.}
    \label{tab:computetime}
\end{table}

Figure~\ref{fig:Dissipativity} exhibits two regimes. For small terminal times, all four schemes behave alike, and $\mathcal{W}_2$ drops from about $55$ at $T=0.4$ to below $5$ at $T=2$. In this regime the initialization error dominates. Its decay is even faster than exponential in $T$, as expected for our schedule, whose forward contraction factor is $\phi(T)=e^{-T^2/4}$, cf.\ Section~\ref{subsec:VerifyingAssumptions}. For $T\gtrsim3$, the error settles at a scheme-dependent plateau, which is the discretization error: approximately $1.10$ for EM and $1.96$ for EI at $h=0.1$, and $0.45$ for SRA1 and $1.92$ for Chang's method at $h=0.5$. All plateaus lie clearly above the floor of our $\mathcal{W}_2$ measure (cf.\ Table~\ref{tab:computetime}). Notably, SRA1 attains a lower plateau than EM with a five times larger step size, i.e.\ with $4$ instead of $10$ NFEs per unit of time. Most importantly, the plateaus do not grow with $T$: between $T=4$ and $T=20$, the errors of EM and EI change by less than $1\%$, and between $T=5$ and $T=10$, the same holds for SRA1 and Chang's method. This is stronger than what Theorem~\ref{thm:MainDissipative} guarantees, and could be explained by Remark \ref{rem:StepwiseConstants}. In contrast, that the admissible step size has to shrink with $T$, as reflected by the threshold $h_0(T)$ in Theorem~\ref{thm:MainDissipative}, is genuine: In our experiments, this is visible for SRA1 and Chang's method with $h=0.5$, which diverge at $T=20$, while EM and EI with $h=0.1$ remain stable over the tested range.

\subsection{Evaluating CIFAR-10}
\label{sec:EvaluatingCifar}
\begin{figure}[]
    \centering
    \sbox{\LatPlotBox}{%
    \setlength{\textwidth}{0.48\textwidth}%
    \setlength{\linewidth}{\textwidth}%
\begin{tikzpicture}

\definecolor{darkgray176}{RGB}{176,176,176}
\definecolor{darkorange25512714}{RGB}{255,127,14}
\definecolor{lightgray204}{RGB}{204,204,204}
\definecolor{steelblue31119180}{RGB}{31,119,180}

\begin{axis}[
legend cell align={left},
ticklabel style={font=\scriptsize},
width=\textwidth,    
height=5cm,
log basis x={10},
log basis y={10},
tick align=outside,
tick pos=left,
x grid style={darkgray176},
xmajorgrids,
xmin=0.0158865646944856, xmax=2.51785082358833,
xminorgrids,
xmode=log,
xtick style={color=black},
xtick={0.0001,0.001,0.01,0.1,1,10,100},
xticklabels={
  \(\displaystyle {10^{-4}}\),
  \(\displaystyle {10^{-3}}\),
  \(\displaystyle {10^{-2}}\),
  \(\displaystyle {10^{-1}}\),
  \(\displaystyle {10^{0}}\),
  \(\displaystyle {10^{1}}\),
  \(\displaystyle {10^{2}}\)
},
y grid style={darkgray176},
ymajorgrids,
ymin=2.84709425938586, ymax=225.55385406444,
yminorgrids,
ymode=log,
ytick style={color=black},
ytick={0.01,0.1,1,10,100,1000},
yticklabels={
  \(\displaystyle {10^{-2}}\),
  \(\displaystyle {10^{-1}}\),
  \(\displaystyle {10^{0}}\),
  \(\displaystyle {10^{1}}\),
  \(\displaystyle {10^{2}}\),
  \(\displaystyle {10^{3}}\)
}
]
\addplot[semithick, Chang_color, mark=x, mark size=2, mark options={solid}]
table {%
2 47.3065758802365
0.8 23.673214822867
0.666666666666667 18.0466743821854
0.571428571428571 14.4207105324115
0.5 11.7709189588701
0.4 8.29455715121928
0.266666666666667 4.79642615639717
0.16 3.82242118384868
0.08 3.7955341382379
0.04 3.81332285736015
0.02 3.82822726328553
};
\addplot [semithick, EM_color, mark=o, mark size=2, mark options={solid}]
table {%
2 184.900872770603
0.8 37.3110279665712
0.666666666666667 29.7583645358376
0.571428571428571 24.5437087144386
0.5 20.9974079501988
0.4 16.0689451490327
0.266666666666667 9.71249971902157
0.16 5.40426333322872
0.08 3.56850624513945
0.04 3.47306680312933
0.02 3.59653941705914
};

\addplot[semithick, EI_color, mark=triangle, mark size=2, mark options={solid}]
table {%
2 2963.73608398438
0.8 621.70166015625
0.666666666666667 191.931228637695
0.571428571428571 66.7673873901367
0.5 39.8491668701172
0.4 22.6286869049072
0.266666666666667 12.9211254119873
0.16 7.49002075195312
0.08 4.24918794631958
0.04 3.53092980384827
0.02 3.56267809867859
};

\addplot[semithick,SRA1_color, mark=square, mark size=2, mark options={solid}]
table {%
2.00000   92.19729614
0.80000   37.44721603
0.66667   29.86562538
0.57143   24.42803764
0.50000   20.23633003
0.40000   14.45434475
0.26667   7.63691616
0.16000   4.40042257
0.08000   3.88121343
0.04000   3.82593703
0.02000   3.83825779

};

\def\xlinearclipmin{0.07}
\def\xlinearclipmax{0.8}
\def\ylinearclipmin{0.1}
\def\ylinearclipmax{100}
\begin{scope}
  \clip (axis cs:\xlinearclipmin,\ylinearclipmin) rectangle (axis cs:\xlinearclipmax,\ylinearclipmax);

  \def\yoffset{15}

  \addplot[
      thick,
      red,
      domain=0.001:\xlinearclipmax,
      samples=2
  ] { \yoffset*0.20831919060095*((x)/0.04) };
  \label{graph:Featureh1}
\end{scope}

\def\xonepointfiveclipmin{0.5}
\def\xonepointfiveclipmax{1.5}
\def\yonepointfiveclipmin{0.1}
\def\yonepointfiveclipmax{100}
\begin{scope}
  \clip (axis cs:\xonepointfiveclipmin,\yonepointfiveclipmin) rectangle (axis cs:\xonepointfiveclipmax,\yonepointfiveclipmax);

  \def\yoffset{0.5}

  \addplot[
      thick,
      blue,
      domain=0.001:\xonepointfiveclipmax,
      samples=2
  ] { \yoffset*0.20831919060095*((x)/0.04)^(1.5) };
  \label{graph:Featureh1.5}
\end{scope}

\end{axis}
\end{tikzpicture}%
    }
    \begin{tabular}{cc}
    \begin{subfigure}[t]{0.48\textwidth}
        \centering
\begin{tikzpicture}

\definecolor{darkgray176}{RGB}{176,176,176}
\definecolor{lightgray204}{RGB}{204,204,204}
\definecolor{steelblue31119180}{RGB}{31,119,180}

\begin{axis}[
legend cell align={left},
ticklabel style={font=\scriptsize},
width=\textwidth,    
height=5cm, 
log basis x={10},
log basis y={10},
tick align=outside,
tick pos=left,
x grid style={darkgray176},
xmajorgrids,
xmin=0.00588656469448564, xmax=2.698565695347128,
xminorgrids,
xmode=log,
xtick style={color=black},
xtick={0.0001,0.001,0.01,0.1,1,10},
xticklabels={
  \(\displaystyle {10^{-4}}\),
  \(\displaystyle {10^{-3}}\),
  \(\displaystyle {10^{-2}}\),
  \(\displaystyle {10^{-1}}\),
  \(\displaystyle {10^{0}}\),
  \(\displaystyle {10^{1}}\)
},
y grid style={darkgray176},
ymajorgrids,
ymin=100, ymax=4000,
yminorgrids,
ymode=log,
ytick style={color=black},
ytick={1,10,100,1000,10000,100000},
yticklabels={
  \(\displaystyle {10^{0}}\),
  \(\displaystyle {10^{1}}\),
  \(\displaystyle {10^{2}}\),
  \(\displaystyle {10^{3}}\),
  \(\displaystyle {10^{4}}\),
  \(\displaystyle {10^{5}}\)
}
]
\addplot [semithick, Chang_color, mark=x, mark size=2, mark options={solid}]
table {%
2 2937.52661132812
0.8 2568.51147460938
0.666666666666667 2328.36279296875
0.571428571428571 2083.8583984375
0.5 1843.86730957031
0.4 1417.94079589844
0.266666666666667 730.255859375
0.16 287.950500488281
0.08 142.344955444336
0.04 162.773574829102
0.02 181.348114013672
0.01 183.485885620117
};
\label{graph:HO}

\addplot [semithick, EM_color, mark=o, mark size=2, mark options={solid}]
table {%
2 3312.80932617188
0.8 2699.42333984375
0.666666666666667 2489.7392578125
0.571428571428571 2297.564453125
0.5 2105.1513671875
0.4 1781.09118652344
0.266666666666667 1175.37902832031
0.16 598.337768554688
0.08 272.890930175781
0.04 219.498443603516
0.02 185.43537902832
0.01 173.456832885742
};
\label{graph:EM}

\addplot[semithick, EI_color, mark=triangle, mark size=2, mark options={solid}]
table {%
2 3613.03662109375
0.8 3560.06787109375
0.666666666666667 3380.78686523438
0.571428571428571 2821.73706054688
0.5 2214.58935546875
0.4 1680.62744140625
0.266666666666667 1180.091796875
0.16 733.313354492188
0.08 321.427429199219
0.04 150.26203918457
0.02 129.116287231445
0.01 135.520202636719
};
\label{graph:EI}

\addplot [semithick,SRA1_color, mark=square, mark size=2, mark options={solid}]
table {%
2.00000   3405.85107422
0.80000   2879.30029297
0.66667   2714.20654297
0.57143   2553.05566406
0.50000   2385.16210938
0.40000   2085.68457031
0.26667   1418.80126953
0.16000   613.03222656
0.08000   155.53099060
0.04000   165.31933594
0.02000   167.66433716
0.01000   186.42512512

};
\label{graph:SRA1}

\def\xlinearclipmin{0.03}
\def\xlinearclipmax{0.8}
\def\ylinearclipmin{0.1}
\def\ylinearclipmax{2000}
\begin{scope}
  \clip (axis cs:\xlinearclipmin,\ylinearclipmin) rectangle (axis cs:\xlinearclipmax,\ylinearclipmax);

  \def\yoffset{500}

  \addplot[
      thick,
      red,
      domain=0.001:\xlinearclipmax,
      samples=2
  ] { \yoffset*((x)/0.05) };
  \label{graph:Featureh1}
\end{scope}

\def\xonepointfiveclipmin{0.2}
\def\xonepointfiveclipmax{1.5}
\def\yonepointfiveclipmin{0.1}
\def\yonepointfiveclipmax{2000}
\begin{scope}
  \clip (axis cs:\xonepointfiveclipmin,\yonepointfiveclipmin) rectangle (axis cs:\xonepointfiveclipmax,\yonepointfiveclipmax);

  \def\yoffset{15}

  \addplot[
      thick,
      blue,
      domain=0.001:\xonepointfiveclipmax,
      samples=2
  ] { \yoffset*((x)/0.04)^(1.5) };
  \label{graph:Featureh1.5}
\end{scope}
\end{axis}

\end{tikzpicture}
        \caption{Problem \ref{prob:CIFAR10}: PIX $h$ vs $\mathcal{W}_2$}
        \label{fig:PixW2}
    \end{subfigure}
    &
    \begin{subfigure}[t]{0.48\textwidth}
        \centering
\begin{tikzpicture}

\definecolor{darkgray176}{RGB}{176,176,176}
\definecolor{lightgray204}{RGB}{204,204,204}
\definecolor{steelblue31119180}{RGB}{31,119,180}

\begin{axis}[
legend cell align={left},
ticklabel style={font=\scriptsize},
width=\textwidth,    
height=5cm, 
log basis x={10},
log basis y={10},
tick align=outside,
tick pos=left,
x grid style={darkgray176},
xmajorgrids,
xmin=0.00588656469448564, xmax=2.698565695347128,
xminorgrids,
xmode=log,
xtick style={color=black},
xtick={0.0001,0.001,0.01,0.1,1,10},
xticklabels={
  \(\displaystyle {10^{-4}}\),
  \(\displaystyle {10^{-3}}\),
  \(\displaystyle {10^{-2}}\),
  \(\displaystyle {10^{-1}}\),
  \(\displaystyle {10^{0}}\),
  \(\displaystyle {10^{1}}\)
},
y grid style={darkgray176},
ymajorgrids,
ymin=3, ymax=600,
yminorgrids,
ymode=log,
ytick style={color=black},
ytick={1,10,100,1000,10000},
yticklabels={
  \(\displaystyle {10^{0}}\),
  \(\displaystyle {10^{1}}\),
  \(\displaystyle {10^{2}}\),
  \(\displaystyle {10^{3}}\),
  \(\displaystyle {10^{4}}\)
}
]
\addplot [semithick, Chang_color, mark=x, mark size=2, mark options={solid}]
table {%
2 429.678804182703
0.8 456.551438340933
0.666666666666667 437.70629945293
0.571428571428571 415.482655095093
0.5 393.369957576787
0.4 348.636261642533
0.266666666666667 248.443422284222
0.16 176.344433087373
0.08 110.286351032415
0.04 54.4828444112551
0.02 17.1096964650802
0.01 4.959536484283376
};

\addplot [semithick,  EM_color, mark=o, mark size=2, mark options={solid}]
table {%
2 400.066657905955
0.8 449.102332821699
0.666666666666667 435.323988817761
0.571428571428571 419.696361148202
0.5 401.003689205573
0.4 367.339974976541
0.266666666666667 293.870803472294
0.16 223.287503366157
0.08 142.377942394162
0.04 80.0198686280442
0.02 37.8889804603344
0.01 14.077513472033104
};

\addplot [semithick, EI_color, mark=triangle, mark size=2, mark options={solid}]
table {%
2 458.06445514707
0.8 377.502552701174
0.666666666666667 316.690671717084
0.571428571428571 297.166503987396
0.5 284.695749087358
0.4 271.958292537685
0.266666666666667 241.547984579159
0.16 195.388689865142
0.08 126.490871078701
0.04 74.9603165691726
0.02 41.3148640570359
0.01 19.36569228066014
};

\addplot [semithick, SRA1_color, mark=square, mark size=2, mark options={solid}]
table {%
           2               4.902853e+02
         0.8                4.771725e+02
    0.666667               4.701428e+02
    0.571429               4.618271e+02
         0.5               4.511877e+02
         0.4               4.269840e+02
    0.266667              3.552266e+02
        0.16               2.334411e+02
        0.08               1.415001e+02
        0.04              7.997167e+01
        0.02              3.396082e+01
        0.01              8.727343e+00
};

\end{axis}

\end{tikzpicture}
        \caption{Problem \ref{prob:CIFAR10}: PIX $h$ vs FID}
        \label{fig:PixFID}
    \end{subfigure}
    \\[5ex]
    \begin{subfigure}[t]{0.48\textwidth}
        \centering
        \usebox{\LatPlotBox}
        \caption{Problem \ref{prob:LatentCIFAR10}: LAT $h$ vs $\mathcal{W}_2$}
        \label{fig:LatW2}
    \end{subfigure}
    &
    \begin{subfigure}[t]{0.48\textwidth}
        \centering
        \begin{minipage}[b][\ht\LatPlotBox][c]{\linewidth}
            \centering
            \setlength{\tabcolsep}{4pt}
            \begin{adjustbox}{max width=\linewidth}
            \begin{tabular}{|c|c|c|c|c|c|}
            \hline
                \multicolumn{2}{|c|}{} & EM & EI & Chang & SRA1 \\
            \hline \hline
                \multirow{2}{*}{\subref{fig:PixW2}} & $\omega_{LS}$ & $1.10$ & $1.07$ & $1.64$ & $1.50$ \\
                \cline{2-6}
                 & $h$ range & $0.08$-$0.6$ & $0.04$-$0.6$ & $0.16$-$0.6$ & $0.08$-$0.5$ \\
            \hline
                \multirow{2}{*}{\subref{fig:PixFID}} & $\omega_{LS}$ & $0.79$ & $0.65$ & $1.03$ & $0.91$ \\
                \cline{2-6}
                 & $h$ range & \multicolumn{4}{c|}{$0.01$-$0.6$} \\
            \hline
                \multirow{2}{*}{\subref{fig:LatW2}} & $\omega_{LS}$ & $1.05$ & $1.02$ & $1.50$ & $1.46$ \\
                \cline{2-6}
                 & $h$ range & $0.08$-$0.8$ & $0.08$-$0.4$ & $0.4$-$0.8$ & $0.2$-$0.8$ \\
            \hline
            \end{tabular}
            \end{adjustbox}
        \end{minipage}
        \caption{Least squares slopes $\omega_{LS}$}
        \label{tab:ConvergencerateCIFAR10}
    \end{subfigure}
    \end{tabular}
    \caption{Simulations for Problem \ref{prob:CIFAR10} and Problem \ref{prob:LatentCIFAR10}. We show the empirical performance of the EM method (\ref{graph:EM}), the EI method (\ref{graph:EI}), Chang's method (\ref{graph:HO}), and the SRA1 method (\ref{graph:SRA1}), as well as the reference lines $\mathcal{O}(h)$ (\ref{graph:Featureh1}) and $O(h^{3/2})$ (\ref{graph:Featureh1.5}). In the graphs (a) and (b), we generate in PIXel space, while for graph (c) we generate in LATent space, as described for Problem \ref{prob:LatentCIFAR10}. Table (d) lists the least squares slopes of the curves in (a)-(c), together with the range of $h$ used for the fit.}
    \label{fig:StepsizeVsW2andFID_CIFAR10}
\end{figure}
In Figure \ref{fig:StepsizeVsW2andFID_CIFAR10}, we evaluate Problem \ref{prob:CIFAR10} and Problem \ref{prob:LatentCIFAR10}. The left panel reports performance for Problem \ref{prob:CIFAR10} in terms of the $\mathcal{W}_2$ distance, while the middle panel shows the FID values for Problem \ref{prob:CIFAR10}. The right panel presents the $\mathcal{W}_2$ performance of the analyzed discretization methods for Problem \ref{prob:LatentCIFAR10}. In this setup, samples are generated directly in the $2048$-dimensional latent space in which the FID is computed and are evaluated using our $\mathcal{W}_2$ measure, effectively recovering the square root of the FID values. In Table \ref{tab:ConvergencerateCIFAR10}, we depict the least squares convergence order of the slopes depicted in Figure \ref{fig:StepsizeVsW2andFID_CIFAR10}. The training setup is similar to the one used in toy problems with the only difference being the amount of training done, i.e. we use the \emph{ncsn++} network architecture and train it for 2000 epochs for Problem \ref{prob:CIFAR10}, and 1000 epochs for Problem \ref{prob:LatentCIFAR10}.

Figure \ref{fig:StepsizeVsW2andFID_CIFAR10} demonstrates that the higher-order methods retain their faster convergence for Problem \ref{prob:CIFAR10}, as evidenced by their steeper slope compared to EM and EI, even though we reach the floor of our measure before this faster convergence translates into a compute-time advantage. Specifically, Chang's method reaches the floor in 150 NFEs, and SRA1 reaches the floor in 100 NFEs, while EI attains comparable quality in 100 NFEs. The advantage regarding the steeper slope is less pronounced when performance is measured using the FID score, although the gap widens as the step size shrinks. This is not unexpected, since the relationship between distances in pixel space and those in latent space is not well understood, and our theoretical predictions for Problem \ref{prob:CIFAR10} hold only in pixel space. Turning to Problem \ref{prob:LatentCIFAR10}, which removes this mismatch between generation and measure space by generating directly in the measure space, the advantage of the higher-order scheme is again clearly visible. Nonetheless, the situation mirrors the first case: Chang's method reaches the floor in 150 NFEs, and SRA1 in 200 NFEs, while EM requires only 100 NFEs. Lowering this floor, e.g.\ by changing the measure or by comparing more samples, should translate the steeper rate directly into faster wall-clock convergence for the considered higher order methods relative to EM and EI. Note that in our case, comparing more samples is not possible, as the dataset CIFAR10 only contains $50\,000$ samples.
\subsection{Evaluating ImageNet}
\begin{figure}[]
    \centering
    \begin{tabular}{cc}

    \begin{subfigure}[t]{0.48\textwidth}
        \centering
\begin{tikzpicture}

\definecolor{darkgray176}{RGB}{176,176,176}
\definecolor{lightgray204}{RGB}{204,204,204}
\definecolor{steelblue31119180}{RGB}{31,119,180}

\begin{axis}[
legend cell align={left},
ticklabel style={font=\scriptsize},
width=\textwidth,    
height=5cm, 
log basis x={10},
log basis y={10},
tick align=outside,
tick pos=left,
x grid style={darkgray176},
xmajorgrids,
xmin=01, xmax=82,
xminorgrids,
xmode=log,
xtick style={color=black},
xtick={0.0001,0.001,0.01,0.1,1,10},
xticklabels={
  \(\displaystyle {10^{-4}}\),
  \(\displaystyle {10^{-3}}\),
  \(\displaystyle {10^{-2}}\),
  \(\displaystyle {10^{-1}}\),
  \(\displaystyle {10^{0}}\),
  \(\displaystyle {10^{1}}\)
},
y grid style={darkgray176},
ymajorgrids,
ymin=2, ymax=3e+04,
yminorgrids,
ymode=log,
ytick style={color=black},
ytick={1,10,100,1000,10000},
yticklabels={
  \(\displaystyle {10^{0}}\),
  \(\displaystyle {10^{1}}\),
  \(\displaystyle {10^{2}}\),
  \(\displaystyle {10^{3}}\),
  \(\displaystyle {10^{4}}\)
}
]
\addplot [semithick, Chang_color, mark=x, mark size=2, mark options={solid}]
table {%
  2.16914 3.95535
  4.30881 4.04462
  8.49896 4.62566
  16.5172 8.75171
  26.4413 8.49294
  37.5848 14.2805
  45.0078 115.659
  49.8167 156.533
  55.5917 493.609
  62.4722 1600.4
  79.998 23304.1
};

\addplot [semithick,  EM_color, mark=o, mark size=2, mark options={solid}]
table {%

  1.08824 5.59636
  2.16914 8.1706
  4.30881 15.8127
  8.49896 24.4763
  16.5172 22.1928
  26.4413 19.5833
  37.5848 27.1388
  45.0078 416.182
  49.8167 546.266
  55.5917 1220.42
  62.4722 2409.05
  79.998 15254.9
};

\addplot [semithick, SRA1_color, mark=square, mark size=2, mark options={solid}]
table {%
  2.16914 3.91493
  4.30881 3.88056
  8.49896 3.83714
  16.5172 4.03823
  26.4413 13.2619
  37.5848 36.8082
  45.0078 269.788
  49.8167 1057.64
  55.5917 4018.62
  62.4722 11600.1
  79.998 27335.1
};

\addplot [semithick, EI_color, mark=triangle, mark size=2, mark options={solid}]
table {%
  2.16914 11.0933
  4.30881 21.5747
  8.49896 27.5776
  16.5172 27.6196
  26.4413 29.3286
  37.5848 34.2581
  45.0078 42.1791
  49.8167 67.3064
  55.5917 123.535
  62.4722 173.382
  79.998 8767.42
};

\def\xlinearclipmin{1.2}
\def\xlinearclipmax{8}
\def\ylinearclipmin{5}
\def\ylinearclipmax{5000}
\begin{scope}
  \clip (axis cs:\xlinearclipmin,\ylinearclipmin) rectangle (axis cs:\xlinearclipmax,\ylinearclipmax);

  \def\yoffset{3}

  \addplot[
      thick,
      red,
      domain=0.001:\xlinearclipmax,
      samples=2
  ] { \yoffset*0.20831919060095*((x)/0.04) };
  \label{graph:Featureh1}
\end{scope}

\def\xonepointfiveclipmin{25}
\def\xonepointfiveclipmax{70}
\def\yonepointfiveclipmin{3}
\def\yonepointfiveclipmax{1000}
\begin{scope}
  \clip (axis cs:\xonepointfiveclipmin,\yonepointfiveclipmin) rectangle (axis cs:\xonepointfiveclipmax,\yonepointfiveclipmax);

  \def\yoffset{0.001}

  \addplot[
      thick,
      blue,
      domain=0.001:\xonepointfiveclipmax,
      samples=2
  ] { \yoffset*0.20831919060095*((x)/0.04)^(1.5) };
  \label{graph:Featureh1.5}
\end{scope}
\end{axis}

\end{tikzpicture}
        \caption{Problem \ref{prob:ImageNet64}: $h_{\max}$ vs $\mathcal{W}_2$}
        \label{fig:Img64W2}
    \end{subfigure} &
    \begin{subfigure}[t]{0.48\textwidth}
        \centering
\begin{tikzpicture}

\definecolor{darkgray176}{RGB}{176,176,176}
\definecolor{lightgray204}{RGB}{204,204,204}
\definecolor{steelblue31119180}{RGB}{31,119,180}

\begin{axis}[
legend cell align={left},
ticklabel style={font=\scriptsize},
width=\textwidth,    
height=5cm, 
log basis x={10},
log basis y={10},
tick align=outside,
tick pos=left,
x grid style={darkgray176},
xmajorgrids,
xmin=1, xmax=82,
xminorgrids,
xmode=log,
xtick style={color=black},
xtick={0.0001,0.001,0.01,0.1,1,10},
xticklabels={
  \(\displaystyle {10^{-4}}\),
  \(\displaystyle {10^{-3}}\),
  \(\displaystyle {10^{-2}}\),
  \(\displaystyle {10^{-1}}\),
  \(\displaystyle {10^{0}}\),
  \(\displaystyle {10^{1}}\)
},
y grid style={darkgray176},
ymajorgrids,
ymin=2, ymax=3e+04,
yminorgrids,
ymode=log,
ytick style={color=black},
ytick={1,10,100,1000,10000},
yticklabels={
  \(\displaystyle {10^{0}}\),
  \(\displaystyle {10^{1}}\),
  \(\displaystyle {10^{2}}\),
  \(\displaystyle {10^{3}}\),
  \(\displaystyle {10^{4}}\)
}
]
\addplot [semithick, Chang_color, mark=x, mark size=2, mark options={solid}]
table {%
  
  2.16914 5.69969
  4.30881 5.58651
  8.49896 5.23863
  16.5172 2.85536
  26.4413 6.97872
  37.5848 28.255
  45.0078 150.894
  49.8167 334.772
  55.5917 800.81
  62.4722 2313.35
  79.998 26961.4
};

\addplot [semithick,  EM_color, mark=o, mark size=2, mark options={solid}]
table {%
  
  2.16914 9.6294
  4.30881 13.0123
  8.49896 17.6933
  16.5172 28.6462
  26.4413 53.6065
  37.5848 162.643
  45.0078 364.667
  49.8167 687.816
  55.5917 1740.25
  62.4722 3580.51
  79.998 17659.2
};

\addplot [semithick, SRA1_color, mark=square, mark size=2, mark options={solid}]
table {%
  2.16914 5.70546
  4.30881 5.63701
  8.49896 5.59257
  16.5172 5.92822
  26.4413 13.9537
  37.5848 67.5866
  45.0078 329.561
  49.8167 867.504
  55.5917 2928.27
  62.4722 9245.17
  79.998 31625
};

\addplot [semithick, EI_color, mark=triangle, mark size=2, mark options={solid}]
table {%
  
  2.16914 11.5222
  4.30881 15.6475
  8.49896 23.533
  16.5172 44.92
  26.4413 85.883
  37.5848 160.298
  45.0078 230.047
  49.8167 283.129
  55.5917 352.152
  62.4722 472.604
  79.998 10170.3
};

\def\xlinearclipmin{5}
\def\xlinearclipmax{20}
\def\ylinearclipmin{10}
\def\ylinearclipmax{5000}
\begin{scope}
  \clip (axis cs:\xlinearclipmin,\ylinearclipmin) rectangle (axis cs:\xlinearclipmax,\ylinearclipmax);

  \def\yoffset{1.5}

  \addplot[
      thick,
      red,
      domain=0.001:\xlinearclipmax,
      samples=2
  ] { \yoffset*0.20831919060095*((x)/0.04) };
  \label{graph:Featureh1}
\end{scope}

\def\xonepointfiveclipmin{8}
\def\xonepointfiveclipmax{25}
\def\yonepointfiveclipmin{8}
\def\yonepointfiveclipmax{1000}
\begin{scope}
  \clip (axis cs:\xonepointfiveclipmin,\yonepointfiveclipmin) rectangle (axis cs:\xonepointfiveclipmax,\yonepointfiveclipmax);

  \def\yoffset{0.01}

  \addplot[
      thick,
      blue,
      domain=0.001:\xonepointfiveclipmax,
      samples=2
  ] { \yoffset*0.20831919060095*((x)/0.04)^(1.5) };
  \label{graph:Featureh1.5}
\end{scope}

\end{axis}
\end{tikzpicture}
        \caption{Problem \ref{prob:ImageNet512}: $h_{\max}$ vs $\mathcal{W}_2$}
        \label{fig:Img512W2}
    \end{subfigure}\\[5ex]
    \begin{subfigure}[t]{0.48\textwidth}
        \centering
\begin{tikzpicture}

\definecolor{darkgray176}{RGB}{176,176,176}
\definecolor{lightgray204}{RGB}{204,204,204}
\definecolor{steelblue31119180}{RGB}{31,119,180}

\begin{axis}[
legend cell align={left},
ticklabel style={font=\scriptsize},
width=\textwidth,    
height=5cm, 
log basis x={10},
log basis y={10},
tick align=outside,
tick pos=left,
x grid style={darkgray176},
xmajorgrids,
xmin=1, xmax=82,
xminorgrids,
xmode=log,
xtick style={color=black},
xtick={0.0001,0.001,0.01,0.1,1,10},
xticklabels={
  \(\displaystyle {10^{-4}}\),
  \(\displaystyle {10^{-3}}\),
  \(\displaystyle {10^{-2}}\),
  \(\displaystyle {10^{-1}}\),
  \(\displaystyle {10^{0}}\),
  \(\displaystyle {10^{1}}\)
},
y grid style={darkgray176},
ymajorgrids,
ymin=1, ymax=500,
yminorgrids,
ymode=log,
ytick style={color=black},
ytick={1,10,100,1000,10000},
yticklabels={
  \(\displaystyle {10^{0}}\),
  \(\displaystyle {10^{1}}\),
  \(\displaystyle {10^{2}}\),
  \(\displaystyle {10^{3}}\),
  \(\displaystyle {10^{4}}\)
}
]
\addplot [semithick, Chang_color, mark=x, mark size=2, mark options={solid}]
table {%
  2.16914 1.6598
  4.30881 1.85592
  8.49896 7.0984
  16.5172 83.7273
  26.4413 282.521
  37.5848 408.78
  45.0078 413.244
  49.8167 395.437
  55.5917 401.837
  62.4722 427.725
  79.998 419.778
};

\addplot [semithick,  EM_color, mark=o, mark size=2, mark options={solid}]
table {%
 1.08824 4.26674
  2.16914 15.9099
  4.30881 76.4821
  8.49896 210.608
  16.5172 435.816
  26.4413 458.513
  37.5848 343.955
  45.0078 330.859
  49.8167 396.014
  55.5917 429.209
  62.4722 441.473
  79.998 419.709
};

\addplot [semithick, SRA1_color, mark=square, mark size=2, mark options={solid}]
table {%
  2.16914 1.6897
  4.30881 1.6909
  8.49896 1.70663
  16.5172 3.01412
  26.4413 94.1392
  37.5848 384.55
  45.0078 381.803
  49.8167 424.061
  55.5917 432.495
  62.4722 430.164
  79.998 420.598
};

\addplot [semithick, EI_color, mark=triangle, mark size=2, mark options={solid}]
table {%
  2.16914 30.972
  4.30881 122.007
  8.49896 254.725
  16.5172 450.244
  26.4413 475.022
  37.5848 448.014
  45.0078 401.344
  49.8167 368.097
  55.5917 373.836
  62.4722 401.227
  79.998 419.84
};
\label{graph:EI}

\end{axis}

\end{tikzpicture}
        \caption{Problem \ref{prob:ImageNet64}: $h_{\max}$ vs FID}
        \label{fig:Img64FID}
    \end{subfigure}&
    \begin{subfigure}[t]{0.48\textwidth}
        \centering
\begin{tikzpicture}

\definecolor{darkgray176}{RGB}{176,176,176}
\definecolor{lightgray204}{RGB}{204,204,204}
\definecolor{steelblue31119180}{RGB}{31,119,180}

\begin{axis}[
legend cell align={left},
ticklabel style={font=\scriptsize},
width=\textwidth,    
height=5cm, 
log basis x={10},
log basis y={10},
tick align=outside,
tick pos=left,
x grid style={darkgray176},
xmajorgrids,
xmin=1, xmax=82,
xminorgrids,
xmode=log,
xtick style={color=black},
xtick={0.0001,0.001,0.01,0.1,1,10},
xticklabels={
  \(\displaystyle {10^{-4}}\),
  \(\displaystyle {10^{-3}}\),
  \(\displaystyle {10^{-2}}\),
  \(\displaystyle {10^{-1}}\),
  \(\displaystyle {10^{0}}\),
  \(\displaystyle {10^{1}}\)
},
y grid style={darkgray176},
ymajorgrids,
ymin=1, ymax=500,
yminorgrids,
ymode=log,
ytick style={color=black},
ytick={1,10,100,1000,10000},
yticklabels={
  \(\displaystyle {10^{0}}\),
  \(\displaystyle {10^{1}}\),
  \(\displaystyle {10^{2}}\),
  \(\displaystyle {10^{3}}\),
  \(\displaystyle {10^{4}}\)
}
]
\addplot [semithick, Chang_color, mark=x, mark size=2, mark options={solid}]
table {%
  2.16914 2.83941
  4.30881 2.86043
  8.49896 3.86959
  16.5172 29.0692
  26.4413 202.131
  37.5848 380.636
  45.0078 379.892
  49.8167 388.746
  55.5917 391.978
  62.4722 399.658
  79.998 410.122
};

\addplot [semithick,  EM_color, mark=o, mark size=2, mark options={solid}]
table {%
  2.16914 9.28495
  4.30881 43.3117
  8.49896 168.288
  16.5172 353.469
  26.4413 399.363
  37.5848 390.258
  45.0078 413.207
  49.8167 405.414
  55.5917 461.099
  62.4722 339.836
  79.998 409.561
};

\addplot [semithick, SRA1_color, mark=square, mark size=2, mark options={solid}]
table {%
  2.16914 2.86233
  4.30881 2.84974
  8.49896 2.87072
  16.5172 3.02416
  26.4413 37.8561
  37.5848 394.042
  45.0078 382.574
  49.8167 390.072
  55.5917 420.949
  62.4722 403.808
  79.998 410.878
};

\addplot [semithick, EI_color, mark=triangle, mark size=2, mark options={solid}]
table {%
  2.16914 17.174
  4.30881 76.7916
  8.49896 217.336
  16.5172 383.15
  26.4413 393.586
  37.5848 380.869
  45.0078 434.22
  49.8167 456.827
  55.5917 458.732
  62.4722 440.248
  79.998 408.494
};

\end{axis}

\end{tikzpicture}
        \caption{Problem \ref{prob:ImageNet512}: $h_{\max}$ vs FID}
        \label{fig:Img512FID}
    \end{subfigure}
\end{tabular}

    \caption{Simulations for Problems \ref{prob:ImageNet64} and \ref{prob:ImageNet512}, measured in $\mathcal{W}_2$ in the top row, and in FID in the bottom row. We show the empirical performance of the EM method (\ref{graph:EM}), the EI method (\ref{graph:EI}), Chang's method (\ref{graph:HO}), and the SRA1 method (\ref{graph:SRA1}), as well as the reference lines $\mathcal{O}(h)$ (\ref{graph:Featureh1}) and $O(h^{3/2})$ (\ref{graph:Featureh1.5}).} 
    \label{fig:ImageNet}
\end{figure}
\begin{table}[!htp]
    \centering
    \begin{tabular}{|c|c|c|c|c||c|}
    \hline
          &EM &  EI & Chang& SRA1& measure floor\\ \hline
         Problem \ref{prob:ImageNet64}, $\mathcal{W}_2$ (Fig. \ref{fig:Img64W2})& $>400$ & $>200$ & $300$& $50$ &$3.9$\\
         Problem \ref{prob:ImageNet512}, $\mathcal{W}_2$ (Fig. \ref{fig:Img512W2})& $>200$ & $>200$ & $150$& $50$&$5.6$\\
         Problem \ref{prob:ImageNet64}, FID (Fig. \ref{fig:Img64FID})& $>400$ & $>200$ & $600$& $100$ & $1.68$\\
         Problem \ref{prob:ImageNet512}, FID (Fig. \ref{fig:Img512FID})& $>200$ & $>200$ & $300$& $50$&$2.84$\\
         \hline
    \end{tabular}
    \caption{Number of neural function evaluations (NFEs) required by each method to reach the error floor given in the last column. Entries $>n$ indicate that the floor was not reached within the tested range of $n$ NFEs.}
    \label{tab:computetimeImageNet}
\end{table}
In Figure \ref{fig:ImageNet}, we evaluate Problems \ref{prob:ImageNet64} and \ref{prob:ImageNet512} in terms of $\mathcal{W}_2$ and FID, and Table \ref{tab:computetimeImageNet} lists the number of NFEs each method requires to reach the respective error floor. On this benchmark, the higher-order methods yield a substantial compute-time advantage: SRA1 reaches the floor in both metrics within $100$ NFEs for Problem \ref{prob:ImageNet64} and within $50$ NFEs for Problem \ref{prob:ImageNet512}, whereas neither EM nor EI reaches it within the tested range. Relative to EM, this amounts to a speed-up of at least a factor of four on both problems. As expected from Remark \ref{rem:EIParameterization}, EI performs similarly to EM here, since for $f\equiv0$ it coincides with EM up to the exact integration of $g^2$ and $g$. Chang's method also reaches the floor, but requires two to four times as many steps as SRA1, despite having the same strong order $\gamma=1.5$. Its error constant is thus considerably larger, which is consistent with its higher plateau in Figure \ref{fig:Dissipativity}. In practice, the choice of tableau therefore matters considerably, which underlines the value of the large design space covered by Theorem \ref{lem:SRA}.

In contrast to Problems \ref{prob:CIFAR10} and \ref{prob:LatentCIFAR10}, we refrain from reporting least-squares slopes: only the $\mathcal{W}_2$ curves of Problem \ref{prob:ImageNet512} exhibit a regime in which the predicted orders are approximately visible, whereas the remaining curves do not show a stable slope over the tested range. For the FID, this is expected for the same reason as for Problem \ref{prob:CIFAR10}, namely the mismatch between generation and measurement space. For $\mathcal{W}_2$, the EDM2 setting differs from that of Problems \ref{prob:CIFAR10} and \ref{prob:LatentCIFAR10} in several respects: we use the pretrained networks and the strongly non-uniform time grid without any adaptation, and ImageNet is a considerably more complex distribution than CIFAR-10. Since Problems \ref{prob:ImageNet64} and \ref{prob:ImageNet512} share the same codebase, sampler, and time grid, the fact that the orders become visible only in Problem \ref{prob:ImageNet512}, where both generation and $\mathcal{W}_2$ measurement take place in the latent space, suggests that the regularity of the target distribution plays a role: the latent distribution of the KL-regularized Stable Diffusion autoencoder \cite{rombach2022high} is plausibly more regular than the distribution of ImageNet images in pixel space.

Finally, these results contrast with the observation that out-of-the-box higher-order SDE solvers do not pay off in practice \cite{jolicoeur2021gotta, karras2022elucidating}: without any adaptation of the EDM2 codebase, its noise schedule, or its time discretization, SRA1 is at least four times faster than EM on both problems. This confirms, on a considerably larger benchmark, the expectation formulated for Problem \ref{prob:CIFAR10}, namely that the steeper convergence rate of the higher-order methods translates into a genuine compute-time advantage.

\section*{Code Availability}
The code used in the numerical experiments is available at
\url{https://github.com/emanuelpfarr/SSGDM}.

\section*{Declarations}
Throughout this work generative AI was used to polish formulation and correct grammar. 

\section*{Acknowledgments}
This work was partially supported by the Alexander von Humboldt foundation.

\appendix
\section{Proof of Theorem \ref{rem:GeneralStrongSchemeDefi}}
\begin{proof}
    We want to show that the remainders of the shape defined in Definition \ref{def:generalstrongscheme} also fulfill this bound. Writing $M_k:=\sum_{n=0}^k R_n^M$ and $S_k:=\sum_{n=0}^k R_n^D$ and using $\|M_k+S_k\|^2\leq 2\|M_k\|^2+2\|S_k\|^2$, it suffices to bound the suprema of $M_k$ and $S_k$ separately. Since $K=T/h$ we have
    \[
        \mathbb E\Big[\sup_{0\leq k\leq K-1}\|S_k\|^2\Big]
        \leq K\,\mathbb E\Big[\sum_{n=0}^{K-1}\|R_n^D\|^2\Big]
        \leq K^2\max_{0\leq k\leq K-1}\mathbb E\big[\|R_k^D\|^2\big]
        \leq T^2 C_R\,h^{2\gamma}.
    \]
    For $M_k$, set $\mathcal G_k:=\mathcal F_{t_{k+1}}$. Since $\BackwardNormalstartingNetworkDiscretizedGeneralStrong_n$ is $\mathcal F_{t_n}$-measurable and $R_n^M$ is a function of $\BackwardNormalstartingNetworkDiscretizedGeneralStrong_n$ together with the driving Wiener path on $[t_n,t_{n+1}]$, each $R_n^M$, and hence $M_k=\sum_{n=0}^kR_n^M$, is $\mathcal F_{t_{k+1}}$-measurable, i.e.\ $(M_k)$ is $(\mathcal G_k)$-adapted. Moreover, $R_{k+1}^M$ is by definition centered with respect to $\mathcal F_{t_{k+1}}=\mathcal G_k$, so
    \[
        \mathbb E[M_{k+1}\mid \mathcal G_k] = M_k + \mathbb E[R_{k+1}^M\mid \mathcal F_{t_{k+1}}] = M_k,
    \]
    i.e.\ $(M_k)$ is a discrete-time martingale with respect to $(\mathcal G_k)$. Since, for $i<j$, $R_i^M$ is $\mathcal F_{t_{i+1}}$-measurable and $\mathcal F_{t_{i+1}}\subseteq\mathcal F_{t_j}$, the tower property together with the centering property $\mathbb E[R_j^M\mid\mathcal F_{t_j}]=0$ gives
    \[
        \mathbb E\big[\langle R_i^M,R_j^M\rangle\big]
        = \mathbb E\Big[\big\langle R_i^M,\mathbb E[R_j^M\mid\mathcal F_{t_j}]\big\rangle\Big] = 0
        \qquad (i<j),
    \]
    so all cross terms in $\|M_{K-1}\|^2=\sum_{i,j=0}^{K-1}\langle R_i^M,R_j^M\rangle$ vanish in expectation, and
    \[
        \mathbb E\big[\|M_{K-1}\|^2\big]=\sum_{n=0}^{K-1}\mathbb E\big[\|R_n^M\|^2\big]\leq K\,C_R\,h^{2\gamma+1}.
    \]
    Doob's $L^2$ maximal inequality then yields
    \[
        \mathbb E\Big[\sup_{0\leq k\leq K-1}\|M_k\|^2\Big]
        \leq 4\,\mathbb E\big[\|M_{K-1}\|^2\big]
        \leq 4K\,C_R\,h^{2\gamma+1}
        \leq 4\,C_R\,T\,h^{2\gamma},
    \]
    again using $K=T/h$. Combining the two bounds,
    \[
        \mathbb E\Bigg[\sup_{0\leq k \leq K-1}\Big\|\sum_{n=0}^k R_n\Big\|^2\Bigg]
        \leq 2\,T^2C_R h^{2\gamma}+8\,C_RTh^{2\gamma} = Ch^{2\gamma},
    \]
    which recovers the bound of Kloeden and Platen.
\end{proof}
\section{Proof of Theorems \ref{thm:Main} and \ref{thm:MainDissipative}}
\label{sec:ProofOfMainTheorem}
\label{sec:proof}
\subsection{Proof Idea}
The proof works by using $p_0 = \mathcal{L}\left(\Backward_T\right)$ and splitting the error in all three error types using the triangle inequality
\begin{align*}
\mathcal{W}_2\left(\mathcal{L}\left(\BackwardNormalstartingNetworkDiscretized_K\right), p_0\right) \leq& \mathcal{W}_2\left(\mathcal{L}\left(\BackwardNormalstarting_T\right), \mathcal{L}\left(\Backward_T\right)\right)+\mathcal{W}_2\left(\mathcal{L}\left(\BackwardNormalstartingNetwork_T\right), \mathcal{L}\left(\BackwardNormalstarting_T\right)\right)+ \mathcal{W}_2\left(\mathcal{L}\left(\BackwardNormalstartingNetworkDiscretized_{K}\right), \mathcal{L}\left(\BackwardNormalstartingNetwork_T\right)\right).\end{align*}
The initialization error is bounded by Assumption~\ref{ass:data}. It remains to estimate the score-matching and discretization errors. We treat these terms separately and use the inequality $\mathcal{W}_2(\mathcal{L}(X),\mathcal{L}(Y))\leq \|X-Y\|_{L_2}$ for arbitrary random variables $X$ and $Y$. Usually, we phrase bounds in terms of $\mathbb{E}[\|\cdot\|^2]$, and take the square root at the very end to achieve $L_2$ bounds. Throughout all proofs, we take care to keep the signs explicit, so we can insert dissipativity later on. 

To bound the score matching error, we first need to transfer the error from our assumption, which is taken on the trajectory of the true backward process $\Backward_t$, to the trajectory of $\BackwardNormalstarting_t$. Then, we can use their shared driving process (cf. Discussion around \eqref{eq:BackwardNormalstartingScoreSDE}) to bound the $L_2$ difference of $\BackwardNormalstarting_t$ and $\BackwardNormalstartingNetwork_t$, keeping the scalar products induced from the cross terms of $\|\BackwardNormalstarting_t-\BackwardNormalstartingNetwork_t\|^2$ explicit, so we can insert dissipativity later on. In the general version, we estimate both cross terms using Young's inequality, and conclude by a standard Gronwall argument. 

To bound the discretization errors, we rely on Lemmas that are close to standard facts about strong solutions of SDEs, multiple Itô-integrals and discretization methods. However, we deliberately change the style of proof in most of the cases, to get a better $T$ dependency, which is suppressed in most of the standard literature and is thus not handled with great care, and, perhaps more importantly, to insert the dissipative case easily, without needing to rewrite the whole proof. Specifically, after citing some known results from Kloeden and Platen \cite{kloeden2010numerical}, we reprove the known fact that the second moment of a strong solution of an SDE stays bounded, but keeping the induced cross terms resulting from the multidimensional Itô formula explicit. We go on to prove two basic facts about the expectation of a measurable random variable and Itô integrals, which we use in the two main lemmas to bound the cross terms resulting from the difference $\|\BackwardNormalstartingNetwork_{t_k}-\BackwardNormalstartingNetworkDiscretized_k\|^2$. The remaining terms of the differences are bounded by Young's inequality, and basic facts about multiple Itô integrals. After bounding the error of a single iterate at time step $t_k$, we get the final error by using a discrete Gronwall type argument. 

With this setup, the dissipative case only relies on inserting the condition at the prepared places.

\subsection{General Case}
Throughout this subsection, we impose Assumptions~\ref{ass:Solution}-\ref{ass:DriftAndDiffusion}. To relax the notation, we suppress the $T$ argument in the constants $M_g=M_g(T)$, $ C_\mathcal{A}=C_\mathcal{A}(T)$, $L_s=L_s(T)$, and $L_{s_\theta}=L_{s_\theta}(T)$ which, as noted in the assumptions in Section \ref{sec:assumptions}, depend at most polynomial on $T$. We continue the suppression of $T$ in the constants defined in this section, which will also depend polynomial on $T$, as they are combinations of $M_g$ and $C_\mathcal{A}$. Thus, we define $C$ (without index) to be a generic constant independent of $h$ and $\varepsilon$, and growing at most polynomial in $T$ which may vary from line to line. Whenever a Constant is not depended on $T$, we explicitly say so. Furthermore, without loss of generality, we assume that $h\leq 1$. Since $h=T/K$, this holds automatically once $K\geq T$, i.e. once the discretization is at least as fine as one step per unit of time. This assumption is used repeatedly below, whenever a higher power of $h$ is bounded by a lower one. Furthermore, throughout the remainder of the paper we simplify notation by writing
\[
  C_{\mathcal{A}_\gamma} := C_{\mathcal{A}_\gamma \cup \mathcal{B}(\mathcal{A}_\gamma )}
\]
for the Lipschitz constants from Assumption~\ref{ass:CoefficientFunctions} associated with the
index sets $\mathcal{A}$ and $\mathcal{B}(\mathcal{A})$ of the Itô--Taylor and general strong
schemes (Definitions~\ref{def:ItôTaylorScheme} and~\ref{def:generalstrongscheme}). Because we use it repeatedly throughout, we restate Young's inequality:
\begin{rem}
    For an arbitrary real number $\delta >0$, Young's inequality implies 
    \[2\langle x,y\rangle\leq \frac{1}{\delta}\|x\|^2 + \delta \|y\|^2.\]
\end{rem}

We begin by showing a lemma which allows us to bound the propagated $L_2$ error introduced by the score network. For this purpose, we first need to quantify how the score matching error along the real backward process $\Backward_t$, which is the real training objective, is related to the score matching error along the normal starting trajectory $\BackwardNormalstarting_t$.
\begin{lem}[Transfer of the score-matching bound]
\label{lem:ScoreErrorTransfer}
Under our imposed Assumptions,
\[
    \sup_{t\in[0,T]}\mathbb E\big[\|s(T-t,\BackwardNormalstarting_t)-s_\theta(T-t,\BackwardNormalstarting_t)\|^2\big]
    \;\le\; 2\varepsilon^2 + 2\big(L_{s_\theta}+L_s\big)^2 B_{\text{init}}(T)^2 \;=:\; \tilde\varepsilon^2.
\]
\end{lem}
\begin{proof}
Fix $t\in[0,T]$ and let $(X,\bar X)$ be a $\mathcal W_2$-optimal coupling of
$\mathcal L(\Backward_t)$ and $\mathcal L(\BackwardNormalstarting_t)$, which exists since both laws have finite second moments. Since the left-hand side of the claim depends only on $\mathcal L(\BackwardNormalstarting_t)$,
\begin{align*}
    \mathbb E\big[\|(s-s_\theta)(T-t,\BackwardNormalstarting_t)\|^2\big]
    &= \mathbb E\big[\|(s-s_\theta)(T-t,\bar X)\|^2\big]\\
    &\le 2\,\mathbb E\big[\|(s-s_\theta)(T-t,X)\|^2\big]
      + 2(L_{s_\theta}+L_s)^2\,\mathbb E\big[\|X-\bar X\|^2\big]\\
    &\le 2\varepsilon^2 + 2(L_{s_\theta}+L_s)^2\,
      \mathcal W_2^2\big(\mathcal L(\Backward_t),\mathcal L(\BackwardNormalstarting_t)\big)
     \le 2\varepsilon^2 + 2(L_{s_\theta}+L_s)^2B_{\text{init}}(T)^2,
\end{align*}
using $X\sim\mathcal L(\Backward_t)$ with Assumption~\ref{ass:ScoreMatchingErrorBound},
and Assumption~\ref{ass:data}.
\end{proof}
\begin{lem} (Propagated score matching error)
    \label{lem:ScoreMatchingError3}
    Under our imposed assumptions, we can control the score matching error via
    \[
         \Big\|\BackwardNormalstartingNetwork_{T}-\BackwardNormalstarting_{T}\Big\|_{L_2}
        \le C_1e^{\text{Poly}(T)}\varepsilon + C_2B_\text{init}(T)e^{\text{Poly}(T)},
    \]
    with some constants $C_1, C_2>0$, which are defined in the proof and depend on $T$ at most polynomial.
\end{lem}
\begin{proof}
    We set $e_t = \BackwardNormalstartingNetwork_t-\BackwardNormalstarting_t$ and write the drift for $\BackwardNormalstartingNetwork$ shorthand as $\mu(t,x) = f(T-t)x + g^2(T-t)s_\theta(T-t,x)$. We note, the only difference between the drift of $\BackwardNormalstartingNetwork$ and $\BackwardNormalstarting$ is that the score function $s$ in $\BackwardNormalstarting$ is replaced by the score network $s_\theta$. Because we assume the driving noise process is shared, the stochastic part vanishes from the difference. Therefore, we have 
    \[e_t = \int_0^t \big[\mu(u, \BackwardNormalstartingNetwork_u) - \mu(u,\BackwardNormalstarting_u) + g^2(T-u)\big((s_\theta-s)(T-u, \BackwardNormalstarting_u)\big)\big] \d u \]
    where we added and subtracted $s_\theta(T-u, \BackwardNormalstarting_u)$. Differentiating yields
    \[\frac{\partial }{\partial t}e_t = \mu(t, \BackwardNormalstartingNetwork_t) - \mu(t,\BackwardNormalstarting_t) + g^2(T-t)\Big((s_\theta-s)(T-t, \BackwardNormalstarting_t)\Big).\]
    By the linear growth of the coefficient functions and the bounded score matching difference from Lemma~\ref{lem:ScoreErrorTransfer}, together with the uniform second moment bounds of the strong solutions (cf.\ Assumption~\ref{ass:Solution} and Lemma~\ref{lem:MomentEstimateOfStrongSolution}), we have
    \[
        \sup_{t\in[0,T]}\mathbb E\Big[\|e_t\|^2+\big\|\tfrac{\partial}{\partial t}e_t\big\|^2\Big]<\infty,
    \]
    so that, by dominated convergence, $a(t):=\mathbb E[\|e_t\|^2]$ is differentiable with
    \[
        a'(t) = 2\,\mathbb E\Big[\big\langle e_t,\ \tfrac{\partial}{\partial t}e_t\big\rangle\Big].
    \]
    Carrying out the calculation yields
    \begin{equation}
    \label{eq:ScoreDriftInnerProduct}
    \frac{\partial }{\partial t} a(t)=2 \mathbb{E}\big[\langle e_t,  \mu(t,\BackwardNormalstartingNetwork_t) - \mu(t, \BackwardNormalstarting_t)\rangle\big]+ 2 \mathbb{E}\big[\big\langle e_t ,g^2(T-t)\big((s_\theta-s)(T-t,\BackwardNormalstarting_t)\big)\big\rangle\big].\end{equation} 
    We investigate the first and second term separately. We bound the first term using Young's inequality with $\delta=1$ and the Lipschitz assumption from Assumption~\ref{ass:CoefficientFunctions}, which yields 
    \begin{align*} 2\mathbb{E}\big[\langle e_t,  \mu(t,\BackwardNormalstartingNetwork_t) - \mu(t, \BackwardNormalstarting_t)\rangle\big]
    \leq\;& \mathbb{E}[\|e_t\|^2] + \mathbb{E}[ \|\mu(t,\BackwardNormalstartingNetwork_t) - \mu(t, \BackwardNormalstarting_t)\|^2 ]\\
    \leq\;& \underbrace{(1+C_{\mathcal{A}}^2)}_{=:2\Lambda}\mathbb{E}[\|e_t\|^2].
    \end{align*}
    We estimate the second term once more using Young's inequality with $\delta=1$ to get
    \begin{align*}2\mathbb{E}\Big[\Big\langle e_t, g^2(T-t)\Big((s_\theta-s)(T-t, \BackwardNormalstarting_t)\Big)  \Big\rangle\Big]
    \leq\;& \mathbb{E}[\|e_t\|^2] + g^4(T-t)\mathbb{E}[\|(s_\theta-s)(T-t, \BackwardNormalstarting_t)\|^2]\\
    \leq\;& \mathbb{E}[\|e_t\|^2] + M_g^4\tilde \varepsilon^2.
    \end{align*}
    Using again the Definition of $a(t)$, we get 
    \[a'(t)\leq (2\Lambda +1)a(t)+M_g^4\tilde \varepsilon^2.\]
    Using the Gronwall inequality, we get 
    \[a(T)\leq M_g^4T e^{(2\Lambda+1)T}\tilde \varepsilon^2\leq M_g^4e^{(2\Lambda+2)T}\tilde \varepsilon^2 = 2M_g^4e^{(2\Lambda + 2)T}\varepsilon^2 + 2M_g^4(L_{s_\theta}+L_s)^2B_\text{init}(T)^2e^{(2\Lambda + 2)T}.\]
    Taking the square root and using $\sqrt{x+y}\leq \sqrt{x}+\sqrt{y}$ proves the claim with $C_1=\sqrt{2}M_g^2$, $C_2 = \sqrt{2}M_g^2(L_{s_\theta}+L_s)$, all depending at most polynomial on $T$, and $e^{(\Lambda+1)T}=e^{\text{Poly}(T)}$, since $\Lambda$ is polynomial in $T$. 
\end{proof}
Note that in the non dissipative case, the above lemma can be proven in a more compact way by using Minkwoski's integral inequality, but we are choosing the present route to make the comparison between the dissipative and non dissipative case more clear. 

To bound the discretization error, we employ a similar proof strategy to \cite{kloeden2010numerical}. We first prove a lemma that bounds the error of an Itô--Taylor approximation scheme, and then proceed by bounding the error of an \SchemeClassShortHand, where we use the result of the Itô--Taylor scheme. For this purpose, we restate a couple of basic facts about stochastic integrals and strong solutions of SDEs.
\begin{rem}
    Kloeden and Platen \cite{kloeden2010numerical} state the results below for scalar-valued integrands. Since we defined $I_\alpha[h]_{s,t}$ for vector-valued $h = (h^1,\ldots,h^d)^\top$ componentwise, every scalar estimate extends immediately by applying it to each component $h^k$ and summing, using $\|h\|^2 = \sum_{k=1}^d |h^k|^2$.
\end{rem}
\begin{lem} (Bound on multiple Itô integrals)
    \label{lem:SecondMomentEstimate}
    Let $\alpha \neq \nu$ be a multi-index, $g \in \mathcal{H}_\alpha$,
    $h > 0$. Then
    \[
        \mathbb{E}\Bigl[\sup_{t_n \leq s \leq t_{n+1}}
            \|I_\alpha[g(\,\cdot\,)]_{t_n,s}\|^2\Bigr]
        \leq
        4^{l(\alpha)-n(\alpha)}\,h^{l(\alpha)+n(\alpha)}\,
        \frac{1}{l(\alpha)!}\,
        \mathbb{E}\Bigl[\sup_{t_n \leq s \leq t_{n+1}}\|g(s)\|^2\Bigr].
    \]
\end{lem}
\begin{proof}
    For each component $m = 1,\ldots,d$, the scalar version of this estimate is \cite[Lemma~5.7.4]{kloeden2010numerical}. Since $I_\alpha[g]_{t_n,s} = (I_\alpha[g^1]_{t_n,s},\ldots,I_\alpha[g^d]_{t_n,s})^\top$ by definition, summing over components gives the stated vector-valued inequality.
\end{proof}
\begin{lem} (Vanishing expectations of stochastic Itô integrals)
    \label{lem:VanishingFirstMoment}
    Let $\alpha \neq \nu$ be a multi-index with $l(\alpha)\neq n(\alpha)$, $g \in \mathcal{H}_\alpha$,
    and $t_n \leq t \leq t_{n+1}$. Then
    \[
        \mathbb{E}\bigl[I_\alpha[g(\,\cdot\,)]_{t_n,t}
            \mid \mathcal{F}_{t_n}\bigr] = 0
        \quad \text{a.s.}
    \]
\end{lem}
\begin{proof}
    For each component $m$, the scalar result is
    \cite[Lemma~5.7.1]{kloeden2010numerical}. The vector-valued statement
    follows componentwise.
\end{proof}

The following lemma provides a proof for the known fact that strong solutions of SDEs stay bounded, whose method is similar to that of Kloeden and Platen \cite[Theorem 4.5.4]{kloeden2010numerical}, but we highlight the place where dissipativity is inserted later on. 
\begin{lem} (Bounded moments of solutions)
    \label{lem:MomentEstimateOfStrongSolution}
    Let $n\in\N$, and suppose $X_t$ is the solution of an SDE with Lipschitz continuous and linearly growing drift $\mu(t,x)$ with constant $C_0(T)$, and with diffusion coefficient $g(t)I_d$ with $\sup_{0\leq t \leq T}|g(t)|\leq M_g(T)$, where $M_g(T)$ and $C_0(T)$ depend at most polynomial on $T$. Furthermore, we impose $\mathbb{E}[\|X_0\|^{2n}]=C_{d,n}(T)d$, where $C_{d,n}$ grows at most polynomial in $T$.
    \begin{enumerate}[label=(\alph*)]
        \item \label{lem:MomentEstimateOfStrongSolution:a}
        (fixed-time moments) The function $m_n(t):=\mathbb{E}[\|X_t\|^{2n}]$ is continuously differentiable with
        \begin{equation}
            \label{eq:MomentInnerProduct}
            m_n'(t) = 2n\mathbb E [\|X_t\|^{2n-2}\langle X_t,\mu(t,X_t)\rangle] + g^2(t)n(d+2n-2)\mathbb{E}[\|X_t\|^{2n-2}]
            \leq C_0'\big(1+m_n(t)\big),
        \end{equation}
        with $C_0' := 4nC_0+n(d+2n-2)M_g^2$ growing at most polynomial in $T$. Consequently,
        \[
            \sup_{0\le t\le T}m(t)\le (1+C_{d,n}(T)d)\,e^{\text{Poly}(T)}.
        \]
        \item \label{lem:MomentEstimateOfStrongSolution:b}
        (local supremum) For every $k$,
        \[
            \mathbb{E}\Big[\sup_{t_k\le s\le t_{k+1}}\|X_s\|^{2}\Big] \le C_{\text{sup}X}\big(1+m_1(t_k)\big)=\O(e^{\text{Poly}(T)}),
        \]
        with $C_{\text{sup}X}:=3+6C_0^2e^{C_0'}+12dM_g^2=\O(e^{Poly(T)})$, and independent of $k$.
    \end{enumerate}
\end{lem}
\begin{proof}
    \ref{lem:MomentEstimateOfStrongSolution:a} We apply the multi-dimensional Itô formula to $\varphi(x):=\|x\|^{2n}$, for which $D\varphi(x)u=2n\|x\|^{2n-2}\langle x,u\rangle$ and $\sum_{i=1}^d\varphi(x)[e_i,e_i]=2n(d+2n-2)\|x\|^{2n-2}$. Inserting these calculations yields
    \[
        \d\|X_t\|^{2n} = 2n\|X_t\|^{2n-2}\langle X_t,\mu(t,X_t)\rangle\,\d t + g^2(t)n(d+2n-2)\|X_t\|^{2n-2}\d t + g(t)2n\|X_t\|^{2n-2}\langle X_t,dW_t\rangle.
    \] 
    In integrated form,
    \begin{equation}
    \label{eq:MomentIntegrated}
    \begin{array}{ll}
    \|X_t\|^{2n} &= \|X_0\|^{2n} + \int_0^t\Big(2n\|X_u\|^{2n-2}\langle X_u,\mu(u,X_u)\rangle + g^2(u)n(d+2n-2)\|X_u\|^{2n-2}\Big)\d u \\[1ex]
    &+ \underbrace{\int_0^tg(u)2n\|X_u\|^{2n-2}\langle X_u,dW_u\rangle}_{=:M_t},
    \end{array}
    \end{equation}
    where $M_t$ has vanishing expectation. Because of the Lipschitz, the linear growth and the bounded second moment of $X_0$  assumptions, $X_t$ has a strong solution which satisfies $\sup_{u\le T}\mathbb{E}[\|X_u\|^{2n}]<\infty$. Thus, using dominated convergence in \eqref{eq:MomentIntegrated} yields
    \[
        m_n(t) = m_n(0) + \int_0^t\Big(2n\,\mathbb{E}\big[\|X_u\|^{2n-2}\langle X_u,\mu(u,X_u)\rangle\big] + g^2(u)n(d+2n-2)\mathbb{E}[\|X_u\|^{2n-2}]\Big)\d u.
    \]
    The integrand is continuous in $u$, so $m$ is continuously differentiable and the equality in \eqref{eq:MomentInnerProduct} follows. For the inequality, we estimate pointwise, using the Cauchy--Schwarz inequality, the linear growth of $\mu$, $\|x\|^{2n-2}\le1+\|x\|^{2n}$, and $(1+\|x\|)\|x\|^m\leq 2(1+\|x\|^{m+1})$,
    \begin{align*}
        &\|X_t\|^{2n-2}\langle X_t,\mu(t,X_t)\rangle \le \|X_t\|^{2n-2}\|X_t\|\,C_0(1+\|X_t\|)\le 2C_0\big(1+\|X_t\|^{2n}\big),
        \\[1ex]
        &g^2(t)n(d+2n-2)\|X_t\|^{2n-2}\le n(d+2n-2)M_g^2\|X_t\|^{2n-2} \leq n(d+2n-2)M_g^2\big(1+\|X_t\|^{2n}\big),
    \end{align*}
    and take expectations. The supremum statement then follows from Gronwall's inequality: $m_n'(t)\le C_0'(1+m_n(t))$ implies $1+m_n(t)\le(1+m_n(0))e^{C_0't}$, and since the right-hand side is increasing in $t$ and $m_n(0)=C_{d,n}d$,
    \[
        \sup_{0\le t\le T}m_n(t)\le\sup_{0\le t\le T}(1+m_n(t))\le(1+C_{d,n}d)\,e^{C_0'T}.
    \]

    \ref{lem:MomentEstimateOfStrongSolution:b} Writing $X_s=X_{t_k}+\int_{t_k}^s\mu(u,X_u)\,\d u+\int_{t_k}^sg(u)\,\d W_u$ for $s\in[t_k,t_{k+1}]$ and using $\|a+b+c\|^2\le3(\|a\|^2+\|b\|^2+\|c\|^2)$,
    \[
        \sup_{t_k\le s\le t_{k+1}}\|X_s\|^2
        \le 3\|X_{t_k}\|^2 + 3\sup_{t_k\le s\le t_{k+1}}\Big\|\int_{t_k}^s\mu(u,X_u)\,\d u\Big\|^2 + 3\sup_{t_k\le s\le t_{k+1}}\Big\|\int_{t_k}^sg(u)\,\d W_u\Big\|^2.
    \]
    For the drift integral, the supremum over $s$ is bounded by extending the integration domain and applying the Cauchy--Schwarz inequality on $[t_k,t_{k+1}]$:
    \[
        \sup_{t_k\le s\le t_{k+1}}\Big\|\int_{t_k}^s\mu(u,X_u)\,\d u\Big\|^2
        \le \Big(\int_{t_k}^{t_{k+1}}\|\mu(u,X_u)\|\,\d u\Big)^2
        \le h\int_{t_k}^{t_{k+1}}\|\mu(u,X_u)\|^2\,\d u.
    \]
    For the stochastic integral, the supremum of the squared norm of a martingale is controlled by Doob's $L^2$ maximal inequality followed by the Itô isometry:
    \[
        \mathbb{E}\Big[\sup_{t_k\le s\le t_{k+1}}\Big\|\int_{t_k}^sg(u)\,\d W_u\Big\|^2\Big]
        \le 4\,\mathbb{E}\Big[\Big\|\int_{t_k}^{t_{k+1}}g(u)\,\d W_u\Big\|^2\Big]
        = 4d\int_{t_k}^{t_{k+1}}g^2(u)\,\d u
        \le 4dM_g^2h.
    \]
    Taking expectations in the first display, inserting these two estimates, and using the linear growth of $\mu$, i.e.\ $\mathbb{E}[\|\mu(u,X_u)\|^2]\le 2C_0^2(1+m(u))$, we obtain
    \begin{equation}
    \label{eq:BoundedMoments(b)LastStep}
        \mathbb{E}\Big[\sup_{t_k\le s\le t_{k+1}}\|X_s\|^2\Big]
        \le 3\,m_1(t_k) + 6C_0^2h^2\Big(1+\sup_{t_k\le u\le t_{k+1}}m_1(u)\Big) + 12dM_g^2h.
    \end{equation}
    Finally, applying the Gronwall estimate of part~\ref{lem:MomentEstimateOfStrongSolution:a} on the interval $[t_k,t_{k+1}]$ with initial value $m_n(t_k)$ gives $\sup_{t_k\le u\le t_{k+1}}m_1(u)\le e^{C_0'h}(1+m_1(t_k))$, and inserting this together with $h\le1$ yields the claim with $C_{\text{sup}X}:=3+6C_0^2e^{C_0'}+12dM_g^2$.
\end{proof}

We note that Assumption \ref{ass:CoefficientFunctions} ensures that the coefficient functions fulfill $f_\alpha \in \mathcal H_\alpha$, and that we can apply Lemma \ref{lem:SecondMomentEstimate} to the SDE $\BackwardNormalstartingNetwork_t$. Furthermore, Assumption \ref{ass:CoefficientFunctions} ensures the existence of the Lipschitz constants in Lemma \ref{lem:MomentEstimateOfStrongSolution}, and the definition of the SDE \eqref{eq:BackwardNormalstartingScoreSDE}, which starts at a Normal with known variance, ensures the second moment bound condition on the initial value. We proceed by proving two lemmas that control the expectation of a scalar product that includes stochastic integrals. The first shows that, crucially, due to Lemma \ref{lem:VanishingFirstMoment}, all genuinely stochastic terms vanish. The second controls the surviving, purely deterministic terms in the case of path-dependent integrands.
\begin{lem} (Vanishing of stochastic scalar products)
    \label{lem:StochasticTermsVanish}
    Let $e_n$ be an $\mathcal F_{t_n}$-measurable random variable and let $\alpha$ be a multi-index with $l(\alpha)\neq n(\alpha)$. Then, for $g\in\mathcal H_\alpha$,
    \[\mathbb{E}\Big[\langle e_n,I_\alpha[g(\cdot)]_{t_n,t_{n+1}}\rangle\Big]=0.\]
\end{lem}
\begin{proof}
    Using the tower property of conditional expectation, and the fact that $e_n$ is $\mathcal F_{t_n}$-measurable, we get
    \begin{align*}
        \mathbb{E}\Big[\langle e_n,I_\alpha[g(\cdot)]_{t_n,t_{n+1}}\rangle\Big] &= \mathbb{E}\Big[\mathbb{E}[\langle e_n,I_\alpha[g(\cdot)]_{t_n,t_{n+1}}\rangle|\mathcal F_{t_n}]\Big]\\
        &=\mathbb{E}\Big[\langle e_n,\mathbb{E}[I_\alpha[g(\cdot)]_{t_n,t_{n+1}}|\mathcal F_{t_n}]\rangle\Big]=0,
    \end{align*}
    according to Lemma \ref{lem:VanishingFirstMoment}, which applies since $l(\alpha)\neq n(\alpha)$ means $\alpha$ contains at least one stochastic integrator.
\end{proof}
\begin{lem} (Control of deterministic scalar products)
    \label{lem:DeterministicCrossTerms}
    Let $e_n$ be an $\mathcal F_{t_n}$-measurable random variable and let $\alpha$ be a multi-index with $l(\alpha)= n(\alpha)$. Then, for $g\in\mathcal H_\alpha$,
    \[2\mathbb{E}\Big[\langle e_n,I_\alpha[g(\cdot)]_{t_n,t_{n+1}}\rangle\Big] \leq h \mathbb{E}[\|e_n\|^2] +\frac{1}{l(\alpha)!}h^{2l(\alpha)-1}\mathbb{E}\Big[\sup_{t_n\leq s\leq t_{n+1}}\|g(s)\|^2\Big].\]
\end{lem}
\begin{proof}
    Using Young's inequality with $\delta= 1/h$ and Lemma~\ref{lem:SecondMomentEstimate}, noting that $4^{l(\alpha)-n(\alpha)}=1$ since $l(\alpha)=n(\alpha)$, we compute
    \begin{align*}
        2\mathbb{E}\Big[\langle e_n,I_\alpha[g(\cdot)]_{t_n,t_{n+1}}\rangle\Big] 
        &\leq  h \mathbb{E}[\|e_n\|^2] + \frac1h \mathbb{E}[\|I_\alpha[g(\cdot)]_{t_n,t_{n+1}}\|^2]\\
        &\leq h \mathbb{E}[\|e_n\|^2] +\frac{1}{l(\alpha)!}h^{2l(\alpha)-1}\mathbb{E}\Big[\sup_{t_n\leq s\leq t_{n+1}}\|g(s)\|^2\Big].
    \end{align*}
\end{proof}
For frozen integrands, i.e.\ $g(\cdot)\equiv u$ for an $\mathcal F_{t_n}$-measurable random variable $u$, we will not need Lemma \ref{lem:DeterministicCrossTerms}: in this case the deterministic multiple integrals can be computed exactly, $I_\alpha[u]_{t_n,t_{n+1}}=u\,\frac{h^{l(\alpha)}}{l(\alpha)!}$ for $l(\alpha)=n(\alpha)$, and the corresponding scalar products can be estimated directly. This sharper, direct treatment is carried out inside the proofs of Lemmas \ref{lem:ItôTaylorApprox} and \ref{lem:GeneralStrongScheme} below. Now, we state one of two main lemmas of this section.
\begin{lem} (Discretization error of Itô--Taylor schemes)
    \label{lem:ItôTaylorApprox}
    Let $\BackwardNormalstartingNetworkDiscretizedItoTaylor_0,\BackwardNormalstartingNetworkDiscretizedItoTaylor_1,...,\BackwardNormalstartingNetworkDiscretizedItoTaylor_K$ be the result of a strong Itô-Taylor approximation scheme for the SDE \eqref{eq:BackwardNormalstartingScoreSDE} with order $\gamma$ on time-steps $0=t_0<t_1<...<t_K=T$ (cf. Definition \ref{def:ItôTaylorScheme}). Then, we have
    \[\|\BackwardNormalstartingNetwork_{T}-\BackwardNormalstartingNetworkDiscretizedItoTaylor_{K}\|_{L_2}\leq Ce^{\text{Poly}(T)}h^\gamma.\]
\end{lem}
The proof works in the following way: We first estimate the squared Euclidean error induced in a single step by expanding the norm square of the difference of the scheme and the Itô--Taylor expansion. We estimate each contributing factor by using the Lipschitz and linear growth conditions of the coefficient functions $f_\alpha$, using the second moment bound of the stochastic integrals from Lemma \ref{lem:SecondMomentEstimate}, the second moment bound of a strong solution at time $t$ from Lemma \ref{lem:MomentEstimateOfStrongSolution}, as well as the bound on the scalar product of a stochastic integral from Lemmas \ref{lem:StochasticTermsVanish}  and \ref{lem:DeterministicCrossTerms}. Then, using a discrete Gronwall type argument, we can control the error over the whole time interval. 
\begin{proof}
    We make use of the Itô-Taylor expansion
    \begin{align*}
    \BackwardNormalstartingNetwork_{t_{k+1}}=&\BackwardNormalstartingNetwork_{t_k}+\sum_{\alpha \in \mathcal A_\gamma \backslash \{\nu\}}I_\alpha[f_\alpha(t_k,\BackwardNormalstartingNetwork_{t_k})]_{t_k,t_{k+1}}+\sum_{\alpha \in \mathcal B(\mathcal A_\gamma)}I_\alpha[f_\alpha(\cdot,\BackwardNormalstartingNetwork_{\cdot})]_{t_k,t_{k+1}}.
    \end{align*}
    Setting $e_n=\BackwardNormalstartingNetwork_{t_n}-\BackwardNormalstartingNetworkDiscretizedItoTaylor_n$,
    we obtain
    \begin{align*}
    e_{k+1}=&e_k+\underbrace{\sum_{\alpha \in \mathcal A_\gamma\backslash\{\nu\}}I_\alpha[f_\alpha(t_k,\BackwardNormalstartingNetwork_{t_k})-f_\alpha(t_k,\BackwardNormalstartingNetworkDiscretizedItoTaylor_k)]_{t_k,t_{k+1}}}_{D_k:=}+\underbrace{\sum_{\alpha \in \mathcal B(\mathcal A_\gamma)}I_\alpha[f_\alpha(\cdot,\BackwardNormalstartingNetwork_{\cdot})]_{t_k,t_{k+1}}}_{V_k:=}.
    \end{align*}
    Expanding the square yields
    \begin{align*}
    \mathbb E[\|e_{k+1}\|^2]=&\mathbb E[\|e_k\|^2]+\mathbb E[\|D_k\|^2]+\mathbb E[\|V_k\|^2]+2\mathbb E[\langle e_k,D_k\rangle]+2\mathbb E[\langle e_k,V_k\rangle]+2\mathbb E[\langle D_k,V_k\rangle].
    \end{align*}
    Before bounding each term, we summarize the argument:
    \begin{description}
        \item[{$\mathbb E[\|D_k\|^2]$.}] Bounded by $O(h)\,\mathbb E[\|e_k\|^2]$, via the Lipschitz continuity of $f_\alpha$ (Assumption~\ref{ass:CoefficientFunctions}) and Lemma~\ref{lem:SecondMomentEstimate}.
        \item[{$\mathbb E[\|V_k\|^2]$.}] Bounded by $O(h^{2\gamma+1})$, via the linear growth of $f_\alpha$, Lemma~\ref{lem:SecondMomentEstimate}, and the moment bound of Lemma~\ref{lem:MomentEstimateOfStrongSolution}.
        \item[{$2\mathbb E[\langle e_k,D_k\rangle]$ and $2\mathbb E[\langle e_k,V_k\rangle]$.}] By Lemma~\ref{lem:StochasticTermsVanish}, every multi-index with $l(\alpha)\neq n(\alpha)$ contributes exactly zero. For $D_k$, whose integrands are frozen at $t_k$, the surviving deterministic terms are computed exactly and estimated directly. For $V_k$, whose integrands are path-dependent, the surviving terms are estimated via Lemma~\ref{lem:DeterministicCrossTerms}. This leaves an $O(h)\,\mathbb E[\|e_k\|^2]+O(h^{2\gamma+1})$ contribution.
        \item[{$2\mathbb E[\langle D_k,V_k\rangle]$.}] Bounded by $\mathbb E[\|D_k\|^2]+\mathbb E[\|V_k\|^2]$ via Young's inequality, contributing nothing beyond the first two bounds.
    \end{description}
    Combining all five bounds yields a one-step recursion $a_{k+1}\le(1+Ch)a_k+Ch^{2\gamma+1}e^{C_UT}$ for $a_k:=\mathbb E[\|e_k\|^2]$, to which we then apply a discrete Gronwall argument. In detail:
    \begin{align*}
    \mathbb E[\|D_k\|^2]\leq^{\footnotemark} |\mathcal A_\gamma\backslash\{\nu\}|\sum_{\alpha \in \mathcal A_\gamma\backslash\{\nu\}}F_k^\alpha,
    \end{align*}
    \footnotetext{In principle it is possible to improve upon this estimate: for $\alpha,\beta\in\mathcal A_\gamma\setminus\{\nu\}$ whose nonzero components correspond to \emph{different} Wiener directions, conditioning on $\mathcal F_{t_k}$ and using the independence of the Wiener components gives $\mathbb E[\langle I_\alpha[u]_{t_k,t_{k+1}},I_\beta[v]_{t_k,t_{k+1}}\rangle]=0$. Expanding $\mathbb E[\|D_k\|^2]$ as a double sum and discarding these vanishing cross terms would replace the prefactor $|\mathcal A_\gamma\setminus\{\nu\}|$ by the maximal number of multi-indices sharing a single Wiener direction, without affecting the asymptotic rate $h^\gamma$. However, we omit this refinement, as it only improves upon the specific constant, not the rate.}
    where
    \begin{align*}
    F_k^\alpha=&\mathbb E\Big[\|I_\alpha[f_\alpha(t_k,\BackwardNormalstartingNetwork_{t_k})-f_\alpha(t_k,\BackwardNormalstartingNetworkDiscretizedItoTaylor_k)]_{t_k,t_{k+1}}\|^2\Big]\\
    \leq&4^{l(\alpha)-n(\alpha)}h^{l(\alpha)+n(\alpha)}\frac{1}{l(\alpha)!}\mathbb E\Big[\sup_{t_k\leq s\leq t_{k+1}}\|f_\alpha(t_k,\BackwardNormalstartingNetwork_{t_k})-f_\alpha(t_k,\BackwardNormalstartingNetworkDiscretizedItoTaylor_k)\|^2\Big]\\
    =&4^{l(\alpha)-n(\alpha)}h^{l(\alpha)+n(\alpha)}\frac{1}{l(\alpha)!}\mathbb E\Big[\|f_\alpha(t_k,\BackwardNormalstartingNetwork_{t_k})-f_\alpha(t_k,\BackwardNormalstartingNetworkDiscretizedItoTaylor_k)\|^2\Big]\\
    \leq&4^{l(\alpha)-n(\alpha)}h^{l(\alpha)+n(\alpha)}\frac{1}{l(\alpha)!}C_{\mathcal{A}_\gamma}^2\mathbb E[\|\BackwardNormalstartingNetwork_{t_k}-\BackwardNormalstartingNetworkDiscretizedItoTaylor_k\|^2]\\
    \leq&Ch\mathbb E[\|e_k\|^2].
    \end{align*}
    Above, the supremum over $s$ in the second line becomes trivial in the third, since the integrand does not depend on $s$, only on the fixed time $t_k$. We then used the Lipschitz continuity of $f_\alpha$ and the fact that $l(\alpha)+n(\alpha)\geq1$ for $\alpha\in\mathcal A_\gamma\backslash\{\nu\}$, together with $h\le1$. Therefore,
    \begin{align*}
    \mathbb E[\|D_k\|^2]\leq Ch\mathbb E[\|e_k\|^2].
    \end{align*}
    We can estimate the second additional term with
    \begin{align*}
    \mathbb E[\|V_k\|^2]\leq |\mathcal B(\mathcal A_\gamma)|\sum_{\alpha \in \mathcal B(\mathcal A_\gamma)}U_k^\alpha,
    \end{align*}
    where
    \begin{align*}
    U_k^\alpha=&\mathbb E\Big[\|I_\alpha[f_\alpha(\cdot,\BackwardNormalstartingNetwork_{\cdot})]_{t_k,t_{k+1}}\|^2\Big]\\
    \leq&4^{l(\alpha)-n(\alpha)}h^{l(\alpha)+n(\alpha)}\frac{1}{l(\alpha)!}\mathbb E\Big[\sup_{t_k\leq s\leq t_{k+1}}\|f_\alpha(s,\BackwardNormalstartingNetwork_s)\|^2\Big]\\
    \leq&4^{l(\alpha)-n(\alpha)}h^{l(\alpha)+n(\alpha)}\frac{1}{l(\alpha)!}C_{\mathcal{A}_\gamma}^2\mathbb E\Big[\sup_{t_k\leq s\leq t_{k+1}}(1+\|\BackwardNormalstartingNetwork_s\|^2)\Big]\\
    \leq&2Ch^{l(\alpha)+n(\alpha)}\mathbb E[1+\sup_{t_k\leq s\leq t_{k+1}}\|\BackwardNormalstartingNetwork_s\|^2]\\
    \leq&Ch^{l(\alpha)+n(\alpha)}e^{C_UT},
    \end{align*}
    with $C_{U}$ growing at most polynomial in $T$ and stemming from Lemma \ref{lem:MomentEstimateOfStrongSolution}\ref{lem:MomentEstimateOfStrongSolution:b}. Above, we again used Lemma \ref{lem:SecondMomentEstimate}, this time genuinely needing the $s$-dependent supremum, since $f_\alpha(\cdot,\BackwardNormalstartingNetwork_\cdot)$ is evaluated along the path rather than frozen at $t_k$, and the linear growth condition of $f_\alpha$. Since $l(\alpha)\geq \gamma+1$ when $l(\alpha)=n(\alpha)$ and $l(\alpha)+n(\alpha)\geq 2\gamma+1$ for $\alpha\in\mathcal B(\mathcal A_\gamma)$, we obtain
    \begin{equation}
        \label{eq:UnalphaEstimate}
    U_k^\alpha\leq Ch^{2\gamma+1}e^{C_{U}T}.
    \end{equation}
    Hence,
    \begin{align*}
    \mathbb E[\|V_k\|^2]\leq Ch^{2\gamma+1}e^{C_{U}T}.
    \end{align*}
    We proceed by estimating the terms that include scalar products, beginning with $2\mathbb E[\langle e_k,D_k\rangle]$. By Lemma \ref{lem:StochasticTermsVanish}, all summands with $l(\alpha)\neq n(\alpha)$ vanish. The remaining multi-indices in $\mathcal A_\gamma\backslash\{\nu\}$ with $l(\alpha)=n(\alpha)$ are the purely deterministic ones, i.e.\ $\alpha\in\{(0),(0,0),\ldots\}$ up to length $\lfloor\gamma+0.5\rfloor$, and since the integrands of $D_k$ are frozen at $t_k$, the corresponding integrals can be computed exactly: writing $u_{\alpha,k}:=f_\alpha(t_k,\BackwardNormalstartingNetwork_{t_k})-f_\alpha(t_k,\BackwardNormalstartingNetworkDiscretizedItoTaylor_k)$, we have $I_\alpha[u_{\alpha,k}]_{t_k,t_{k+1}}=u_{\alpha,k}\frac{h^{l(\alpha)}}{l(\alpha)!}$, and hence
    \begin{align*}
    2\mathbb E[\langle e_k,D_k\rangle] &= \sum_{\substack{\alpha \in \mathcal A_\gamma\backslash\{\nu\}\\ l(\alpha)=n(\alpha)}}\frac{2h^{l(\alpha)}}{l(\alpha)!}\,\mathbb E\big[\langle e_k,u_{\alpha,k}\rangle\big]\\
    &\leq \sum_{\substack{\alpha \in \mathcal A_\gamma\backslash\{\nu\}\\ l(\alpha)=n(\alpha)}}\frac{2h^{l(\alpha)}}{l(\alpha)!}\,C_{\mathcal A_\gamma}\,\mathbb E[\|e_k\|^2]
    \leq C h \mathbb E[\|e_k\|^2],
    \end{align*}
    where we used the Cauchy--Schwarz inequality, the Lipschitz continuity of $f_\alpha$, and $l(\alpha)\geq1$ together with $h\leq1$. We single out the dominant term of this sum for later reference: for $\alpha=(0)$, the contribution is exactly
    \begin{equation}
        \label{eq:DriftCrossTerm}
        2h\mathbb{E}[\langle e_k, u_{(0),k}\rangle]=2h\,\mathbb E\big[\big\langle e_k,\ \mu(t_k,\BackwardNormalstartingNetwork_{t_k})-\mu(t_k,\BackwardNormalstartingNetworkDiscretizedItoTaylor_k)\big\rangle\big],
    \end{equation}
    which we have estimated by its worst case $2C_{\mathcal A_\gamma}h\,\mathbb E[\|e_k\|^2]$, discarding the sign of the inner product. For $2\mathbb E[\langle e_k,V_k\rangle]$, the integrands $f_\alpha(\cdot,\BackwardNormalstartingNetwork_\cdot)$ are path-dependent, so the deterministic integrals cannot be computed exactly. Instead, applying Lemma \ref{lem:StochasticTermsVanish} to the indices with $l(\alpha)\neq n(\alpha)$ and Lemma \ref{lem:DeterministicCrossTerms} to the remaining ones yields
    \begin{align*}
        2\mathbb E[\langle e_k,V_k\rangle]&\leq |\mathcal B(\mathcal A_\gamma)_{l=n}|\,h \,\mathbb E[\|e_k\|^2] + \sum_{\alpha \in \mathcal B(\mathcal A_\gamma)_{l=n}}\frac{h^{2l(\alpha)-1}}{l(\alpha)!}\,\mathbb E\Big[\sup_{t_k\leq s\leq t_{k+1}}\|f_\alpha(s,\BackwardNormalstartingNetwork_s)\|^2\Big]\\
        &\leq  Ch \mathbb E[\|e_k\|^2]+ Ch^{2\gamma +1}e^{C_{U}T}.
    \end{align*}
    To get the last inequality above we use the same argument as in the derivation of \eqref{eq:UnalphaEstimate}, particularly that $l(\alpha)\geq \gamma+1$ when $l(\alpha)=n(\alpha)$ for $\alpha\in\mathcal B(\mathcal A_\gamma)$.
    
    Finally, we use Young's inequality with $\delta=1$ to get
    \begin{align*}
    2\mathbb E[\langle D_k,V_k\rangle]\leq \mathbb E[\|D_k\|^2]+\mathbb E[\|V_k\|^2].
    \end{align*}
    Hence,
    \begin{align*}
    \mathbb E[\|e_{k+1}\|^2]\leq& \mathbb E[\|e_{k}\|^2] + \underbrace{Ch\mathbb E[\|e_{k}\|^2]}_{\mathbb E[\|D_k\|^2]} + \underbrace{Ch^{2\gamma +1}e^{C_{U}T}}_{\mathbb E[\|V_k\|^2]}+\underbrace{Ch\mathbb E[\|e_{k}\|^2]}_{2\mathbb E[\langle e_k, D_k\rangle]}\\
    &+ \underbrace{Ch\mathbb E[\|e_{k}\|^2] + Ch^{2\gamma+1}e^{C_{U}T}}_{2\mathbb E[\langle e_k,V_k\rangle]} + \underbrace{Ch\mathbb E[\|e_{k}\|^2] + Ch^{2\gamma +1}e^{C_{U}T}}_{2\mathbb E[\langle D_k,V_k\rangle]}
    \end{align*}
    which yields
    \begin{align*}
    \mathbb E[\|e_{k+1}\|^2]\leq (1+Ch)\mathbb E[\|e_k\|^2]+Ch^{2\gamma+1}e^{C_{U}T}.
    \end{align*}
    Setting $a_k=\mathbb E[\|e_k\|^2]$, we obtain
    \begin{align*}
    a_{k+1}\leq (1+Ch)a_k+Ch^{2\gamma+1}e^{C_{U}T}.
    \end{align*}
    Since $a_0=0$, a discrete Gronwall argument yields
    \begin{align*}
    a_K\leq Ch^{2\gamma+1}e^{C_{U}T}\sum_{j=0}^{K-1}(1+Ch)^j\leq Ch^{2\gamma+1}e^{C_{U}T}\frac{(1+Ch)^K-1}{Ch}.
    \end{align*}
    Using $(1+Ch)^K\leq e^{CKh}=e^{CT}$, we obtain
    \begin{align*}
    a_K\leq Ch^{2\gamma}e^{(C_{U}+C)T}.
    \end{align*}
    Thus, defining $C_\text{exp}:=\frac{1}{2}(C_U+C)$, taking square root and inserting $\sqrt{a_K}=\sqrt{\mathbb E[\|e_K\|^2]}$ yields
    \begin{align*}
    \|\BackwardNormalstartingNetwork_{T}-\BackwardNormalstartingNetworkDiscretizedItoTaylor_K\|_{L_2}\leq Ch^\gamma e^{C_\text{exp}T}.
    \end{align*}
    The claim now follows from the fact that both $C$ and $C_\text{exp}$ are polynomial in $T$.\end{proof}

Below we state our lemma for the error induced by \SchemeClassShortHand. The proof relies on first estimating the difference of a Itô--Taylor scheme and a strong scheme, and then use Lemma \ref{lem:ItôTaylorApprox} to finally estimate the $L_2$ error of a strong scheme. Controlling the difference to the Itô--Taylor schemes uses very similar arguments as the proof of Lemma \ref{lem:ItôTaylorApprox}. Note, that we are assuming the constant of the remainder in the \SchemeClassShortHand $C_R$ to be growing at most exponentially with $T$, i.e. $C_R(T)=\O(e^{\text{Poly(T)}})$.
\begin{lem} (Discretization error of \SchemeClassShortHand)
    \label{lem:GeneralStrongScheme}
    Let $\BackwardNormalstartingNetworkDiscretizedGeneralStrong_0,\BackwardNormalstartingNetworkDiscretizedGeneralStrong_1,...,\BackwardNormalstartingNetworkDiscretizedGeneralStrong_K$ be the result of an \SchemeClassShortHand with order $\gamma$ on time-steps $0=t_0<t_1<...<t_K=T$ (cf. Definition \ref{def:generalstrongscheme}). Then, we have
    \[\|\BackwardNormalstartingNetwork_{T}-\BackwardNormalstartingNetworkDiscretizedGeneralStrong_{K}\|_{L_2}\leq Ce^{\text{Poly}(T)} h^\gamma.\]
\end{lem}
\begin{proof}
    We first prove that the difference between a strong Itô--Taylor scheme $\BackwardNormalstartingNetworkDiscretizedItoTaylor$ and an \SchemeClassShortHand satisfies
    \[\|\BackwardNormalstartingNetworkDiscretizedItoTaylor_K - \BackwardNormalstartingNetworkDiscretizedGeneralStrong_K\|_{L_2}\leq Ce^{\text{Poly}(T)} h^\gamma.\]
     For this purpose, we analyze the difference at timestep $t_{k+1}$ in a similar manner to the proof of Lemma \ref{lem:ItôTaylorApprox}. For convenience, we restate the update rule of the Itô--Taylor Scheme and the \SchemeClassShortHand:
     \begin{align*}
         \BackwardNormalstartingNetworkDiscretizedItoTaylor_{k+1}&=\BackwardNormalstartingNetworkDiscretizedItoTaylor_k + \sum_{\alpha \in \mathcal A_\gamma \backslash \{\nu \}}I_{\alpha}[f_\alpha(t_k, \BackwardNormalstartingNetworkDiscretizedItoTaylor_k)]_{t_k, t_{k+1}}\\
         \BackwardNormalstartingNetworkDiscretizedGeneralStrong_{k+1}&= \BackwardNormalstartingNetworkDiscretizedGeneralStrong_k + \sum_{\alpha \in \mathcal A_\gamma \backslash \{\nu \}} I_\alpha[g_{\alpha,k}]_{t_k, t_{k+1}}+R_k.
     \end{align*}
     We define $e_k :=\BackwardNormalstartingNetworkDiscretizedItoTaylor_{k} - \BackwardNormalstartingNetworkDiscretizedGeneralStrong_{k}$, $D_k =\sum_{\alpha \in \mathcal A_\gamma\backslash\{\nu\}}I_\alpha[f_\alpha(t_k,\BackwardNormalstartingNetworkDiscretizedItoTaylor_{k})-g_{\alpha,k}]_{t_k,t_{k+1}} $. Then we have
     \begin{align*}
         \mathbb{E}[\|e_{k+1}\|^2]=& \mathbb{E}[\|e_k\|^2]+\mathbb{E}[\|D_k\|^2]+\mathbb{E}[\|R_k\|^2]+2\mathbb{E}[\langle e_k, D_k\rangle]+2\mathbb{E}[\langle e_k,R_k\rangle]+ 2\mathbb{E}[\langle D_k,R_k\rangle]
     \end{align*}
     The argument mirrors that of Lemma~\ref{lem:ItôTaylorApprox}, with $V_k$ replaced by $R_k$. We highlight only the places the argument genuinely differs:
    \begin{description}
        \item[{$\mathbb E[\|D_k\|^2]$.}] As in Lemma~\ref{lem:ItôTaylorApprox}, but now picking up an additional $O(h^{2\gamma+1})$ term from the approximation condition on $g_{\alpha,k}$ in Definition~\ref{def:generalstrongscheme}.
        \item[{$\mathbb E[\|R_k\|^2]$.}] Bounded directly by the two-part remainder condition of Definition~\ref{def:generalstrongscheme}.
        \item[{$2\mathbb E[\langle e_k,D_k\rangle]$.}] As in Lemma~\ref{lem:ItôTaylorApprox}, again with the same additional $O(h^{2\gamma+1})$ term.
        \item[{$2\mathbb E[\langle e_k,R_k\rangle]$.}] Splits via $R_k=R_k^M+R_k^D$: the $R_k^M$-part vanishes \emph{exactly}, since $e_k$ is $\mathcal F_{t_k}$-measurable and $R_k^M$ is centered given $\mathcal F_{t_k}$. The $R_k^D$-part is handled by Young's inequality, exactly as the cross terms in Lemma~\ref{lem:ItôTaylorApprox}.
        \item[{$2\mathbb E[\langle D_k,R_k\rangle]$.}] As in Lemma~\ref{lem:ItôTaylorApprox}, via Young's inequality.
     \end{description}
     We start by estimating $\mathbb{E}[\|D_k\|^2]$. Splitting the sum into
     \[\mathbb{E}[\|D_k\|^2] \leq |\mathcal A_\gamma\backslash \{\nu\}|\sum_{\alpha \in \mathcal A_\gamma\backslash\{\nu\}}\mathbb{E}[\|I_\alpha[f_\alpha(t_k,\BackwardNormalstartingNetworkDiscretizedItoTaylor_{k})-g_{\alpha,k}]_{t_k,t_{k+1}}\|^2]\]
     and using Lemma \ref{lem:SecondMomentEstimate}, we get
     \begin{align*}
         &\mathbb{E}[\|I_\alpha[f_\alpha(t_k,\BackwardNormalstartingNetworkDiscretizedItoTaylor_{k})-g_{\alpha,k}]_{t_k,t_{k+1}}\|^2]\\
         \leq \;&4^{l(\alpha)-n(\alpha)}h^{l(\alpha)+n(\alpha)}\frac{1}{l(\alpha)!}\mathbb{E}\Big[\Big\|f_\alpha(t_k,\BackwardNormalstartingNetworkDiscretizedItoTaylor_{k})-g_{\alpha,k}\Big\|^2\Big]\\
         = \;&4^{l(\alpha)-n(\alpha)}h^{l(\alpha)+n(\alpha)}\frac{1}{l(\alpha)!}\mathbb{E}\Big[\Big\|f_\alpha(t_k,\BackwardNormalstartingNetworkDiscretizedItoTaylor_{k})- f_\alpha (t_k, \BackwardNormalstartingNetworkDiscretizedGeneralStrong_k) +f_\alpha (t_k, \BackwardNormalstartingNetworkDiscretizedGeneralStrong_k)-g_{\alpha,k}\Big\|^2\Big]\\
         \leq &\; C h^{l(\alpha)+n(\alpha)}\mathbb{E}[\|e_k\|^2]+Ch^{l(\alpha)+n(\alpha)+2\gamma-\phi(\alpha)}\\
         \leq &\;C h\mathbb{E}[\|e_k\|^2]+Ch^{2\gamma +1},
     \end{align*}
     where the first summand uses the Lipschitz continuity of $f_\alpha$ as before, and the second uses the approximation condition of Definition~\ref{def:generalstrongscheme}, $\mathbb E[\max_k\|g_{\alpha,k}-f_\alpha(t_k,\BackwardNormalstartingNetworkDiscretizedGeneralStrong_k)\|^2]\le Ch^{2\gamma-\phi(\alpha)}$. The final line uses $l(\alpha)+n(\alpha)\geq 1$ for all $\alpha \in \mathcal A_\gamma\backslash \{\nu\}$ (as before) together with $\phi(\alpha)\le l(\alpha)+n(\alpha)-1$, which holds with equality when $l(\alpha)\ne n(\alpha)$ and strictly when $l(\alpha)=n(\alpha)$. Hence $l(\alpha)+n(\alpha)+2\gamma-\phi(\alpha)\ge 2\gamma+1$, and $h\le1$ lets us bound the (possibly larger) exponent by $h^{2\gamma+1}$. Together, we get
     \[\mathbb{E}[\|D_k\|^2] \leq  Ch\mathbb{E}[\|e_k\|^2]+Ch^{2\gamma +1}.\]
     We now estimate $\mathbb E[\|R_k\|^2]$ using the decomposition $R_k=R_k^M+R_k^D$ from Definition~\ref{def:generalstrongscheme}. By $\|R_k^M+R_k^D\|^2\leq 2\|R_k^M\|^2+2\|R_k^D\|^2$ and $h\leq 1$,
     \[
         \mathbb{E}[\|R_k\|^2] \leq 2C_Rh^{2\gamma+1}+2C_Rh^{2\gamma+2}\leq 4C_Rh^{2\gamma+1}.
     \]
    For the term $2\mathbb{E}[\langle e_k, D_k\rangle]$, Lemma \ref{lem:StochasticTermsVanish} removes all summands with $l(\alpha)\neq n(\alpha)$, and since the integrands $f_\alpha(t_k,\BackwardNormalstartingNetworkDiscretizedItoTaylor_k)-g_{\alpha,k}$ are $\mathcal F_{t_k}$-measurable, the surviving deterministic integrals are again exact: with $w_{\alpha,k}:=f_\alpha(t_k,\BackwardNormalstartingNetworkDiscretizedItoTaylor_{k})-g_{\alpha,k}$,
     \begin{align*}
         2\mathbb{E}[\langle e_k, D_k\rangle] &= \sum_{\substack{\alpha \in \mathcal A_\gamma\backslash\{\nu\}\\ l(\alpha)=n(\alpha)}}\frac{2h^{l(\alpha)}}{l(\alpha)!}\,\mathbb E\big[\langle e_k,w_{\alpha,k}\rangle\big]\\
         &\leq \sum_{\substack{\alpha \in \mathcal A_\gamma\backslash\{\nu\}\\ l(\alpha)=n(\alpha)}}\frac{h^{l(\alpha)}}{l(\alpha)!}\Big(\frac{l(\alpha)!}{h^{l(\alpha)-1}}\,\mathbb E[\|e_k\|^2] + \frac{h^{l(\alpha)-1}}{l(\alpha)!}\,\mathbb E[\|w_{\alpha,k}\|^2]\Big)\\
         &\leq \sum_{\substack{\alpha \in \mathcal A_\gamma\backslash\{\nu\}\\ l(\alpha)=n(\alpha)}}\Big(h\,\mathbb E[\|e_k\|^2]
   + \tfrac{C_{\mathcal A_\gamma}^2h^{2l(\alpha)-1}}{(l(\alpha)!)^2}\,\mathbb E[\|e_k\|^2] + \tfrac{2h^{2l(\alpha)-1}}{(l(\alpha)!)^2}\,\mathbb E[\|f_\alpha(t_k,\BackwardNormalstartingNetworkDiscretizedGeneralStrong_k)-g_{\alpha,k}\|^2]\Big)\\
         &= Ch \mathbb{E}[\|e_k\|^2] + Ch^{2\gamma +1},
     \end{align*}
     where we used Young's inequality on each summand with $\delta_\alpha=\frac{h^{l(\alpha)-1}}{l(\alpha)!}$, split $w_{\alpha,k}$ into the Lipschitz difference $f_\alpha(t_k,\BackwardNormalstartingNetworkDiscretizedItoTaylor_k)-f_\alpha(t_k,\BackwardNormalstartingNetworkDiscretizedGeneralStrong_k)$ and the tolerance term $f_\alpha(t_k,\BackwardNormalstartingNetworkDiscretizedGeneralStrong_k)-g_{\alpha,k}$ as in the estimate of $\mathbb E[\|D_k\|^2]$ above, and used $l(\alpha)\geq1$, $h\leq1$, and $2l(\alpha)-1+2\gamma-\phi(\alpha)=2\gamma+1$ for the deterministic indices. As in Lemma \ref{lem:ItôTaylorApprox}, we want to single out the first summand of the sum for later reference:
     \begin{equation}
         \label{eq:DriftCrossTermGeneralStrong}
         2h\mathbb{E}[\langle e_k, w_{(0) ,k}\rangle]=2h\Big(\mathbb{E}[\langle e_k, \mu(t_k, \BackwardNormalstartingNetworkDiscretizedItoTaylor_k)-\mu(t_k,\BackwardNormalstartingNetworkDiscretizedGeneralStrong_k) \rangle] + \mathbb{E}[\langle e_k, \mu(t_k, \BackwardNormalstartingNetworkDiscretizedGeneralStrong_k)-g_{(0),k} \rangle]\Big).
     \end{equation}
    We estimate the term $2\mathbb{E}[\langle e_k,R_k\rangle]$ by again splitting $R_k=R_k^M+R_k^D$. Since $\BackwardNormalstartingNetworkDiscretizedItoTaylor_k$ and $\BackwardNormalstartingNetworkDiscretizedGeneralStrong_k$, and hence $e_k$, are $\mathcal F_{t_k}$-measurable, the tower property together with $\mathbb E[R_k^M\mid\mathcal F_{t_k}]=0$ gives
     \[
         2\mathbb E[\langle e_k,R_k^M\rangle] = 2\mathbb E\big[\langle e_k,\mathbb E[R_k^M\mid\mathcal F_{t_k}]\rangle\big]=0.
     \]
     For the remaining part we use Young's inequality with $\delta=1/h$ as before,
     \[
         2\mathbb E[\langle e_k,R_k^D\rangle] \leq h \mathbb{E}[\|e_k\|^2] + \frac1h\mathbb{E}[\|R_k^D\|^2]\leq h \mathbb{E}[\|e_k\|^2] + C_R h^{2\gamma +1},
     \]
     so that altogether
     \[
         2\mathbb{E}[\langle e_k,R_k\rangle] \leq h \mathbb{E}[\|e_k\|^2] + Ch^{2\gamma +1}.
     \]
     A similar argument, now using the estimate on $\mathbb E[\|R_k\|^2]$ derived above, yields
     \[2\mathbb{E}[\langle D_k,R_k\rangle] \leq \mathbb{E}[\|D_k\|^2]+\mathbb{E}[\|R_k\|^2]\le Ch\mathbb{E}[\|e_k\|^2]+Ch^{2\gamma+1}.\]
     Putting everything together, we arrive at
     \begin{align*}
         \mathbb{E}[\|e_{k+1}\|^2]\leq& \mathbb{E}[\|e_k\|^2] + \underbrace{Ch\mathbb{E}[\|e_k\|^2]+Ch^{2\gamma +1}}_{\mathbb{E}[\|D_k\|^2]} + \underbrace{4C_Rh^{2\gamma+1}}_{\mathbb{E}[\|R_k\|^2]}+ \underbrace{Ch \mathbb{E}[\|e_k\|^2] + Ch^{2\gamma +1}}_{2\mathbb{E}[\langle e_k, D_k\rangle]} \\
         &\qquad+ \underbrace{h \mathbb{E}[\|e_k\|^2] + C h^{2\gamma +1}}_{2 \mathbb{E}[\langle e_k, R_k\rangle]}+ \underbrace{Ch\mathbb{E}[\|e_k\|^2]+Ch^{2\gamma +1}}_{2\mathbb{E}[\langle D_k, R_k\rangle ]}\\
         \leq & (1+Ch) \mathbb{E}[\|e_k\|^2] + Ce^{\text{Poly}(T)}h^{2\gamma +1}.
     \end{align*}
     The $e^{\text{Poly}(T)}$ in the second term factor stems from the constants $C$ and $C_R$ of Definition \ref{def:generalstrongscheme}. Using the same discrete Gronwall argument as in the proof of Lemma \ref{lem:ItôTaylorApprox} and taking the square root, we get 
     \begin{align*}
         \|\BackwardNormalstartingNetworkDiscretizedItoTaylor_K - \BackwardNormalstartingNetworkDiscretizedGeneralStrong_K\|_{L_2}\leq e^{\text{Poly}(T)}Ch^\gamma.
     \end{align*}
     With this result and Lemma \ref{lem:ItôTaylorApprox}, we get
     \begin{align*}
         \|\BackwardNormalstartingNetwork_{T}-\BackwardNormalstartingNetworkDiscretizedGeneralStrong_K\|_{L_2}\leq \|\BackwardNormalstartingNetwork_{T}-\BackwardNormalstartingNetworkDiscretizedItoTaylor_K\|_{L_2} + \|\BackwardNormalstartingNetworkDiscretizedItoTaylor_K - \BackwardNormalstartingNetworkDiscretizedGeneralStrong_K\|_{L_2}\leq Ce^{\text{Poly}(T)}h^\gamma.
     \end{align*}
\end{proof}

Now, we can prove the main theorem in the general case:
\begin{proof}[\textbf{Proof of Theorem \ref{thm:Main}.}]
    First we note that $\text{Poly}(T)e^{\text{Poly}(T)}=\O(e^{\text{Poly}(T)})$. Then, splitting the error into
    \begin{align*}
    \mathcal{W}_2\left(\mathcal{L}\left(\BackwardNormalstartingNetworkDiscretized_K\right), p_0\right) \leq& \mathcal{W}_2\left(\mathcal{L}\left(\BackwardNormalstarting_T\right), \mathcal{L}\left(\Backward_T\right)\right)+\mathcal{W}_2\left(\mathcal{L}\left(\BackwardNormalstartingNetwork_T\right), \mathcal{L}\left(\BackwardNormalstarting_T\right)\right)+ \mathcal{W}_2\left(\mathcal{L}\left(\BackwardNormalstartingNetworkDiscretized_{K}\right), \mathcal{L}\left(\BackwardNormalstartingNetwork_T\right)\right),\end{align*}
    and using Lemma \ref{lem:ScoreMatchingError3}, as well as Lemma \ref{lem:GeneralStrongScheme} and the assumption on the bounded initialization error, we get the claim.
\end{proof}

\begin{rem}[Relation to Milstein's fundamental theorem]
\label{rem:MilsteinRelation}
Milstein's fundamental theorem \cite[Section~1.1, Thm.~1.1]{milstein2004stochastic} infers global mean-square order $p_2-\tfrac12$ from two bounds on the one-step error, i.e.\ on the deviation of one step of the scheme from the exact solution started at the same point: a bound of order $h^{p_1}$ on its mean and of order $h^{p_2}$ on its mean-square norm, where $p_2\ge\tfrac12$ and $p_1\ge p_2+\tfrac12$. In our analysis, the role of the exact solution is taken by the Itô--Taylor scheme of Definition~\ref{def:ItôTaylorScheme}: in the proof of Lemma~\ref{lem:GeneralStrongScheme}, we compare an \SchemeClassShortHand with the Itô--Taylor scheme rather than with the exact solution. The remainder conditions of Definition~\ref{def:generalstrongscheme} are chosen precisely such that the remainders satisfy Milstein's bounds with $(p_1,p_2)=(\gamma+1,\gamma+\tfrac12)$, i.e.
\[
    \max_{0\le k\le K-1}\big\|\mathbb E[R_k\mid\mathcal F_{t_k}]\big\|_{L_2}\le C h^{\gamma+1},
    \qquad
    \max_{0\le k\le K-1}\|R_k\|_{L_2}\le C h^{\gamma+1/2},
\]
which follows from $\mathbb E[R_k\mid\mathcal F_{t_k}]=\mathbb E[R_k^D\mid\mathcal F_{t_k}]$ and $h\le1$. The tolerance conditions on the coefficients $g_{\alpha,k}$ ensure the same orders for the remaining part of the one-step difference. 
\end{rem}

\begin{rem}[Alternative decomposition]
\label{rem:AlternativeDecomposition}
The score-matching and initialization errors can also be separated via the process $\BackwardNetwork$ from \eqref{eq:BackwardNetworkScoreSDE} instead of $\BackwardNormalstarting$, i.e.\ by comparing $\BackwardNormalstartingNetwork$ with $\BackwardNetwork$, and $\BackwardNetwork$ with $\Backward$. Assumption~\ref{ass:ScoreMatchingErrorBound} then applies directly, without the transfer of Lemma~\ref{lem:ScoreErrorTransfer}, whereas the initialization error is propagated by the learned instead of the true reverse dynamics, so that its decay relies on the regularity of $s_\theta$, e.g.\ Assumption~\ref{ass:Dissipativity}, rather than on the data distribution. As Assumption~\ref{ass:Dissipativity} is in turn typically inherited from strong log-concavity of the forward marginals, both routes rest on the same structural condition in practice. We use the present decomposition, as it allows existing bounds on the initialization error to be inserted directly.
\end{rem}

\subsection{Dissipative Case}
\label{subsec:Dissipative}
Throughout this subsection we additionally impose Assumption~\ref{ass:Dissipativity}. In particular, we assume the constants $C_R$ and $C$ of Definition \ref{def:generalstrongscheme} to grow at most polynomial in $T$. The proofs are adaptations from the non dissipative case, and we particularly revisit four labeled displays \eqref{eq:ScoreDriftInnerProduct}, \eqref{eq:MomentInnerProduct}, \eqref{eq:DriftCrossTerm}, and \eqref{eq:DriftCrossTermGeneralStrong}. Every estimate not explicitly restated below is taken over from the general case unchanged. We will use the following elementary comparison fact: if $a\colon[0,T]\to[0,\infty)$ is differentiable with $a'(t)\le-\lambda a(t)+b$ for some $b\ge0$, then $(e^{\lambda t}a(t))'\le be^{\lambda t}$, and integrating yields
\begin{equation}
    \label{eq:ODEComparison}
    a(t)\le a(0)e^{-\lambda t}+\tfrac b\lambda\big(1-e^{-\lambda t}\big)\le\max\big\{a(0),\tfrac b\lambda\big\}\qquad\text{for all }t\in[0,T].
\end{equation}

\begin{lem}(Propagated score matching error, dissipative case)
    \label{lem:ScoreMatchingErrorDissipative}
    In the setting of Lemma~\ref{lem:ScoreMatchingError3}, impose additionally Assumption~\ref{ass:Dissipativity}. Then
    \[
        \Big\|\BackwardNormalstartingNetwork_{T}-\BackwardNormalstarting_{T}\Big\|_{L_2}\le \frac{\sqrt2 M_g^2}{\lambda}\Big(\varepsilon+(L_{s_\theta}+L_s)B_\text{init}(T)\Big),
    \]
    where $M_g=M_g(T)$ grows at most polynomial in $T$.
\end{lem}
\begin{proof}
    We restate \eqref{eq:ScoreDriftInnerProduct}:
    \[
        a'(t)=2\,\mathbb{E}\big[\langle e_t,\mu(t,\BackwardNormalstartingNetwork_t)-\mu(t,\BackwardNormalstarting_t)\rangle\big]+2\,\mathbb{E}\big[\langle e_t,g^2(T-t)(s_\theta-s)(T-t,\BackwardNormalstarting_t)\rangle\big].
    \]
    For the first term, Assumption~\ref{ass:Dissipativity} gives directly
    \[
        2\,\mathbb{E}\big[\langle e_t,\mu(t,\BackwardNormalstartingNetwork_t)-\mu(t,\BackwardNormalstarting_t)\rangle\big]\le-2\lambda\,a(t),
    \]
    replacing the worst-case bound $2\Lambda a(t)$. For the second term we repeat the estimate of the general case, but choose the Young parameter as $1/\lambda$ instead of $1$:
    \[
        2\,\mathbb{E}\big[\langle e_t,g^2(T-t)(s_\theta-s)(T-t,\BackwardNormalstarting_t)\rangle\big]\le\lambda\,a(t)+\tfrac{M_g^4}{\lambda}\tilde\varepsilon^2,
    \]
    using Lemma~\ref{lem:ScoreErrorTransfer}. Hence $a'(t)\le-\lambda a(t)+\tfrac{M_g^4}{\lambda}\tilde\varepsilon^2$ with $a(0)=0$, and \eqref{eq:ODEComparison} yields $a(T)\le M_g^4\tilde\varepsilon^2/\lambda^2$. Taking the square root and splitting $\tilde\varepsilon$ as in the general case proves the claim.
\end{proof}

\begin{lem}(Uniform moment bound, dissipative case)
    \label{lem:MomentEstimateDissipative}
    In the setting of Lemma~\ref{lem:MomentEstimateOfStrongSolution}, impose additionally Assumption~\ref{ass:Dissipativity}. Then, for every $n\in\mathbb N$ with $\mathbb{E}[\|X_0\|^{2n}]<\infty$,
    \[
        \sup_{0\le t\le T}m_n(t)\;\le\;M_{\infty,n}(T),
    \]
    where $M_{\infty,0}(T):=1$ and, recursively,
    \[
        M_{\infty,n}(T):=\max\Big\{\mathbb{E}[\|X_0\|^{2n}],\ \tfrac1\lambda\Big(\tfrac{C_{\mathcal A_\gamma}(T)^{2n}}{\lambda^{2n-1}}+n(d+2n-2)M_g(T)^2\,M_{\infty,n-1}\Big)\Big\} \quad \text{growing at most polynomial in $T$. }
    \]
    In particular, $M_{\infty,1}(T)=\max\big\{C_{d,1}(T)d,\ \tfrac1\lambda\big(\tfrac{C_{\mathcal A_\gamma}(T)^2}{\lambda}+dM_g(T)^2\big)\big\}$. Moreover, part~\ref{lem:MomentEstimateOfStrongSolution:b} of Lemma~\ref{lem:MomentEstimateOfStrongSolution} holds with the $T$-polynomial right-hand side $C_{\text{sup,poly}X}(T)(1+M_{\infty,1}(T))$, where $C_{\text{sup,poly}X}(T):=3+6C_\mathcal{A}(T)^2 + 12dM_g(T)^2$.
\end{lem}
\begin{proof}
    In the proof, we suppress $T$ dependence as all constants grow at most polynomial with $T$. We restate \eqref{eq:MomentInnerProduct}, with the convention $m_0\equiv1$:
    \[
        m_n'(t) = 2n\,\mathbb{E}\big[\|X_t\|^{2n-2}\langle X_t,\mu(t,X_t)\rangle\big] + n(d+2n-2)\,g^2(t)\,m_{n-1}(t).
    \]
    Instead of the worst-case linear-growth estimate, we now insert Assumption~\ref{ass:Dissipativity}: adding and subtracting $\mu(t,0)$ and using the linear growth of $f_{(0)}=\mu$ at $x=0$, i.e.\ $\|\mu(t,0)\|\le C_{\mathcal A_\gamma}$, we obtain for all $x\in\R^d$
    \[
        \langle x,\mu(t,x)\rangle=\langle x-0,\mu(t,x)-\mu(t,0)\rangle+\langle x,\mu(t,0)\rangle\le-\lambda\|x\|^2+C_{\mathcal A_\gamma}\|x\|.
    \]
    Multiplying by $2n\|x\|^{2n-2}\ge0$, inserting $x=X_t$, taking expectations and using $g^2\le M_g^2$ gives
    \begin{equation}
        \label{eq:MomentDissipativeRaw}
        m_n'(t)\le-2n\lambda\,m_n(t)+2nC_{\mathcal A_\gamma}\,\mathbb{E}[\|X_t\|^{2n-1}]+n(d+2n-2)M_g^2\,m_{n-1}(t).
    \end{equation}
    To absorb the middle term, we use the scalar Young inequality $ab\le\frac{a^r}{r}+\frac{b^s}{s}$ with $r=\frac{2n}{2n-1}$ and $s=2n$. For $y\ge0$ and $\delta>0$,
    \[
        y^{2n-1}=\big(\delta^{\frac1{2n}}y\big)^{2n-1}\,\delta^{-\frac{2n-1}{2n}}
        \le\frac{2n-1}{2n}\,\delta\,y^{2n}+\frac{1}{2n}\,\delta^{-(2n-1)}.
    \]
    Choosing $\delta:=\lambda/C_{\mathcal A_\gamma}$ and multiplying by $2nC_{\mathcal A_\gamma}$ yields
    \[
        2nC_{\mathcal A_\gamma}\,\mathbb{E}[\|X_t\|^{2n-1}]\le(2n-1)\lambda\,m_n(t)+\frac{C_{\mathcal A_\gamma}^{2n}}{\lambda^{2n-1}},
    \]
    and inserting this into \eqref{eq:MomentDissipativeRaw},
    \begin{equation}
        \label{eq:MomentDissipativeRecursion}
        m_n'(t)\le-\lambda\,m_n(t)+\frac{C_{\mathcal A_\gamma}^{2n}}{\lambda^{2n-1}}+n(d+2n-2)M_g^2\,m_{n-1}(t).
    \end{equation}
    Note that $\mathbb{E}[\|X_0\|^{2n}]<\infty$ implies $\mathbb{E}[\|X_0\|^{2k}]<\infty$ for all $k\le n$ by Jensen's inequality, so that Lemma~\ref{lem:MomentEstimateOfStrongSolution}~\ref{lem:MomentEstimateOfStrongSolution:a} applies to every $m_k$, $k\le n$. We now proceed by induction over $n$.

    \emph{Base case $n=1$.} Since $m_0\equiv1$, \eqref{eq:MomentDissipativeRecursion} reads $m_1'\le-\lambda m_1+\frac{C_{\mathcal A_\gamma}^2}{\lambda}+dM_g^2$, and the comparison inequality~\eqref{eq:ODEComparison} together with $m_1(0)=C_{d,1}d$ gives $\sup_{t}m_1(t)\le M_{\infty,1}$.

    \emph{Induction step $n-1\to n$.} Assume $\sup_{t}m_{n-1}(t)\le M_{\infty,n-1}$. Then the right-hand side of \eqref{eq:MomentDissipativeRecursion} is bounded by $-\lambda m_n(t)+b_n$ with the constant
    \[
        b_n:=\frac{C_{\mathcal A_\gamma}^{2n}}{\lambda^{2n-1}}+n(d+2n-2)M_g^2\,M_{\infty,n-1},
    \]
    and \eqref{eq:ODEComparison} together with $m_n(0)=\mathbb{E}[\|X_0\|^{2n}]$ gives $\sup_tm_n(t)\le\max\{\mathbb{E}[\|X_0\|^{2n}],b_n/\lambda\}=M_{\infty,n}$.

    For the second claim, we repeat all the calculations up until \eqref{eq:BoundedMoments(b)LastStep}, but then insert $\sup_um(u)\leq M_{\infty,1}$ directly instead of using a Gronwall argument, which gives the polynomial constant $3+6C_\mathcal{A}(T)^2 + 12dM_g(T)^2$.
\end{proof}

\begin{lem}(Discretization error Itô--Taylor schemes, dissipative case)
    \label{lem:ItôTaylorDissipative}
    In the setting of Lemma~\ref{lem:ItôTaylorApprox}, impose additionally Assumption~\ref{ass:Dissipativity}, and let $h\le \hat h_0(T):=\min\{1,\lambda/(2\hat C(T)), 1/\lambda\}$, where $\hat C$ is the constant defined in the proof. Then
    \[
        \|\BackwardNormalstartingNetwork_{T}-\BackwardNormalstartingNetworkDiscretizedItoTaylor_{K}\|_{L_2}\le \text{Poly}(T)h^\gamma.
    \]
\end{lem}
\begin{proof}
    The proof is an adaptation to the proof of Lemma \ref{lem:ItôTaylorApprox}, and large parts remain the same. We begin at the decomposition 
    \[
        \mathbb E[\|e_{k+1}\|^2]=\mathbb E[\|e_k\|^2]+\mathbb E[\|D_k\|^2]+\mathbb E[\|V_k\|^2]+2\mathbb E[\langle e_k,D_k\rangle]+2\mathbb E[\langle e_k,V_k\rangle]+2\mathbb E[\langle D_k,V_k\rangle].
    \]
    We now estimate each term separately:
    \begin{description}
        \item[{$\mathbb E[\|D_k\|^2]$.}] We diverge from the non-dissipative case only at the estimate 
        \begin{align*}
             F_k^\alpha=&\mathbb E\Big[\|I_\alpha[f_\alpha(t_k,\BackwardNormalstartingNetwork_{t_k})-f_\alpha(t_k,\BackwardNormalstartingNetworkDiscretizedItoTaylor_k)]_{t_k,t_{k+1}}\|^2\Big]
        \leq4^{l(\alpha)-n(\alpha)}h^{l(\alpha)+n(\alpha)}\frac{1}{l(\alpha)!}C_{\mathcal{A}_\gamma}^2\mathbb E[\|\BackwardNormalstartingNetwork_{t_k}-\BackwardNormalstartingNetworkDiscretizedItoTaylor_k\|^2].
        \end{align*}
        Specifically, rather than estimating the term with $l(\alpha)\geq 1$, we notice that, since the diffusion coefficient is spatially independent, $f_{(j)}(t,x)=g(t)e_j$ (cf.\ Section~\ref{sec:ExplicitSchemes} and \eqref{eq:coefficientfunctions1}) does not depend on $x$ in our setting, so the differences vanish exactly for this index. Every remaining $\alpha\in\mathcal A_\gamma\setminus\{\nu\}$ with $\alpha \notin [d]$ satisfies $l(\alpha)+n(\alpha)\ge2$, yielding the improved bound \footnote{We want to note that this step would be possible in the general case, but it only would improve the constant $C_\text{exp}$, not improve the rate, which is why we omitted this step earlier to keep the presentation as simple as possible.}
        \[
            \mathbb E[\|D_k\|^2]\le \hat C_1 h^2\,\mathbb E[\|e_k\|^2].
        \]
        \item[{$\mathbb E[\|V_k\|^2]$.}] We depart from the non-dissipative case at the estimate
        \begin{align*}
        U_k^\alpha=&\mathbb E\Big[\|I_\alpha[f_\alpha(\cdot,\BackwardNormalstartingNetwork_{\cdot})]_{t_k,t_{k+1}}\|^2\Big]
        \leq2Ch^{l(\alpha)+n(\alpha)}\mathbb E[1+\sup_{t_k\leq s\leq t_{k+1}}\|\BackwardNormalstartingNetwork_s\|^2].
        \end{align*}
        Now, rather than using Lemma \ref{lem:MomentEstimateOfStrongSolution}, we use the dissipative version, i.e. Lemma \ref{lem:MomentEstimateDissipative}, to arrive at the bound 
        \begin{align*}
        \mathbb E[\|V_k\|^2]\leq C h^{2\gamma+1}.
        \end{align*}
        \item[{$2\mathbb E[\langle e_k,D_k\rangle]$.}] Departure happens only with the term $\alpha=(0)$, which we singled out in \eqref{eq:DriftCrossTerm}, where we insert Assumption~\ref{ass:Dissipativity}:
        \[
            2h\,\mathbb E\big[\langle e_k,\mu(t_k,\BackwardNormalstartingNetwork_{t_k})-\mu(t_k,\BackwardNormalstartingNetworkDiscretizedItoTaylor_k)\rangle\big]\le -2\lambda h\mathbb E[\|e_k\|^2],
        \]
        while the remaining deterministic indices have $l(\alpha)\ge2$ and contribute at most $\hat C_2 h^2\,\mathbb E[\|e_k\|^2]$.
        
        \item[{$2\mathbb E[\langle e_k,V_k\rangle]$.}] We repeat the estimate of Lemma~\ref{lem:DeterministicCrossTerms} with Young parameter $\delta:=\tfrac{2}{\lambda h}$ in place of $1/h$, and use Lemma~\ref{lem:MomentEstimateDissipative} in place of Lemma~\ref{lem:MomentEstimateOfStrongSolution}, obtaining
        \[
            2\mathbb E[\langle e_k,V_k\rangle]\le\tfrac\lambda2 h\,\mathbb E[\|e_k\|^2]+C_V h^{2\gamma+1},
        \]
        with $C_V$ polynomial in $T$.
        
        \item[{$2\mathbb E[\langle D_k,V_k\rangle]$.}] Stays exactly the same: bounded by $\mathbb E[\|D_k\|^2]+\mathbb E[\|V_k\|^2]$ via Young's inequality, contributing nothing beyond the first two bounds.
    \end{description}
    Setting $\hat C:=2\hat C_1+\hat C_2$ and collecting terms,
    \[
        a_{k+1}\le\big(1-2\lambda h+\tfrac\lambda2h+\hat C h^2\big)a_k+C h^{2\gamma+1}
        \le(1-\lambda h)\,a_k+C h^{2\gamma+1},
    \]
    where $a_k=\mathbb{E}[\|e_k\|^2]$, and the last step uses $\hat C h^2\le\tfrac\lambda2 h$ for $h\le \hat h_0$. Since $a_0=0$ and $0\le1-\lambda h<1$, the geometric series gives, uniformly in $K$,
    \[
        a_K\le C h^{2\gamma+1}\sum_{j\ge0}(1-\lambda h)^j=\frac{C}{\lambda}\,h^{2\gamma},
    \]
    and taking square roots, and respecting the polynomial $T$ dependence of the constants yields the claim.
\end{proof}

\begin{lem}(Discretization error of \SchemeClassShortHand, dissipative case)
    \label{lem:GeneraslStrongSchmeDissipative}
    In the setting of Lemmas~~\ref{lem:GeneralStrongScheme}, impose additionally Assumption~\ref{ass:Dissipativity}, and let $h\leq h_0: = \min\{\hat h_0, \bar h_0\}$, with $\bar h_0(T):=\min\{1,\lambda/(2\bar C(T)), 1/\lambda\}$, where $\bar C$ is the constant defined in the proof. Then
    \[
        \|\BackwardNormalstartingNetwork_{T}-\BackwardNormalstartingNetworkDiscretizedGeneralStrong_{K}\|_{L_2}\le \text{Poly}(T)h^\gamma.
    \]
\end{lem}
\begin{proof}
    Like in the case of the earlier Lemma, the proof is an adaptation to the proof of the non-dissipative Lemma \ref{lem:GeneralStrongScheme}, and large parts remain the same. We begin at the decomposition 
    \begin{align*}
         \mathbb{E}[\|e_{k+1}\|^2]=& \mathbb{E}[\|e_k\|^2]+\mathbb{E}[\|D_k\|^2]+\mathbb{E}[\|R_k\|^2]+2\mathbb{E}[\langle e_k, D_k\rangle]+2\mathbb{E}[\langle e_k,R_k\rangle]+ 2\mathbb{E}[\langle D_k,R_k\rangle]
     \end{align*}
     We highlight only the places the argument genuinely differs:
    \begin{description}
        \item[{$\mathbb E[\|D_k\|^2]$.}] The only thing that differs from the non-dissipative case is that in the bound 
        \begin{align*}
         \mathbb{E}[\|I_\alpha[f_\alpha(t_k,\BackwardNormalstartingNetworkDiscretizedItoTaylor_{k})-g_{\alpha,k}]_{t_k,t_{k+1}}\|^2]\leq  C h^{l(\alpha)+n(\alpha)}\mathbb{E}[\|e_k\|^2]+Ch^{l(\alpha)+n(\alpha)+2\gamma-\phi(\alpha)},
        \end{align*}
        we drop the term corresponding to the multiindex $(j)$ from the Lipschitz difference, as our diffusion coefficient is spatially independent (cf. the bound for $\mathbb E[\|D_k\|^2]$ in the proof of Lemma \ref{lem:ItôTaylorDissipative}), giving us the bound  
        \[\mathbb E[\|D_k\|^2]\leq \bar C_1 h^2\mathbb E[\|e_k\|^2]+Ch^{2\gamma +1}.\]
        
        \item[{$\mathbb E[\|R_k\|^2]$.}] No difference to the non-dissipative case, remaining at 
        \[\mathbb E[\|R_k\|^2]\leq Ch^{2\gamma +1}.\]
        
        \item[{$2\mathbb E[\langle e_k,D_k\rangle]$.}] This term differs from the non-dissipative case more substantially: We again drop all non deterministic summands and define $w_{\alpha,k}:=f_\alpha(t_k,\BackwardNormalstartingNetworkDiscretizedItoTaylor_{k})-g_{\alpha,k}$, as well as $N_\gamma:=\{\alpha\in\mathcal A_\gamma\setminus\{\nu\}:l(\alpha)=n(\alpha), l(\alpha)>1\}$ and $N:=\#N_\gamma$. We decompose 
         \begin{align*}
         2\mathbb{E}[\langle e_k, D_k\rangle] &= 2h\,\mathbb E[\langle e_k,w_{(0),k}\rangle]+\sum_{\substack{\alpha \in N_\gamma}}\frac{2h^{l(\alpha)}}{l(\alpha)!}\,\mathbb E\big[\langle e_k,w_{\alpha,k}\rangle\big],
         \end{align*}
         and estimate the sum on the right hand side as 
         \begin{align*}
             &\sum_{\substack{\alpha \in N_\gamma}}\frac{2h^{l(\alpha)}}{l(\alpha)!}\,\mathbb E\big[\langle e_k,w_{\alpha,k}\rangle\big] \\\leq&\; \sum_{\substack{\alpha \in N_\gamma}}\frac{h^{l(\alpha)}}{l(\alpha)!}\Big(\frac{\lambda l(\alpha)!}{8N h^{l(\alpha)-1}}\,\mathbb E[\|e_k\|^2] + \frac{8N h^{l(\alpha)-1}}{\lambda l(\alpha)!}\,\mathbb E[\|w_{\alpha,k}\|^2]\Big)\\
             \leq&\; \sum_{\substack{\alpha \in N_\gamma}}\Big(\frac{\lambda}{8N}h\,\mathbb E[\|e_k\|^2]
           + \tfrac{8N C_{\mathcal A_\gamma}^2h^{2l(\alpha)-1}}{\lambda(l(\alpha)!)^2}\,\mathbb E[\|e_k\|^2] + \tfrac{8Nh^{2l(\alpha)-1}}{\lambda(l(\alpha)!)^2}\,\mathbb E[\|f_\alpha(t_k,\BackwardNormalstartingNetworkDiscretizedGeneralStrong_k)-g_{\alpha,k}\|^2]\Big)\\
             \leq&\; \frac{\lambda}{8}h \mathbb{E}[\|e_k\|^2] + \bar C_2 h^2\mathbb{E}[\|e_k\|^2]+Ch^{2\gamma +1},
         \end{align*}
        where we used Young's inequality with parameter $\delta_\alpha=\frac{8N h^{l(\alpha)-1}}{\lambda l(\alpha)!}$, and $l(\alpha)>1$. The term missing from the sum, i.e. the term corresponding to the multiindex $\alpha=(0)$ was singled out in \eqref{eq:DriftCrossTermGeneralStrong}, and we estimate it directly:
        \begin{align*}
            2h\,\mathbb E[\langle e_k,w_{(0),k}\rangle]=\;&2h\,\mathbb E\big[\langle e_k,\mu(t_k,\BackwardNormalstartingNetworkDiscretizedItoTaylor_k)-\mu(t_k,\BackwardNormalstartingNetworkDiscretizedGeneralStrong_k)\rangle\big]+2h\,\mathbb E\big[\langle e_k,f_{(0)}(t_k,\BackwardNormalstartingNetworkDiscretizedGeneralStrong_k)-g_{(0),k}\rangle\big]\\
            \le\;&-2\lambda h\,\mathbb E[\|e_k\|^2]+\tfrac{\lambda}{8}h\,\mathbb E[\|e_k\|^2]+\tfrac{8}{\lambda}h\,\mathbb E[\|f_{(0)}(t_k,\BackwardNormalstartingNetworkDiscretizedGeneralStrong_k)-g_{(0),k}\|^2]\\
            \le\;&\big(-2\lambda+\tfrac{\lambda}{8}\big)h\,\mathbb E[\|e_k\|^2]+\tfrac{8C}{\lambda}h^{2\gamma+1},
        \end{align*}
        using Assumption~\ref{ass:Dissipativity} on the first summand and Young's inequality with parameter $8/\lambda$ on the second. Putting the whole sum together yields 
        \[2\mathbb{E}[\langle e_k, D_k\rangle] \leq \left(-2\lambda h +\frac{\lambda}{4}h\right)\mathbb E[\|e_k\|^2]+\bar C_2 h^2\mathbb E[\|e_k\|^2] + Ch^{2\gamma+1}.\]
        \item[{$2\mathbb E[\langle e_k,R_k\rangle]$.}] Splitting $R_k=R_k^M+R_k^D$. The term $\mathbb E[\langle e_k,R_k^M\rangle]=0 $ like in the non-dissipative case and using Young's inequality with parameter $\delta=4/\lambda h$ yields
         \[
             2\mathbb E[\langle e_k,R_k^D\rangle] \leq \frac{\lambda h}{4} \mathbb{E}[\|e_k\|^2] + \frac{4}{\lambda h}\mathbb{E}[\|R_k^D\|^2]\leq\frac{\lambda }{4}h \mathbb{E}[\|e_k\|^2] + C h^{2\gamma +1},
         \]
         so that altogether
         \[
             2\mathbb{E}[\langle e_k,R_k\rangle] \leq \frac{\lambda }{4}h \mathbb{E}[\|e_k\|^2] + Ch^{2\gamma +1}.
         \]

        \item[{$2\mathbb E[\langle D_k,R_k\rangle]$.}] Exactly as in the non-dissipative case.
     \end{description}
     Setting $\bar C:=2\bar C_1+\bar C_2$ and collecting terms,
    \[
        a_{k+1}\le\big(1-2\lambda h+\tfrac\lambda4h+\tfrac\lambda4h+\bar C h^2\big)a_k+C h^{2\gamma+1}
        \le(1-\lambda h)\,a_k+C h^{2\gamma+1},
    \]
    where $a_k=\mathbb{E}[\|e_k\|^2]$, and the last step uses $\bar C h^2\le\tfrac\lambda2 h$ for $h\le \bar h_0$. Since $a_0=0$ and $0\le1-\lambda h<1$, the geometric series gives, uniformly in $K$,
    \[
        a_K\le C h^{2\gamma+1}\sum_{j\ge0}(1-\lambda h)^j=\frac{C}{\lambda}\,h^{2\gamma},
    \]
    and taking square roots, respecting the polynomial $T$ dependence of the constants, yields the claim. $\|\BackwardNormalstartingNetworkDiscretizedItoTaylor_K - \BackwardNormalstartingNetworkDiscretizedGeneralStrong_K\|_{L_2} \leq \text{Poly(T)} h^{\gamma }$. Mirroring the non-dissipative case, combining this result with Lemma \ref{lem:ItôTaylorDissipative} yields the claim.
     
\end{proof}

Now we can prove the main theorem in the dissipative case.
\begin{proof}[\textbf{Proof of Theorem \ref{thm:MainDissipative}.}]
    Splitting the error by the triangle inequality exactly as in the proof of Theorem~\ref{thm:Main}, and using Assumption~\ref{ass:data}, Lemma~\ref{lem:ScoreMatchingErrorDissipative}, as well as Lemma~\ref{lem:GeneraslStrongSchmeDissipative} for the three terms,
    \[
        \mathcal W_2\big(\mathcal L(\BackwardNormalstartingNetworkDiscretized_K),p_0\big)\le \O\left(\text{Poly}(T)\left(B_\text{init}(T)+\varepsilon+h^\gamma\right)\right).
    \]
\end{proof}

\section{Proving the Order of Schemes in Section \ref{sec:ExplicitSchemes}}
\label{sec:ProofsExplicitSchemes} Recall, we consider SDEs with additive noise, in particular with a general vector-valued drift coefficient $\mu(t,x_t)\in\R^d$ and a spatially independent, scalar diffusion coefficient $g(t)$:
\[
\mathrm dx_t = \mu(t,x_t)\,\mathrm dt + g(t)\,\mathrm dW_t, \qquad x_t\in\R^d,
\]
driven by a $d$-dimensional Wiener process $W_t$, so that the diffusion matrix is $g(t)I_d$. Since we are interested in strong schemes up to order $\gamma = 1.5$, we introduce the corresponding multi-index sets
\begin{align*}
[d] &:= \{(1),\ldots,(d)\},\\
\mathcal A_{1/2} &= \{\nu,(0)\}\cup [d],\\
\mathcal A_{1} &= \mathcal A_{1/2}\cup [d]^2,\\
\mathcal A_{3/2} &= \mathcal A_{1} \cup \{(0,0)\}
\cup \bigl({(0)}\times [d]\bigr)
\cup \bigl([d]\times {(0)}\bigr)
\cup [d]^3.
\end{align*}
One key observation is that the diffusion coefficient is independent of the spatial variable. Consequently, all coefficient functions corresponding to multi-indices containing at least two nonzero entries vanish. We want to note that throughout this section, we assume all the assumptions made in Section \ref{sec:assumptions}, in particular Assumption \ref{ass:CoefficientFunctions}, which gives us the existence of the derivatives used below. Furthermore, we impose Assumption \ref{ass:EI}, whenever we talk about the EI scheme, and Assumption \ref{ass:SRA} whenever we talk about the SRA family. 

A straightforward calculation yields
\begin{equation}
\label{eq:coefficientfunctions1}
\begin{array}{rcl@{\qquad}rcl@{\qquad}rcl}
f_{(0)}
&=& \mu(t,x)
&
f_{(j)}
&=& g(t)e_j
&
f_{(j,0)}
&=& g(t)D\mu(t,x)e_j
\\[1ex]
f_{(0,j)}
&=& g'(t)e_j
&
f_{(j_1,j_2)}
&=& 0
&
f_{(j_1,j_2,j_3)}
&=& 0,
\end{array}
\end{equation}
and
\begin{equation}
    \label{eq:coefficientfunctions2}
    f_{(0,0)}(t,x) =  \partial_t\mu(t,x) + D\mu(t,x)\,\mu(t,x) +\frac12 g(t)^2 \sum_{i=1}^d D^2\mu(t,x)[e_i,e_i],
\end{equation}
where $e_j$ denotes the $j$-th standard basis vector of $\R^d$. In particular, the coefficient functions associated with the multi-indices in $\mathcal A_1\setminus \mathcal A_{1/2}$ vanish identically. Therefore, in the present setting the strong order-$1$ and strong order-$1/2$ schemes coincide, a well-known consequence of the spatial independence of the diffusion coefficient, see \cite[discussion below (10.3.10)]{kloeden2010numerical}.

As can be seen at a glance in \eqref{eq:coefficientfunctions1}, we need, additionally to the stochastic integrals defined in \eqref{eq:StochasticIntegrals}, the stochastic integral
\[\sum_{j=1}^d I_{(0,j)}[1]_{t_k,t_{k+1}}\,e_j\in\R^d, \qquad I_{(0,j)}[1]_{t_k,t_{k+1}}=\int_{t_k}^{t_{k+1}}(t-t_k)\,\d W_t^{(j)}\]
for the analysis of the schemes introduced in section \ref{sec:ExplicitSchemes}. To express this integral in terms of $\Delta W_{t_k}^{(j)}$ and $\Delta Z_{t_k}^{(j)}$, fix $j\in\{1,\ldots,d\}$ and apply the Itô formula to the product $f(t)\,Y_t$, where $f(t):=t-t_k$ is deterministic and $Y_t:=W_t^{(j)}-W_{t_k}^{(j)}$. Since $f$ is of finite variation, its quadratic variation vanishes, so the Itô formula for the product reduces to the classical product rule,
\[
    \d\bigl[(t-t_k)\bigl(W_t^{(j)}-W_{t_k}^{(j)}\bigr)\bigr] = (t-t_k)\,\d W_t^{(j)} + \bigl(W_t^{(j)}-W_{t_k}^{(j)}\bigr)\,\d t.
\]
Integrating both sides over $[t_k,t_{k+1}]$, the left-hand side telescopes to its value at $t_{k+1}$ minus its value at $t_k$ (which is $0$), giving
\[
    (t_{k+1}-t_k)\bigl(W_{t_{k+1}}^{(j)}-W_{t_k}^{(j)}\bigr)
    = \int_{t_k}^{t_{k+1}}(t-t_k)\,\d W_t^{(j)} + \int_{t_k}^{t_{k+1}}\bigl(W_t^{(j)}-W_{t_k}^{(j)}\bigr)\,\d t,
\]
that is,
\[
    h\,\Delta W_{t_k}^{(j)} = I_{(0,j)}[1]_{t_k,t_{k+1}} + \Delta Z_{t_k}^{(j)}, \qquad\text{i.e.}\qquad \sum_{j=1}^d I_{(0,j)}[1]_{t_k,t_{k+1}} = h\,\Delta W_{t_k}-\Delta Z_{t_k}.
\]
This identity allows us to express the remaining stochastic integral entirely in terms of the three integrals defined above. Furthermore, as described in \cite[Chapter 7]{chang1987numerical}, we can simulate $\Delta Z$ in an elementary fashion: a direct computation of the covariance structure of iterated Itô integrals gives, for each $j$,
\begin{equation}
    \label{eq:VarofDeltaZ}
    \mathrm{Cov}\big(\Delta W_{t_k}^{(j)},\Delta Z_{t_k}^{(j)}\big) = \tfrac{h^2}{2}, \qquad
    \mathrm{Var}\big(\Delta Z_{t_k}^{(j)}\big) = \tfrac{h^3}{3}.
\end{equation}
Both $\Delta W_{t_k}$ and $\Delta Z_{t_k}$ are linear functionals of the increments $(W_s-W_{t_k})_{s\in[t_k,t_{k+1}]}$. Hence the vector $(\Delta W_{t_k},\Delta Z_{t_k})\in\mathbb R^{2d}$ is jointly Gaussian and independent of $\mathcal F_{t_k}$. Moreover, the pairs $\big(\Delta W_{t_k}^{(j)},\Delta Z_{t_k}^{(j)}\big)$, $j=1,\ldots,d$, are independent, each pair is a functional of the single component $W^{(j)}$ of the underlying Wiener process, and these components are themselves independent. Setting 
\[
    \Gamma_{t_k}^{(j)} := \frac{2\sqrt3}{h^{3/2}}\Big(\Delta Z_{t_k}^{(j)}-\frac h2\,\Delta W_{t_k}^{(j)}\Big), \qquad j=1,\ldots,d,
\]
we obtain the decomposition
\begin{equation}
    \label{eq:DeltaZDecomposition}
    \Delta Z_{t_k}^{(j)} = \tfrac{h}{2}\Delta W_{t_k}^{(j)} + \tfrac{h^{3/2}}{2\sqrt3}\Gamma_{t_k}^{(j)}, \qquad j=1,\ldots,d,
\end{equation}
where $\Gamma_{t_k}=(\Gamma_{t_k}^{(1)},\ldots,\Gamma_{t_k}^{(d)})^\top \sim \mathcal N(0,I_d)$ is a standard Gaussian vector, independent of $\Delta W_{t_k}$ and of $\mathcal F_{t_k}$. With the notation fixed and practical simulation discussed, we can proceed by proving Theorems \ref{lem:EMScheme} to \ref{lem:SRA}.

\subsection{Technical Lemmas}
We proceed by introducing intermediate results that are used to bound the remainders $R_k$ in the proofs of Theorems~\ref{lem:EIScheme} and~\ref{lem:SRA}. These remainders consist of terms including the coefficients and their derivatives evaluated at (or near) one iterate of the scheme, and polynomial terms depending on the iterate. Since these terms grow linearly in the iterate, bounding the remainders requires bounds on the moments of the iterates themselves, which, for Theorem~\ref{thm:MainDissipative}, additionally have to be polynomial in $T$. The following lemmas provide exactly these bounds, using basic stability and consistency reasoning. Notably however, the notion of consistency used here differs from the one in the proofs of the schemes: there, the scheme is compared with the Itô--Taylor scheme along its own iterates, which requires moment bounds of the iterates. Here, one step of the scheme is compared with the exact solution, which only requires moment bounds of the exact solution. Moreover, a local error of order $h$ suffices, since we only need boundedness of the moments and no convergence rate. The Euler--Maruyama scheme has no remainder and does not require these bounds. 

We use the notation $\bar \Forward_{k+1}=\Psi_k(\bar \Forward_k, \zeta_k)$ to denote the update of a one-step scheme depending on the previous state $\bar \Forward_k$ and the random variables $\zeta_k=(\Delta W_{t_k},\Delta Z_{t_k})$. Throughout, $\Forward_t$ denotes the exact solution of the SDE discretized by the scheme, driven by the same Wiener process, and the scheme is started in its initial value, $\bar\Forward_0=\Forward_{t_0}$.
\begin{lem}[Moment bound of one-step schemes]
\label{lem:BoundedMomentsScheme}
Let $n\in \N$ with $\sup_{0\le t\le T}\|\Forward_t\|_{L_{2n}}<\infty$, and suppose there are $h_0\in (0,1]$, $\kappa \in \R$ and $C_\delta>0$ such that for all $h\leq h_0$ and all $k$:
\begin{enumerate}[label=(H\arabic*)]
    \item\label{H1} (pathwise contraction/expansion) $\|\Psi_k(x,\zeta)-\Psi_k(y,\zeta)\|\leq (1-\kappa h)\|x-y\|$ for all $x,y\in \R^d$ and all values $\zeta$ of the noise
    \item\label{H2} (local defect) $\|\Delta_k\|_{L_{2n}}:=\|\Forward_{t_{k+1}}-\Psi_k(\Forward_{t_k},\zeta_{k})\|_{L_{2n}}\leq C_\delta h$.
\end{enumerate}
    Then, for all $h\le h_0$ and with $e_k:= \Forward_{t_k}-\bar \Forward_k$,
    \begin{enumerate}[label=(\alph*)]
        \item \label{lem:BoundedMomentsSchemeNonDissipative} if $\kappa=-C_\Psi\leq 0$: $\sup_k\|e_k\|_{L_{2n}}\leq C_\delta T e^{C_\Psi T} $ and hence $\sup_k\|\bar \Forward_k\|_{L_{2n}}\leq \sup_t\|\Forward_t\|_{L_{2n}}+C_\delta T e^{C_\Psi T}$
        \item \label{lem:BoundedMomentsSchemeDissipative} if $\kappa>0$ and $\kappa h_0\leq 1$: $\sup_k\|e_k\|_{L_{2n}}\leq \tfrac{C_\delta}{\kappa} $ and hence $\sup_k\|\bar \Forward_k\|_{L_{2n}}\leq \sup_t\|\Forward_t\|_{L_{2n}}+\tfrac{C_\delta}{\kappa}$.
    \end{enumerate}
\end{lem}
\begin{proof}
    Let $a_k := \|e_k\|_{L_{2n}}$. By definition of the scheme, $e_{k+1}=\Forward_{t_{k+1}}-\Psi_k(\bar\Forward_k,\zeta_k)$, so the triangle inequality for $\|\cdot\|_{L_{2n}}$ gives
    \begin{equation}
    \label{eq:MomentBoundOneStepSchemesFirstStep}
        e_{k+1}=\Delta_k+\big(\Psi_k(\Forward_{t_k}, \zeta_k)-\Psi_k(\bar \Forward_{k}, \zeta_k)\big)\quad \Rightarrow\quad a_{k+1}\leq \|\Delta_k\|_{L_{2n}}+\|\Psi_k(\Forward_{t_k}, \zeta_k)-\Psi_k(\bar \Forward_{k}, \zeta_k)\|_{L_{2n}}.
    \end{equation}
    In both cases of the lemma $1-\kappa h\ge0$ (for $\kappa>0$ because $\kappa h\le\kappa h_0\le1$). Hence \ref{H1} yields
    \begin{align*}
        \|\Psi_k(\Forward_{t_k}, \zeta_k)-\Psi_k(\bar \Forward_{k}, \zeta_k)\|_{L_{2n}} &= \big(\mathbb{E}[\|\Psi_k(\Forward_{t_k}, \zeta_k)-\Psi_k(\bar \Forward_{k}, \zeta_k)\|^{2n}]\big)^{\frac{1}{2n}}\\
        &\leq \big(\mathbb{E}[((1-\kappa h)\|e_k\|)^{2n}]\big)^{\frac{1}{2n}}
        = (1-\kappa h)\, a_k.
    \end{align*}
    Inserting this and \ref{H2} into the right-hand side of \eqref{eq:MomentBoundOneStepSchemesFirstStep}, we obtain
    \[a_{k+1}\leq (1-\kappa h)a_k + C_{\delta}h.\]
    Since $\bar\Forward_0=\Forward_{t_0}$, we have $a_0=0$, and iterating this recursion (which is allowed since $1-\kappa h\ge0$) gives
    \[
        a_k\le C_\delta h\sum_{j=0}^{k-1}(1-\kappa h)^j\qquad\text{for all }0\le k\le K.
    \]
    The claims now follow from a discrete Gronwall argument:
    \begin{enumerate}[label=(\alph*)]
        \item If $\kappa=-C_\Psi\le0$, each summand satisfies $(1+C_\Psi h)^j\le(1+C_\Psi h)^K\le e^{C_\Psi hK}=e^{C_\Psi T}$, and there are at most $K=T/h$ summands, so
        \[a_k\leq C_\delta h\,K\,e^{C_\Psi T} = C_\delta Te^{C_\Psi T}.\]
        \item If $\kappa >0$ and $\kappa h\leq 1$, then $0\le1-\kappa h<1$, and bounding the finite sum by the geometric series gives
        \[a_k\leq C_\delta h\sum_{j\geq 0}(1-\kappa h)^j = \frac{C_\delta}{\kappa}.\]
    \end{enumerate}
    Finally, $\bar \Forward_k=\Forward_{t_k}-e_k$ and the triangle inequality once more give $\|\bar\Forward_k\|_{L_{2n}}\le\|\Forward_{t_k}\|_{L_{2n}}+a_k$, which proves the bounds on the iterates.
\end{proof}

\begin{lem}[EI and SRA satisfy \ref{H1} and \ref{H2}]
\label{lem:MomentsEISRA}
Let $n\in\mathbb N$ with $\mathbb{E}[\|\Forward_0\|^{2n}]<\infty$, and let $f,g$ be Lipschitz continuous on $[0,T]$ with
polynomial in $T$ constants.
\begin{enumerate}[label=(\alph*)]
    \item \label{lem:MomentsEISRA_EI}The exponential integrator~\eqref{eq:EI_compact} satisfies
    \ref{H1} with $\kappa=\lambda/2$ under~Assumption~\ref{ass:Dissipativity}, and with
    $\kappa=-C_\Psi$ in general, and \ref{H2} with $C_\delta=C\big(1+\sup_{0\le t\le T}\|\Forward_t\|_{L_{2n}}\big)$,    for all $h\le h_0^\text{EI}$, where $h_0^\text{EI}:=\min\{1,1/F,1/\lambda,\lambda/(2C)\}$ under Assumption~\ref{ass:Dissipativity} and $h_0^\text{EI}:=\min\{1,1/F\}$ in general, and the constants $C_\Psi,C$ grow at most polynomially in $T$.
    \item \label{lem:MomentsEISRA_SRA}The same holds for every member of the SRA family~\eqref{eq:SRA}
    with explicit stages (i.e.\ $A^{(0)}$ strictly lower triangular),
    $c^{(0)},c^{(1)}\in[0,1]^s$, ${\beta^{(1)}}^Te=1$, ${\beta^{(2)}}^Te=0$, $\alpha_i\ge0$, $\alpha^Te=1$, and $h\le h_0^\text{SRA}:=\min\{1,1/\lambda,\lambda/(2C)\}$ under Assumption~\ref{ass:Dissipativity} and $h_0^\text{SRA}:=1$ in general.
\end{enumerate}
Consequently, by Lemma~\ref{lem:BoundedMomentsScheme} together with Lemma~\ref{lem:MomentEstimateOfStrongSolution} resp.\ Lemma~\ref{lem:MomentEstimateDissipative}, the iterates of both schemes satisfy
\[
    M_n':=\sup_{0\le k\le K}\mathbb{E}\big[\|\bar\Forward_k\|^{2n}\big]
    \le\begin{cases}
        \big(M_{\infty,n}^{1/(2n)}+2C_\delta/\lambda\big)^{2n} & \text{under Assumption~\ref{ass:Dissipativity}},\\[0.5ex]
        \big((1+\mathbb{E}[\|\Forward_0\|^{2n}])^{1/(2n)}e^{C_n'T/(2n)}+C_\delta Te^{C_\Psi T}\big)^{2n} & \text{in general},
    \end{cases}
\]
for all $h\le h_0^\text{EI}$ or $h\le h_0^\text{SRA}$. Under Assumption~\ref{ass:Dissipativity} the bound grows at most polynomial in $T$.
\end{lem}
\begin{proof}
In both schemes the noise enters additively and does not depend on the
current state, so it cancels in $\Psi_k(x,\zeta)-\Psi_k(y,\zeta)$, i.e.\ \ref{H1} is a statement about the drift only. We write $L_\mu:=C_{\mathcal A_\gamma}$ for the Lipschitz and linear-growth constant of $\mu=f_{(0)}$ from Assumption~\ref{ass:CoefficientFunctions}, and $C_s$ for the linear-growth constant of the score network, $\|s_\theta(\tau,x)\|\le C_s(1+\|x\|)$, which follows from Assumption~\ref{ass:ScoreFunctionLipschitz} together with $\sup_{\tau}\|s_\theta(\tau,0)\|<\infty$. We use the generic
constant $C$ for quantities depending only on $L_\mu,M_g,L_g,L_{s_\theta},C_s,
\sup|f|$, the Lipschitz constant of $f$, the tableau, $d$ and $n$.

\smallskip\noindent\emph{(a) Exponential integrator.}
Here $\Psi_k(x)=\Phi_kx+a_ks_\theta(\tau_k,x)+b_k$ with $b_k$ independent of
$x$. For $u:=x-y$ and $v:=s_\theta(\tau_k,x)-s_\theta(\tau_k,y)$, so that
$\|v\|\le L_{s_\theta}\|u\|$ by Assumption~\ref{ass:ScoreFunctionLipschitz},
\[
    \|\Psi_k(x)-\Psi_k(y)\|^2=\Phi_k^2\|u\|^2+2\Phi_ka_k\langle u,v\rangle+a_k^2\|v\|^2.
\]
The Lipschitz continuity of $f$ and $g$ yields $\Phi_k=1+f(\tau_k)h+r_\Phi$ and
$a_k=g^2(\tau_k)h+r_a$ with $|r_\Phi|,|r_a|\le Ch^2$ (cf. Theorem~\ref{lem:EIScheme}, this uses $h\leq 1/F$), hence
\[
    \|\Psi_k(x)-\Psi_k(y)\|^2
    \le\|u\|^2+2h\underbrace{\big(f(\tau_k)\|u\|^2+g^2(\tau_k)\langle v,u\rangle\big)}_{B:=}+Ch^2\|u\|^2.
\]
By the identity~\eqref{eq:DissipativityIdentity}, the bracket $B$ equals
$B=\langle\mu(t_k,x)-\mu(t_k,y),u\rangle$. Under Assumption~\ref{ass:Dissipativity} we have $B\le-\lambda\|u\|^2$,
so
\[\|\Psi_k(x)-\Psi_k(y)\|^2\le(1-2\lambda h+Ch^2)\|u\|^2\le(1-\lambda h)\|u\|^2\le(1-\tfrac{\lambda}{2}h)^2\|u\|^2\]
for $h\le h_0^\text{EI}:=\min\{1,1/F,1/\lambda,\lambda/(2C)\}$,
which is \ref{H1} with $\kappa=\lambda/2$, and $\kappa h_0\le\frac12$. In general, $B\le L_\mu\|u\|^2$ due to the Lipschitz continuity of $\mu$ from Assumption~\ref{ass:CoefficientFunctions}, so that for $h\le\min\{1,1/F\}$
\[
    \|\Psi_k(x)-\Psi_k(y)\|^2\le\big(1+2L_\mu h+Ch^2\big)\|u\|^2\le\big(1+(L_\mu+\tfrac C2)h\big)^2\|u\|^2,
\]
which is \ref{H1} with $\kappa=-C_\Psi$, $C_\Psi:=L_\mu+\frac C2$. Thus $C_\Psi$ is the Lipschitz rate of the drift plus the contribution of the second-order term of the step.

For \ref{H2}, the exact solution, driven by the same $W$, reads
\[
    \Forward_{t_{k+1}}=\Phi_k\Forward_{t_k}+\int_{t_k}^{t_{k+1}}\Phi_{t_{k+1},s}\,g^2(T-s)\,s_\theta(T-s,\Forward_s)\,\d s+b_k,
\]
so that $\Delta_k=\int_{t_k}^{t_{k+1}}\Phi_{t_{k+1},s}\,g^2(T-s)\,s_\theta(T-s,\Forward_s)\,\d s-a_ks_\theta(\tau_k,\Forward_{t_k})$. Since $\Phi_{t_{k+1},s}\le e^{Fh}\leq e$ for $h\leq h_0$ and $s\in[t_k,t_{k+1}]$, $g^2\le M_g^2$ and $s_\theta$ has linear growth, Minkowski's integral inequality yields
\[
    \|\Delta_k\|_{L_{2n}}
    \le C\,h\Big(1+\sup_{t_k\le s\le t_{k+1}}\|\Forward_s\|_{L_{2n}}\Big)
    \le C\,h\big(1+\sup_{0\le t\le T}\|\Forward_t\|_{L_{2n}}\big),
\]
and the supremum is bounded by Lemma~\ref{lem:MomentEstimateOfStrongSolution}, and polynomial in $T$ by Lemma~\ref{lem:MomentEstimateDissipative} under Assumption~\ref{ass:Dissipativity}. Note that
no time-regularity of $s_\theta$ is needed here, only linear growth.

\smallskip\noindent\emph{(b) SRA family.}
Write $\mu_i(x):=\mu(t_k+c_i^{(0)}h,x)$ and let $H_i(x)$ denote the stage
values~\eqref{eq:SRA_b} started from $x$. The noise contributions
$\eta_i:=\sum_jB^{(0)}_{ij}g(t_k+c^{(1)}_jh)\Delta Z_{t_k}/h$ in the stages and
$\eta_0:=\sum_i(\beta^{(1)}_i\Delta W_{t_k}+\beta^{(2)}_i\Delta Z_{t_k}/h)\,g(t_k+c^{(1)}_ih)$
in~\eqref{eq:SRA_a} do not depend on $x$. Hence, for $u:=x-y$ and
$u_i:=H_i(x)-H_i(y)$,
\[
    u_i=u+h\sum_{j<i}A^{(0)}_{ij}\big(\mu_j(H_j(x))-\mu_j(H_j(y))\big),
    \qquad
    \Psi_k(x)-\Psi_k(y)=u+h\sum_{i=1}^s\alpha_i\big(\mu_i(H_i(x))-\mu_i(H_i(y))\big).
\]
Since
$\alpha_i\ge0$, expanding the squared norm yields
\[
    \|\Psi_k(x)-\Psi_k(y)\|^2 \le\|u\|^2+2h\sum_{i=1}^s\Big\langle u,\alpha_i\big(\mu_i(H_i(x)) - \mu_i(H_i(y))\big)\Big\rangle
      +h^2\Big(\sum_{i=1}^s\alpha_i\big\|\mu_i(H_i(x)) - \mu_i(H_i(y))\big\|\Big)^2.
\]

Since the stages are explicit, $u_1=u$, and a simple induction over the stages together with the Lipschitz continuity of $\mu$ yields $\|u_i-u\|\le C_Hh\|u\|$ and $\|u_i\|\le(1+C_Hh)\|u\|$ with $C_H$
depending only on $L_\mu$ and $A^{(0)}$ (for $h\le1$). Under Assumption~\ref{ass:Dissipativity}, applied
at time $t_k+c_i^{(0)}h\in[0,T]$, we have
\begin{align*}
    \big\langle u,\mu_i(H_i(x))-\mu_i(H_i(y))\big\rangle
    &=\big\langle u_i,\mu_i(H_i(x))-\mu_i(H_i(y))\big\rangle
      +\big\langle u-u_i,\mu_i(H_i(x))-\mu_i(H_i(y))\big\rangle\\
    &\le-\lambda\|u_i\|^2+C_Hh\,L_\mu(1+C_Hh)\|u\|^2
     \\&\le(-\lambda+Ch)\|u\|^2,
\end{align*}
where we used $\|u_i\|^2\ge(1-2C_Hh)\|u\|^2$, which follows from $\|u_i\|\ge(1-C_Hh)\|u\|$ for $C_Hh\le1$ and holds trivially otherwise.
Inserting these calculations and using $\sum_i\alpha_i=1$, we get
\[
    \|\Psi_k(x)-\Psi_k(y)\|^2\le\|u\|^2+2h\sum_i\alpha_i(-\lambda+Ch)\|u\|^2
      +h^2\Big(\sum_i\alpha_iL_\mu(1+C_Hh)\Big)^2\|u\|^2
    \le(1-2\lambda h+Ch^2)\|u\|^2,
\]
which implies \ref{H1} with $\kappa=\lambda/2$ for $h\le h_0^\text{SRA}:=\min\{1,1/\lambda,\lambda/(2C)\}$. Without Assumption~\ref{ass:Dissipativity} the same computation with $-\lambda$ replaced by
$L_\mu$ gives $\|\Psi_k(x)-\Psi_k(y)\|^2\le(1+2L_\mu h+Ch^2)\|u\|^2$, and hence \ref{H1} with $\kappa=-C_\Psi$, $C_\Psi:=L_\mu+\frac C2$, as in (a). (Note that $\alpha_i\ge0$ is
needed only to pass the dissipativity through the convex combination. It specifically holds for SRA1 and for Chang's method.)

For \ref{H2}, split $\Delta_k=\Delta_k^{\mathrm{drift}}+\Delta_k^{\mathrm{noise}}$ with
\begin{align*}
    \Delta_k^{\mathrm{drift}}
    &:=\int_{t_k}^{t_{k+1}}\mu(s,\Forward_s)\,\d s-h\sum_i\alpha_i\mu_i\big(H_i(\Forward_{t_k})\big),\\
    \Delta_k^{\mathrm{noise}}
    &:=\int_{t_k}^{t_{k+1}}g(s)\,\d W_s
      -\sum_i\Big(\beta^{(1)}_i\Delta W_{t_k}+\tfrac{\beta^{(2)}_i}{h}\Delta Z_{t_k}\Big)g(t_k+c^{(1)}_ih).
\end{align*}
For the drift part, Minkowski's integral inequality and linear growth give
$\|\int_{t_k}^{t_{k+1}}\mu(s,\Forward_s)\,\d s\|_{L_{2n}}\le L_\mu h(1+\sup_t\|\Forward_t\|_{L_{2n}})$, while
$\|H_i(\Forward_{t_k})\|_{L_{2n}}\le C(1+\|\Forward_{t_k}\|_{L_{2n}})$ by induction over the
stages (using $\|\Delta Z_{t_k}\|_{L_{2n}}=\O(h^{3/2})$, see below, and $h\le1$), so that
$\|\Delta_k^{\mathrm{drift}}\|_{L_{2n}}\le Ch(1+\sup_t\|\Forward_t\|_{L_{2n}})$.
For the noise part, using ${\beta^{(1)}}^Te=1$ and ${\beta^{(2)}}^Te=0$,
\[
    \Delta_k^{\mathrm{noise}}
    =\int_{t_k}^{t_{k+1}}\big(g(s)-g(t_k)\big)\d W_s
     -\sum_i\beta^{(1)}_i\big(g(t_k+c^{(1)}_ih)-g(t_k)\big)\Delta W_{t_k}
     -\sum_i\frac{\beta^{(2)}_i}{h}\big(g(t_k+c^{(1)}_ih)-g(t_k)\big)\Delta Z_{t_k}.
\]
Each difference of $g$ is bounded by $L_gh$, since $c^{(1)}\in[0,1]^s$. The first term is a Wiener integral with deterministic integrand, hence Gaussian with covariance $\sigma_k^2I_d$, where $\sigma_k^2=\int_{t_k}^{t_{k+1}}(g(s)-g(t_k))^2\,\d s\le\frac{L_g^2h^3}{3}$. For $\xi\sim\mathcal N(0,I_d)$ let $\gamma_{d,n}:=\|\xi\|_{L_{2n}}$. Then the first term has $L_{2n}$-norm $\sigma_k\gamma_{d,n}\le\frac{L_g\gamma_{d,n}}{\sqrt3}h^{3/2}$, and in the same way \eqref{eq:VarofDeltaZ} gives $\|\Delta W_{t_k}\|_{L_{2n}}=\gamma_{d,n}h^{1/2}$ and $\|\Delta Z_{t_k}\|_{L_{2n}}=\frac{\gamma_{d,n}}{\sqrt3}h^{3/2}$. Together, $\|\Delta_k^{\mathrm{noise}}\|_{L_{2n}}\le Ch^{3/2}\le Ch$, and \ref{H2} holds with $C_\delta=C(1+\sup_t\|\Forward_t\|_{L_{2n}})$, which is finite by Lemma~\ref{lem:MomentEstimateOfStrongSolution}, and bounded polynomial in $T$ by Lemma~\ref{lem:MomentEstimateDissipative} under Assumption~\ref{ass:Dissipativity}.
\end{proof}

\subsection{Main Order-Proofs}
\label{sec:MainOrderProofs}
\begin{proof}[\textbf{Proof for the Euler--Maruyama scheme (Theorem \ref{lem:EMScheme})}]
    We need to verify that the scheme includes all terms in $\mathcal A_1\setminus\{\nu\}$, which in our setup is simply $\{(0)\}\cup[d]$. Thus, for an Itô--Taylor scheme to have order $\gamma=1$, it needs to include the terms 
    \[f_{(0)}=\mu(t,x) \quad \text{and} \quad f_{(j)}=g(t)e_j,\]
    which is precisely \eqref{eq:EMDiscretization}.
\end{proof}

\begin{proof}[\textbf{Proof for the exponential integrator scheme (Theorem \ref{lem:EIScheme})}]
    Recall that in our situation we have $\mathcal A_1\setminus\{\nu\}=\{(0)\}\cup[d]\cup[d]^2$. Since $f_{(j_1,j_2)}\equiv0$ for all $j_1,j_2\in[d]$ by spatial independence of $g$, we set $g_{\alpha,k}:=0$ for these multi-indices, which trivially satisfies the tolerance in Definition~\ref{def:generalstrongscheme}. It therefore suffices to identify $g_{(0),k}$, $g_{(j),k}$ ($j=1,\ldots,d$), and a remainder $R_k=R_k^M+R_k^D$. We define all $g_{\alpha,k}$ as exact matches of the coefficient functions $f_\alpha(t_k,\BackwardNormalstartingNetworkDiscretized_k)$ and move every leftover into the remainder, so that only the remainder has to be bounded. Throughout, we use Assumption~\ref{ass:EI}, write $F:=\sup_{\tau\in[0,T]}|f(\tau)|$, and denote by $L_f$ and $L_g$ the Lipschitz constants of $f$ and $g$, which, like $M_g$ and $C_s$, grow at most polynomially in $T$. Moreover, we let $h\le h_0^\text{EI}$ with $h_0^\text{EI}$ from Lemma~\ref{lem:MomentsEISRA}\ref{lem:MomentsEISRA_EI}, so that in particular $h\le1$ and $Fh\le1$. For $s\in[t_k,t_{k+1}]$, the exponent $\int_s^{t_{k+1}}f(T-u)\,\d u$ of $\Phi_{t_{k+1},s}$ lies in $[-Fh,Fh]$. Using $e^x-1\le xe^x$ for $x\ge0$, where $e$ denotes Euler's number, we therefore have
    \begin{equation}
        \label{eq:EIIntegrandBounds}
        \Phi_{t_{k+1},s}\le e^{Fh}\le e,\qquad |\Phi_{t_{k+1},s}-1|\le e^{Fh}-1\le eFh,\qquad |g(T-s)-g(\tau_k)|\le L_gh.
    \end{equation}
    We start with the stochastic part, which we split as
    \[
        b_k=\int_{t_k}^{t_{k+1}}\Phi_{t_{k+1},s}\,g(T-s)\,\d W_s=\sum_{j=1}^d\underbrace{g(\tau_k)e_j}_{g_{(j),k}:=}\Delta W_{t_k}^{(j)}+\rho_{\text{Stochastic}},
        \qquad
        \rho_{\text{Stochastic}}:=\int_{t_k}^{t_{k+1}}\big(\Phi_{t_{k+1},s}\,g(T-s)-g(\tau_k)\big)\,\d W_s.
    \]
    Writing $\Phi_{t_{k+1},s}\,g(T-s)-g(\tau_k)=(\Phi_{t_{k+1},s}-1)\,g(T-s)+\big(g(T-s)-g(\tau_k)\big)$, \eqref{eq:EIIntegrandBounds} and $|g|\le M_g$ show that the deterministic integrand of $\rho_{\text{Stochastic}}$ satisfies $|\Phi_{t_{k+1},s}\,g(T-s)-g(\tau_k)|\le(eM_gF+L_g)h$. We continue with the drift part. With $\varphi_k:=\int_{t_k}^{t_{k+1}}f(T-u)\,\d u$, we have $\Phi_k=e^{\varphi_k}$, $|\varphi_k|\le Fh\le1$ and $|\varphi_k-f(\tau_k)h|\le\frac{L_f}{2}h^2$. Using $|e^{\varphi}-1-\varphi|\le\frac12\varphi^2e^{|\varphi|}\le\frac e2\varphi^2$ for $|\varphi|\le1$, we obtain
    \[
        |\Phi_k-1-f(\tau_k)h|\le|e^{\varphi_k}-1-\varphi_k|+|\varphi_k-f(\tau_k)h|\le\tfrac12\big(eF^2+L_f\big)h^2.
    \]
    For $a_k$, we write
    \[
        a_k-g^2(\tau_k)h=\int_{t_k}^{t_{k+1}}\Big[(\Phi_{t_{k+1},s}-1)\,g^2(T-s)+\big(g^2(T-s)-g^2(\tau_k)\big)\Big]\d s,
    \]
    and use \eqref{eq:EIIntegrandBounds} together with $|g^2(T-s)-g^2(\tau_k)|\le|g(T-s)-g(\tau_k)|\,|g(T-s)+g(\tau_k)|\le2M_gL_gh$. Altogether,
    \[
        \Phi_k = 1+f(\tau_k)h+r_\Phi, \qquad a_k = g^2(\tau_k)h+r_a,\qquad |r_\Phi|\le c_\Phi h^2,\quad |r_a|\le c_ah^2,
    \]
    with $c_\Phi:=\frac12\big(eF^2+L_f\big)$ and $c_a:=eM_g^2F+2M_gL_g$, so that
    \begin{equation}
        \label{eq:EIMasterExpansion}
        \begin{array}{rl}
             \Phi_k\BackwardNormalstartingNetworkDiscretized_k + a_k\,s_\theta(\tau_k,\BackwardNormalstartingNetworkDiscretized_k)
        &= \BackwardNormalstartingNetworkDiscretized_k + \underbrace{\big[f(\tau_k)\BackwardNormalstartingNetworkDiscretized_k+g^2(\tau_k)s_\theta(\tau_k,\BackwardNormalstartingNetworkDiscretized_k)\big]}_{g_{(0),k}:=}h + \rho_{\text{Drift}} \\[6ex]
        \rho_{\text{Drift}}&:=r_\Phi\BackwardNormalstartingNetworkDiscretized_k+r_a\,s_\theta(\tau_k,\BackwardNormalstartingNetworkDiscretized_k).
        \end{array}
    \end{equation}
    Collecting all the definitions of $g_{\alpha,k}$ made throughout this proof, we can write the scheme as
    \[
        \BackwardNormalstartingNetworkDiscretized_{k+1} = \BackwardNormalstartingNetworkDiscretized_k + I_{(0)}\big[g_{(0),k}\big]_{t_k,t_{k+1}} + \sum_{j=1}^dI_{(j)}\big[g_{(j),k}\big]_{t_k,t_{k+1}} + \underbrace{\rho_{\text{Stochastic}}}_{R_k^M:=} + \underbrace{\rho_{\text{Drift}}}_{R_k^D:=}.
    \]
    From construction, $g_{(0),k}=\mu(t_k,\BackwardNormalstartingNetworkDiscretized_k)=f_{(0)}(t_k,\BackwardNormalstartingNetworkDiscretized_k)$ and $g_{(j),k}=g(\tau_k)e_j=f_{(j)}(t_k,\BackwardNormalstartingNetworkDiscretized_k)$ match the coefficient functions exactly (cf.\ \eqref{eq:coefficientfunctions1} with diffusion coefficient $g(T-t)$). Thus, to prove the claim, we only need to verify the conditions on the remainders $R^M_k$ and $R^D_k$.
    \begin{description}
        \item[\textbf{Bounding $R^M_k$}:] $R^M_k$ is a Wiener integral over $[t_k,t_{k+1}]$ with deterministic integrand. Hence it is independent of $\mathcal F_{t_k}$ and centered, so that $\mathbb{E}[R^M_k|\mathcal{F}_{t_k}]=0$. Furthermore, the Itô isometry and the bound on its integrand yield
        \[\mathbb{E}\left[\|R^M_k\|^2\right]=d\int_{t_k}^{t_{k+1}}\big(\Phi_{t_{k+1},s}\,g(T-s)-g(\tau_k)\big)^2\,\d s\leq d\,(eM_gF+L_g)^2h^3=Ch^{2\gamma+1}|_{\gamma = 1},\]
        which verifies the necessary conditions.
        \item[\textbf{Bounding $R^D_k$}:] By Assumption~\ref{ass:EI}, $\|s_\theta(\tau_k,\BackwardNormalstartingNetworkDiscretized_k)\|\le C_s(1+\|\BackwardNormalstartingNetworkDiscretized_k\|)$, and thus $\|R^D_k\|\le(c_\Phi+c_aC_s)h^2(1+\|\BackwardNormalstartingNetworkDiscretized_k\|)$. With the moment bound $\sup_k\mathbb E[\|\BackwardNormalstartingNetworkDiscretized_k\|^2]\le M_1'$ from Lemma~\ref{lem:MomentsEISRA}, valid for $h\le h_0^\text{EI}$, we obtain
        \[\max_{0\le k\le K-1}\mathbb E\left[\|R^D_k\|^2\right]\le 2(c_\Phi+c_aC_s)^2(1+M_1')\,h^4=Ch^{2\gamma+2}|_{\gamma = 1},\]
        which verifies the necessary conditions.
    \end{description}
    Both conditions of Definition~\ref{def:generalstrongscheme} hold with $\gamma=1$, which proves the claim. All constants depend only on $F$, $L_f$, $L_g$, $M_g$, $C_s$, $d$ and $M_1'$. All of these except $M_1'$ grow at most polynomially in $T$ by Assumptions~\ref{ass:CoefficientFunctions} and~\ref{ass:EI}. Under Assumption~\ref{ass:Dissipativity}, the same holds for $M_1'$ by Lemma~\ref{lem:MomentsEISRA}, so that $C$ and $C_R$ grow at most polynomially in $T$. In general, $M_1'$ is bounded by $e^{CT}$ with $C$ growing at most polynomially in $T$, which is absorbed into $e^{\text{Poly}(T)}$ in Theorem~\ref{thm:Main}.
\end{proof}

\begin{proof}[\textbf{Proof for the explicit SRA Family (Theorem \ref{lem:SRA})}]
    Recall
    \[\mathcal A_{3/2}\setminus\{\nu\}=\{(0)\}\cup[d]\cup\{(j,0):j\in[d]\}\cup\{(0,j):j\in[d]\}\cup\{(0,0)\}\cup[d]^2\cup[d]^3.\]
    Since $f_{(j_1,j_2)}\equiv0$ for all $j_1,j_2\in[d]$ and $f_{(j_1,j_2,j_3)}\equiv0$ for all $j_1,j_2,j_3\in[d]$ by spatial independence of $g$, we set $g_{\alpha,k}:=0$ for these multi-indices, which trivially satisfies the tolerance in Definition~\ref{def:generalstrongscheme}. It therefore suffices to identify $g_{(0),k}$, $g_{(j),k}$, $g_{(j,0),k}$, $g_{(0,j),k}$ ($j=1,\ldots,d$), $g_{(0,0),k}$, and a remainder $R_k=R_k^M+R_k^D$. Throughout, we abbreviate $\mu_k:=\mu(t_k,x_k)$, $g_k:=g(t_k)$, $H_i:=H_i^{(0)}$, and write $\partial_t\mu_k$, $\partial_t^2\mu_k$, $D\mu_k$, $\partial_tD\mu_k$ and $D^2\mu_k$ for the derivatives of $\mu$ (cf. \eqref{eq:DifferentialOperatorDefinition}) evaluated at $(t_k,x_k)$. Furthermore, we use Assumption~\ref{ass:SRA} and write $P_k:=1+\|x_k\|$, which is $\mathcal F_{t_k}$-measurable. By Assumption~\ref{ass:SRA}, $\mu_k$ and all derivatives of $\mu$ up to total order three, evaluated at $(t_k,x_k)$, are bounded in norm by $CP_k$, while $D\mu$ is bounded by the Lipschitz constant of $\mu$. By Lemma~\ref{lem:MomentsEISRA}, the moments $M_n':=\sup_k\mathbb E[\|x_k\|^{2n}]$ are finite for every $n$, so that $\mathbb E[P_k^{2n}]\le2^{2n-1}(1+M_n')$. The generic constant $C$ below may depend on $M_n'$ for $n\le4$, but not on $h$ or $k$. We start with the stochastic part: Using Taylor, we get
    \[g(t_k+c_i^{(1)}h)=g(t_k)+c_{i}^{(1)}hg'(t_k) + \O(h^2),\]
    and thus
    \begin{align*}
        \sum_{i=1}^s\left(\beta_i^{(1)}\Delta W_{t_k}+\frac{\beta_i^{(2)}}{h}\Delta Z_{t_k}\right)g(t_k+c_i^{(1)}h) &= \sum_{i=1}^s\left(\beta_i^{(1)}\Delta W_{t_k}+\frac{\beta_i^{(2)}}{h}\Delta Z_{t_k}\right)(g(t_k)+c_{i}^{(1)}hg'(t_k) + \O(h^2))\\
        &=\sum_{i=1}^s\beta_i^{(1)}g(t_k)\Delta W_{t_k} + \sum_{i=1}^s \frac{\beta_i^{(2)}}{h}g(t_k)\Delta Z_{t_k}\\&\qquad + \sum_{i=1}^s\left(c_i^{(1)}\beta_i^{(1)}g'(t_k) h \Delta W_{t_k} + c_i^{(1)}\beta_i^{(2)}g'(t_k)\Delta Z_{t_k}\right) \\&\qquad + \sum_{i=1}^s \left(\beta_i^{(1)}\Delta W_{t_k}\O(h^2)+\beta_i^{(2)}\Delta Z_{t_k}\O(h)\right)
    \end{align*}
    Inserting the order conditions ${\beta^{(1)}}^Te=1$,${\beta^{(2)}}^Te=0$, ${\beta^{(1)}}^Tc^{(1)}=1$ and ${\beta^{(2)}}^Tc^{(1)}=-1$ yields
    \begin{align*}
        \sum_{i=1}^s\left(\beta_i^{(1)}\Delta W_{t_k}+\frac{\beta_i^{(2)}}{h}\Delta Z_{t_k}\right)g(t_k+c_i^{(1)}h) = \underbrace{g(t_k)}_{g_{(j),k}:=}\Delta W_{t_k} + \underbrace{g'(t_k)}_{g_{(0,j),k}:=}(h\Delta W_{t_k}-\Delta Z_{t_k}) + \rho_{\text{Stochastic}},
    \end{align*}
    where $\rho_{\text{Stochastic}} = \sum_{i=1}^s \left(\beta_i^{(1)}\Delta W_{t_k}\O(h^2)+\beta_i^{(2)}\Delta Z_{t_k}\O(h)\right)$. We continue with expanding the drift:
    \begin{align*}
        \mu(t_k+c_i^{(0)}h,H_i)&=\mu_k+c_i^{(0)}h\partial_t\mu_k+D\mu_k(H_i-x_k) + \underbrace{\frac{1}{2}\left(\partial_t^2\mu_k(c_i^{(0)}h)^2 + 2c_i^{(0)}h\,\partial_tD\mu_k(H_i-x_k)\right)}_{\text{rem}_1:=}\\
        &\qquad + \frac12 D^2\mu_k[H_i-x_k,H_i-x_k] + \O\big((P_k+\|H_i-x_k\|)(h+\|H_i-x_k\|)^3\big),
    \end{align*}
    where the factor $P_k+\|H_i-x_k\|$ in the Taylor remainder stems from the linear growth of the third derivatives of $\mu$, evaluated at an intermediate point between $(t_k,x_k)$ and $(t_k+c_i^{(0)}h,H_i)$. Inserting
    \begin{align*}
        H_i-x_k=\sum_{j=1}^sA_{ij}^{(0)}h\mu_k +\O(h^2P_k+h\max_j\|H_j-x_k\|)+\sum_{j=1}^sB_{ij}^{(0)}g_k\frac{1}{h}\Delta Z_{t_k}  + \O(\|\Delta Z_{t_k}\|)
    \end{align*}
    where we write $a_i$ and $b_i$ for the two sums and $r_i$ for the two $\O$ terms, yields
    \begin{align*}
        \mu(t_k+c_i^{(0)}h,H_i)&=\mu_k+c_i^{(0)}h\partial_t\mu_k+\sum_{j=1}^sA_{ij}^{(0)}h\,D\mu_k\mu_k + \sum_{j=1}^sB_{ij}^{(0)}g_k\frac{1}{h}D\mu_k\Delta Z_{t_k} \\
        &\qquad + \underbrace{\frac{1}{2}\left(\sum_{j_1,j_2=1}^sA_{ij_1}^{(0)}A_{ij_2}^{(0)}h^2D^2\mu_k[\mu_k,\mu_k] + 2\sum_{j_1,j_2=1}^sA_{ij_1}^{(0)}B_{ij_2}^{(0)}g_kD^2\mu_k[\mu_k,\Delta Z_{t_k}]\right)}_{\text{rem}_2:=} \\
        &\qquad + \frac{1}{2}\sum_{j_1,j_2=1}^sB_{ij_1}^{(0)}B_{ij_2}^{(0)}g_k^2\frac{1}{h^2}D^2\mu_k[\Delta Z_{t_k},\Delta Z_{t_k}]\qquad +\text{rem}_1+\text{rem}_3,
    \end{align*}
    where
    \begin{equation}
        \label{eq:rem_3}
        \begin{array}{rl}
             \text{rem}_3&:=\O(h^2P_k+h\max_j\|H_j-x_k\|)+\O(\|\Delta Z_{t_k}\|)+\O\big(P_k(\|a_i\|+\|b_i\|+\|r_i\|)\|r_i\|\big)\\[1ex]
             &\qquad+\O\big((P_k+\|H_i-x_k\|)(h+\|H_i-x_k\|)^3\big).
        \end{array}
    \end{equation}
    Thus, expanding the whole term we get
    \begin{align*}
        \sum_{i=1}^s\alpha_i\mu(t_k+c_i^{(0)}h, H_i)h&=\sum_{i=1}^s\alpha_ih\mu_k + \sum_{i=1}^s\alpha_ic_i^{(0)}h^2\partial_t\mu_k + \sum_{i=1}^s\alpha_i h^2\sum_{j=1}^sA_{ij}^{(0)}D\mu_k\mu_k\\
        &\qquad + \sum_{i=1}^s \alpha_i h \sum_{j=1}^sB_{ij}^{(0)}g_k\frac{1}{h}D\mu_k\Delta Z_{t_k}\\
        &\qquad+\sum_{i=1}^s\alpha_ih\frac{1}{2}\sum_{j_1,j_2=1}^sB_{ij_1}^{(0)}B_{ij_2}^{(0)}g_k^2\frac{1}{h^2}D^2\mu_k[\Delta Z_{t_k},\Delta Z_{t_k}] +\rho_{\text{Drift}},
    \end{align*}
    where $\rho_{\text{Drift}}:=h\sum_{i=1}^s\alpha_i\left(\mathrm{rem}_1+\mathrm{rem}_2+\mathrm{rem}_3\right)$ collects the remainders of all stages.
    Inserting the order conditions, $\alpha^Te=1$, $\alpha^TB^{(0)}e=1$, $\alpha ^TA^{(0)}e=\frac{1}{2}$, and $\alpha^T(B^{(0)}e)^2=\frac32$, as well as $c^{(0)}=A^{(0)}e$ we get
    \begin{align*}
        \sum_{i=1}^s\alpha_i\mu(t_k+c_i^{(0)}h, H_i)h&=h\underbrace{\mu_k}_{g_{(0),k}:=} + \frac12h^2\partial_t\mu_k + \frac{1}{2} h^2D\mu_k\mu_k  + \sum_{j=1}^d\underbrace{g_kD\mu_ke_j}_{g_{(j,0),k}:=} \Delta Z_{t_k}^{(j)}\\
        &\qquad+\underbrace{\frac{3g_k^2}{4h}D^2\mu_k[\Delta Z_{t_k},\Delta Z_{t_k}]}_{T_k:=}+\rho_{\text{Drift}}.
    \end{align*}
    We proceed by calculating the mean of $T_k$. In the representation of $\Delta Z_{t_k}$ in \eqref{eq:DeltaZDecomposition}, we can see that the components are independent of one another, and thus, together with the variance computation of $\Delta Z_{t_k}$ in \eqref{eq:VarofDeltaZ} and the bilinearity of $D^2\mu_k$, we get
    \begin{align*}
        \mathbb E\left[T_k\mid \mathcal{F}_{t_k} \right]&=\frac{3g_k^2}{4h}\sum_{n_1,n_2=1}^dD^2\mu_k[e_{n_1},e_{n_2}] \mathbb E[\Delta Z_{t_k}^{(n_1)}\Delta Z_{t_k}^{(n_2)}\mid \mathcal{F}_{t_k}]&&=\frac{3g_k^2}{4h}\sum_{n=1}^dD^2\mu_k[e_n,e_n] \mathbb E\left[\left(\Delta Z_{t_k}^{(n)}\right)^2\mid \mathcal{F}_{t_k}\right]\\
        &=\frac{3g_k^2}{4h}\sum_{n=1}^dD^2\mu_k[e_n,e_n] \frac{h^3}{3}&& =\frac{1}{4}g_k^2h^2\sum_{n=1}^dD^2\mu_k[e_n,e_n].
    \end{align*}
    We define \[g_{(0,0),k}:=\partial_t\mu_k + D\mu_k\mu_k + \frac{2}{h^2} \mathbb E\left[T_k\mid \mathcal{F}_{t_k} \right] = \partial_t\mu_k + D\mu_k\mu_k + \frac{1}{2}g_k^2\sum_{n=1}^dD^2\mu_k[e_n,e_n]. \]
    Collecting all the definitions of $g_{\alpha,k}$ made throughout this proof, we can write the scheme as
    \begin{align*} x_{k+1}&=x_k + g_{(0),k}h + \sum_{j=1}^dg_{(j),k}\Delta W_{t_k}^{(j)} + \sum_{j=1}^dg_{(j,0),k}\Delta Z_{t_k}^{(j)} + \sum_{j=1}^dg_{(0,j),k}(h\Delta W_{t_k}-\Delta Z_{t_k})^{(j)} + g_{(0,0),k}\frac{h^2}{2}\\
    &\qquad + \underbrace{(T_k-\mathbb E[T_k|\mathcal F_{t_k}])}_{R^M_k:=} + \underbrace{\rho_{\text{Stochastic}}+\rho_{\text{Drift}}}_{R^D_k:=}.
    \end{align*}
    From construction, we can see that all the $g_{(\alpha),k}$ match the corresponding $f_{\alpha}(t_k,x_k)$ exactly (cf. \eqref{eq:coefficientfunctions1} and \eqref{eq:coefficientfunctions2}). Thus, to prove the claim, we only need to verify the conditions on the remainders $R^D_k$ and $R^M_k$. Throughout, we use that $\Delta Z_{t_k}$ is a centered Gaussian vector with independent components of variance $\frac{h^3}{3}$ by \eqref{eq:DeltaZDecomposition} and \eqref{eq:VarofDeltaZ}, independent of $\mathcal F_{t_k}$, so that $\mathbb E[\|\Delta Z_{t_k}\|^{2p}]\leq Ch^{3p}$ for $p\le4$, and $\mathbb E[P_k^{2q}\|\Delta Z_{t_k}\|^{2p}]=\mathbb E[P_k^{2q}]\,\mathbb E[\|\Delta Z_{t_k}\|^{2p}]$.
    \begin{description}
        \item[\textbf{Bounding $R^M_k$}:] By construction, we have $\mathbb{E}[R^M_k|\mathcal{F}_{t_k}]=0$. Furthermore, since subtracting the conditional mean does not increase the second moment, and $\|T_k\|\le\frac{C}{h}P_k\|\Delta Z_{t_k}\|^2$, we can infer
        \[\mathbb{E}\left[\|R^M_k\|^2\right]\leq \mathbb E\left[\|T_k\|^2\right]\le\frac{C}{h^2}\mathbb E[P_k^2]\,\mathbb{E}\left[\|\Delta Z_{t_k}\|^4\right]\le\frac{C}{h^2}h^6=Ch^{2\gamma+1}|_{\gamma = 1.5},\]
        which verifies the necessary conditions.
        \item[\textbf{Bounding $R^D_k$}:] $R^D_k$ consists of the stochastic remainder $\rho_{\text{Stochastic}}$ and the drift remainder $\rho_{\text{Drift}}$, which collects the remainders $\mathrm{rem}_{1}$, $\mathrm{rem}_{2}$ and $\mathrm{rem}_{3}$ made in the expansion of the drift. For the stochastic remainder we get
        \[\mathbb E[\|\rho_{\text{Stochastic}}\|^2]\leq Ch^4\mathbb E[\|\Delta W_{t_k}\|^2]+Ch^2\mathbb E[\|\Delta Z_{t_k}\|^2]=\O(h^5).\]
        For the drift remainder $\rho_{\text{Drift}}=h\sum_{i=1}^s\alpha_i(\mathrm{rem}_1+\mathrm{rem}_2+\mathrm{rem}_3)$, it suffices to show $\mathbb E[\|\mathrm{rem}_l\|^2]\leq Ch^3$ for $l=1,2,3$ and every stage $i$. We begin by bounding the moments of $H_i-x_k$, as they will be used throughout. From \eqref{eq:SRA_b}, the linear growth of $\mu$ and an induction over the explicit stages, we get $\|H_i-x_k\|\leq ChP_k+\frac{C}{h}\|\Delta Z_{t_k}\|$, and thus, for $h\leq1$ and all $p,q\in\N_0$ with $p+q\le4$,
        \begin{equation}
            \label{eq:HMoments}
            \mathbb E\big[P_k^{2q}\|H_i-x_k\|^{2p}\big]\leq C\left(h^{2p}\,\mathbb E[P_k^{2(p+q)}]+h^{-2p}\,\mathbb E[P_k^{2q}]\,\mathbb E[\|\Delta Z_{t_k}\|^{2p}]\right)\leq Ch^p.
        \end{equation}
        Since $\|\mathrm{rem}_1\|\leq CP_k(h^2+h\|H_i-x_k\|)$ and $\|\mathrm{rem}_2\|\leq Ch^2P_k^3+CP_k^2\|\Delta Z_{t_k}\|$, we get
        \begin{align*}
            \mathbb E[\|\mathrm{rem}_1\|^2]\leq Ch^4\mathbb E[P_k^2]+Ch^2\mathbb E[P_k^2\|H_i-x_k\|^2]\leq Ch^3,\\
            \mathbb E[\|\mathrm{rem}_2\|^2]\leq Ch^4\mathbb E[P_k^6]+C\mathbb E[P_k^4]\,\mathbb E[\|\Delta Z_{t_k}\|^2]\leq Ch^3.
        \end{align*}
        For $\mathrm{rem}_3$ (cf. \eqref{eq:rem_3}), we bound its four terms separately. For the first two, $\mathbb E[(h^2P_k+h\max_j\|H_j-x_k\|)^2]\leq Ch^4\mathbb E[P_k^2]+Ch^2\sum_{j=1}^s\mathbb E[\|H_j-x_k\|^2]\leq Ch^3$ and $\mathbb E[\|\Delta Z_{t_k}\|^2]\leq Ch^3$. For the third, we first notice that $\|a_i\|\leq ChP_k$ and $\|b_i\|\le\frac Ch\|\Delta Z_{t_k}\|$, and, inserting the bound on $\|H_j-x_k\|$, $\|r_i\|\leq C(h^2P_k+\|\Delta Z_{t_k}\|)$. Hence, $\mathbb E[P_k^4(\|a_i\|+\|b_i\|+\|r_i\|)^4]\leq C(h^4\mathbb E[P_k^8]+h^{-4}\mathbb E[P_k^4]\,\mathbb E[\|\Delta Z_{t_k}\|^4])\leq Ch^2$ and $\mathbb E[\|r_i\|^4]\leq C(h^8\mathbb E[P_k^4]+\mathbb E[\|\Delta Z_{t_k}\|^4])\leq Ch^6$. Then, the Cauchy--Schwarz inequality yields
        \[\mathbb E\left[P_k^2(\|a_i\|+\|b_i\|+\|r_i\|)^2\|r_i\|^2\right]\leq \mathbb E\left[P_k^4(\|a_i\|+\|b_i\|+\|r_i\|)^4\right]^{\frac12}\mathbb E\left[\|r_i\|^4\right]^{\frac12}\leq Ch\cdot Ch^3=Ch^4,\]
        and for the Taylor remainder, \eqref{eq:HMoments} yields
        \begin{align*}
            &\mathbb E\big[(P_k+\|H_i-x_k\|)^2(h+\|H_i-x_k\|)^6\big]\\
            \leq\;& C\Big(h^6\mathbb E[P_k^2]+\mathbb E[P_k^2\|H_i-x_k\|^6]+h^6\mathbb E[\|H_i-x_k\|^2]+\mathbb E[\|H_i-x_k\|^8]\Big)\\
            \leq\; & Ch^3.
        \end{align*}
        Hence $\mathbb E[\|\mathrm{rem}_3\|^2]\leq Ch^3$, and altogether
        \[\mathbb E[\|\rho_{\text{Drift}}\|^2]\leq Ch^2\sum_{i=1}^s\sum_{l=1}^3\mathbb E[\|\mathrm{rem}_l\|^2]\leq Ch^5.\]
        Combining both remainders, we obtain
        \[\mathbb E[\|R^D_k\|^2]\leq 2\,\mathbb E[\|\rho_{\text{Stochastic}}\|^2]+2\,\mathbb E[\|\rho_{\text{Drift}}\|^2]\leq Ch^5=Ch^{2\gamma+2}|_{\gamma = 1.5},\]
        which verifies the necessary conditions and concludes the proof.
    \end{description}
    All constants above depend on $T$ only through the moments $M_n'$, $n\le4$, and the constants of Assumptions~\ref{ass:CoefficientFunctions} and~\ref{ass:SRA}, which grow at most polynomially in $T$. Under Assumption~\ref{ass:Dissipativity}, the same holds for $M_n'$ by Lemma~\ref{lem:MomentsEISRA}. In general, $M_n'\le e^{\text{Poly}(T)}$, so that $C,C_R=\O(e^{\text{Poly}(T)})$.
\end{proof}

\bibliographystyle{abbrv}
\bibliography{references}
\end{document}